\documentclass[11pt]{article}
\usepackage{amssymb,latexsym,epsf,xypic,amsmath}

\usepackage[mathscr]{euscript}

\DeclareFontFamily{OT1}{rsfs}{}
\DeclareFontShape{OT1}{rsfs}{n}{it}{<-> rsfs10}{}
\DeclareMathAlphabet{\curly}{OT1}{rsfs}{n}{it}

\newcommand{\rt}[1]{\stackrel{#1\,}{\rightarrow}}
\newcommand{\Rt}[1]{\stackrel{#1\,}{\longrightarrow}}
\newcommand\To{\longrightarrow}
\newcommand\ot{\leftarrow}
\newcommand\Ot{\longleftarrow}
\newcommand\into{\hookrightarrow}
\newcommand\bull{{\scriptscriptstyle\bullet}}
\newcommand\udot{^\bull}
\newcommand\R{\mathbf R\hspace{1pt}}
\newcommand\C{\mathbb C\,}
\newcommand\Z{\mathbb Z}
\newcommand\Pee{\mathbb P}
\newcommand\Ft{F^{0,2}_A}
\newcommand\iL{i\Lambda F^{1,1}_A}
\newcommand\db{\bar{\partial}}
\newcommand\A{\curly A}
\newcommand\B{\curly B}

\newcommand\M{\mathcal M}
\newcommand\N{\mathcal N}
\newcommand\E{\mathcal E}
\newcommand\F{\mathcal F}
\newcommand\Ll{\mathcal L} 
\newcommand\KK{\mathscr K}
\newcommand\K{\curly K}
\newcommand\m{\mathfrak m}
\newcommand\n{\mathfrak n}
\newcommand\OO{\mathscr O}
\newcommand\I{\curly I}
\newcommand\II{\mathcal I}
\newcommand\T{\mathcal T}
\newcommand\V{\mathcal V}
\newcommand\X{\mathcal X}

\newcommand\AND{\hspace{-1.3mm}&\hspace{-1.3mm}}
 
\newcommand\xp{\!\times_{\Pee^1}\!}
\newcommand\x{\times}
\newcommand\res{\arrowvert_}

\DeclareMathOperator\tr{tr}
\DeclareMathOperator\id{id}
\DeclareMathOperator\ad{ad}
\DeclareMathOperator\Td{Td}

\DeclareMathOperator\Spec{Spec}
\newcommand\Endo{\mathrm{End}_0\,}
\newcommand\Ext{\mathcal E\hspace{-1pt}xt}
\newcommand\Hom{\mathcal H\hspace{-1pt}om}
\newcommand\Homd{\mathcal H\hspace{-1pt}om^\bull}
\newcommand\take{\,\backslash\,}
\newcommand\otspam{\ot\!\shortmid}
\newcommand\dd{\widehat{\mathcal E\,}}
\newcommand\ve{^\vee}

\makeatletter \@addtoreset{equation}{section} \makeatother

\newtheorem{Lemma}[equation]{Lemma}
\newtheorem{Theorem}[equation]{Theorem}
\newtheorem{Prop}[equation]{Proposition}
\newtheorem{Definition}[equation]{Definition}
\newtheorem{Corollary}[equation]{Corollary}
\newenvironment{Proof}{\noindent\emph{Proof\ }}{\ \ q.e.d.\\}
\newenvironment{Remarks}{\noindent\textbf{Remarks}\ }{\\}
\newenvironment{Remark}{\noindent\textbf{Remark}\ }{\\}

\begin{document}

\title{\textbf{A holomorphic Casson
invariant for Calabi-Yau 3-folds, and bundles on $K3$ fibrations}}
\author{R.\,P. Thomas}
\date{}
\maketitle

\begin{abstract} \noindent
We briefly review the formal picture in which a Calabi-Yau $n$-fold is
the complex analogue of an oriented real $n$-manifold, and a Fano with
a fixed smooth anticanonical divisor is the analogue of a manifold
with boundary, motivating a holomorphic Casson invariant counting
bundles on a Calabi-Yau 3-fold. We develop the deformation theory
necessary to obtain the virtual moduli cycles of
\cite{LT}, \cite{BF} in moduli spaces of stable sheaves whose higher
obstruction groups vanish. This gives, for instance, virtual moduli
cycles in Hilbert schemes of curves in $\Pee^3$,
and Donaldson-- and Gromov-Witten-- like invariants of Fano 3-folds.
It also allows us to define the holomorphic Casson invariant of a
Calabi-Yau 3-fold $X$, prove it is deformation invariant, and
compute it explicitly in some examples. Then we calculate moduli spaces
of sheaves on a general $K3$ fibration $X$, enabling us to compute the
invariant for some ranks and Chern classes, and equate it to
Gromov-Witten invariants of the ``Mukai-dual'' 3-fold for
others. As an example the invariant is shown to distinguish Gross'
diffeomorphic 3-folds. Finally the Mukai-dual 3-fold is shown to be
Calabi-Yau and its cohomology is related to that of $X$.
\end{abstract}

\section{Introduction}

This paper is a continuation of the ideas presented in \cite{DT},
\cite{T1}. There a formal picture was outlined in which the complex
analogue of a real \emph{oriented} $n$-manifold is a \emph{Calabi-Yau}
$n$-fold with a fixed holomorphic $n$-form playing the role of a
``complex orientation'', while a Fano with fixed smooth anticanonical
divisor is the analogue of a manifold with boundary; the boundary
being the (Calabi-Yau) divisor. I have since discovered this picture
was known and used in low dimensions by the Yale school of Frenkel,
Khesin, Todorov and others (see for instance \cite{FK}, \cite{FKT},
\cite{Kh}); we concentrate on the central, three dimensional theory
of \cite{DT}. The delay in publication is due to a complete reworking
of the deformation theory of Section 3. While \cite{T1} and earlier
forms of this paper
initially used an older version of \cite{LT}, then the derived
category language of \cite{BF}, here we present more elementary sheaf
deformation theory that allows us to use the more down-to-earth
definition of obstruction theory in the published form of \cite{LT}
(without mention of the derived categories or $T^1$-lifting of older
drafts of this paper).  This deformation theory is folklore but
scattered and hard to find, and either highly abstract or not done in
enough generality for this application (namely in global families, as
the complex structure on the manifold is allowed to vary, arbitrary
order deformations and obstructions are considered, and the determinant
is fixed). So we are forced to give a full account, for which the
unpublished manuscript \cite{Fr} has been useful.

In Section 2 we review the classical Casson invariant and Chern-Simons
functional, before describing their holomorphic analogues. We give two
formulae for the holomorphic Chern-Simons functional, one in an
algebro-geometric framework, the other illustrating the complex
analogue of a manifold with boundary. We then discuss the holomorphic
Casson invariant, and a version of it that so far has not been made
rigorous but is reviewed as it motivates a number of calculations (and
predicts the correct result for them).

Section 3 discusses the technicalities involved in defining the
invariant via gauge theory, and then tackles them using algebraic
geometry, and the virtual moduli cycle theory of \cite{LT},
\cite{BF}. The result required to apply this machinery to count stable
sheaves is that the tangent-obstruction complex \cite{LT} or cut-off
cotangent complex \cite{BF} (these are dual to each other) of the
moduli space admits a certain two-step locally free resolution. We
obtain one in all the cases we might hope for, namely whenever the
higher obstruction groups Ext$^i_0(\E,\E),\ i\ge3,$ of the sheaves
$\E$ vanish.

This is the case, for instance, for ideal sheaves of curves in
$\Pee^3$, so we obtain a virtual moduli cycle of the correct dimension
in the corresponding Hilbert schemes. The result also applies to a
Calabi-Yau 3-fold, of course, allowing us to make our
definition. Finally in this section we give some examples, including
one conjectured in \cite{DT} that is fitted into this scheme and
worked out in full, by computing all reflexive sheaves of certain
Chern classes on a $(2,4)$ complete intersection in $\Pee^5$.

Section 4 deals with bundles on $K3$ fibrations. The results here may
be of interest to physicists, being the ``F-theory'' dual picture to
the geometry of $K3$-fibred 4-folds. Calculating moduli spaces of
sheaves on threefolds is extremely difficult, and we are forced to
consider only Chern classes satisfying certain constraints, though
this allows us to work on a general $K3$ fibred Calabi-Yau $X$
(without reducible or multiple fibres). Thus we obtain a calculation
of the invariant in some cases (where it is 1), and obtain
Gromov-Witten invariants of the ``Mukai-dual'' 3-fold in others. We
then show this dual 3-fold is Calabi-Yau and determine its cohomology
in terms of that of $X$. Finally we show that these results allow us
to distinguish Gross' diffeomorphic Calabi-Yau 3-folds. \newpage

\noindent \textbf{Acknowledgements.}  I am most grateful to Simon
Donaldson, whose influence is all over this paper. Conversations, help
and encouragement from Tom Bridgeland, Brian Conrad, Akira Ishii,
Adrian Langer, Jun Li, Michael McQuillan, Paul Seidel, Bernd Siebert,
Ivan Smith and Andrey Todorov,
amongst many others, have been invaluable. I would also like to thank
Bob Friedman for allowing me to see the top-secret
never-to-be-published manuscript \cite{Fr},  and H and L for
\cite{HL}, which has been indispensable.

Thanks also to the Institute for Advanced Study, Princeton, Balliol
and Hertford Colleges, Oxford, NSF (grant number DMS 9304580), EPSRC
and Professors Yau and Taubes at Harvard University for support.

\section{The holomorphic Casson invariant}

We begin by describing Taubes' version \cite{Ta} of the Casson
invariant in purely formal terms, ignoring such issues as the
structure group and reducible connections; we shall be able to bypass
these in the holomorphic situation by considering only bundles for
which semistability implies stability (such as bundles with rank and
degree coprime).

For us, then, the Casson invariant of a real 3-manifold $M$ counts
flat connections with structure group $G$ on a fixed vector bundle
$E$.  Formally the curvature $F_A$ of a connection $A$ defines a
closed one-form
\begin{equation} \label{FA}
a\mapsto\frac1{4\pi^2}\int_M\tr(a\wedge F_A), \qquad a\in\Omega^1(\ad
E),
\end{equation}
on the space of connections $\A$. This is gauge invariant, and so
descends to the space of gauge equivalence classes $\B$. Fixing a
basepoint $A_0\in\A$ this one-form is in fact the exterior derivative
of a locally defined function, the Chern-Simons functional:
$$
CS(A)=\frac1{4\pi^2}\int_M\tr\left(\frac12d_{A_0}a\wedge a+\frac13a
\wedge a\wedge a\right), \qquad A=A_0+a,
$$
which is independent of gauge transformations connected to the
identity, and well defined modulo $\Z$ on $\B$. In particular, at a
zero of the one form, i.e. a flat connection, we see that the
deformation complex of a flat connection is self-dual -- this is the
statement that the Hessian of $CS$ is symmetric -- as then are its
cohomology groups $H^i(\ad E;A)\cong H^{3-i}(\ad E;A)^*$ by Poincar\'e
duality.  Therefore the virtual dimension of the moduli space of flat
connections
$$
\sum_{i=0}^3 (-1)^{i+1}\mathrm{\,dim\,}H^i(\ad E;A)
$$
is zero, and we could hope to count them. Formally,
flat connections are the zeroes of the covector field $F_A$ (\ref{FA})
on $\B$, i.e. critical points of $CS$, and we are trying to make sense
of the Euler characteristic of the infinite dimensional space $\B$.

This formal picture translates wholesale onto a Calabi-Yau 3-fold
(which for us means a smooth, compact, K\"ahler 3-fold $X$ with
a trivialisation $\theta\in H^{3,0}$ of the canonical bundle $K_X\cong
\OO_X$). Naively, we replace $x$ by $z$, $d$ by $\db$, Poincar\'e
duality by Serre duality, and integrating against the complex volume
form $\theta$ on $X$ instead of against the orientation of $M$.

So now we consider the space $\A$ of $\db$-operators (or
``half-connections'') on a fixed $C^\infty$-bundle $E\to X$, and the
closed one-form given by $\Ft$:
$$
a\mapsto\frac1{4\pi^2}\int_X\tr(a\wedge\Ft)\wedge\theta, \qquad
a\in\Omega^{0,1}(\ad E).
$$
Again this is gauge invariant and descends to the space of gauge
equivalence classes $\B$. Fixing a basepoint $A_0\in\A$ this one-form
is the exterior derivative of a locally defined holomorphic function,
the holomorphic Chern-Simons functional:
\begin{equation} \label{CS}
CS(A)=\frac1{4\pi^2}\int_X\tr\left(\frac12\db_{A_0}a\wedge
a+\frac13a\wedge a\wedge a\right)\wedge\theta, \qquad A=A_0+a,
\end{equation}
which is independent of gauge transformations connected to the
identity.  The periods under large gauge transformations are more
complicated (and will usually be dense), but this will not concern
us. All statements will be made ``mod periods''; what we should really
do is consider $CS$ to be an element of some Albanese torus and
formulate statements there, as discussed in \cite{T1}, but this would
take us too far afield.  The zeroes of the one-form, i.e. the critical
points of $CS$, are the integrable holomorphic structures on the
bundle $E$ (i.e. $\db^2=\Ft=0$ instead of $d^2=F_A=0$), and the
deformation complex of a holomorphic connection is self-dual -- again
the statement that the Hessian of $CS$ is symmetric -- as then are its
cohomology groups $H^{0,i}(\ad E;A)\cong H^{0,3-i}(\ad E;A)^*$, by
Serre duality and the fixed trivialisation $\theta$ of the canonical
bundle $K_X$.  Therefore the virtual dimension of the moduli space of
holomorphic bundles
$$
\sum_{i=0}^3 (-1)^{i+1}\mathrm{\,dim\,}H^{0,i}(\ad E;A)
$$
is zero, and we could hope to count the bundles to formally compute
the Euler characteristic of $\B$.

Of course to get some kind of \emph{compact} space of objects to count
we have to consider either Hermitian-Yang-Mills connections in the
gauge theory set-up, or stable holomorphic bundles in algebraic
geometry; this is in some sense most of them.

Similar holomorphic analogues of all the main gauge theories also
exist \cite{DT}, \cite{T1} and have also now been studied by
physicists \cite{AOS}, \cite{BKS}, and there is work of Tyurin
\cite{Ty2} on related topics. Also, formally manipulating the
holomorphic analogue of the Chern-Simons path integral using $CS$
(\ref{CS}) gives a holomorphic linking number for complex curves in a
Calabi-Yau manifold \cite{T1} but since I have discovered the more
professional treatment of \cite{FKT}, \cite{KR}, here we concentrate
solely on defining and calculating the holomorphic analogue of the
Casson invariant. Firstly, however, we give two formulae for $CS$ to
make it more familiar and further illustrate the holomorphic analogy.

The first is in an algebro-geometric spirit, and well known on complex
curves as Abel-Jacobi theory: fixing a complex curve $\Sigma$ the
appropriate Chern-Simons functional computes holonomy of
$\db$-operators $A=\db+a$ on a topologically trivial line bundle, for
simplicity. It is
$$
a\mapsto\int_\Sigma a\wedge\omega, \qquad \omega\in H^{1,0}(\Sigma),
$$
modulo periods. Again these periods are dense and the function is only
well defined modulo a discrete lattice when considered as a function
of all $\omega\,$s at once, as an element of the torus
$$
\left(H^{1,0}\right)^*\big/\,H_1(\Sigma;\mathbb Z).
$$
In this case we also have the alternative formula
$$
\qquad\int_\Sigma a\wedge\omega=\int_\gamma\omega\qquad
\mathrm{modulo\ periods},
$$
where $\gamma$ is a path connecting the points that are the zeroes and
poles of a section that is meromorphic with respect to the holomorphic
structure $\db+a$. That is, if the holomorphic structure defined by
$a$ on the line bundle corresponds to the divisor $\sum a_i(p_i)$,
with $a_i\in\mathbb Z$ and $p_i\in\Sigma$, then $\partial\gamma=\sum
a_i(p_i)$. The principle of the result is that the Euler class of the
line bundle is represented holomorphically by both the divisor, which
is $\partial\gamma$, and the curvature $da$, and $\partial$ and $d$
are adjoints.

In a special case there is an analogous formula on a Calabi-Yau
3-fold, for a rank 2 holomorphic bundle $E$. Similarly the principle
is that $d\,CS=p_1\wedge\theta$ and, when det\,$E$ is trivial,
$p_1=c_2$ represents the Euler class of $E$, so we are interested in a
homology class $\Delta$, the analogue of $\gamma$, with boundary the
zero set of a section. (Here and below we use $CS$ to denote both the
functional (\ref{CS}) and the integrand.)

\begin{Prop} \emph{\cite{T1}} Suppose that $A_0$ and $A=A_0+a$ are
integrable ($F^{0,2}=0$) $\db$-operators on $E$ with trivial
determinant, and $(E,A),\ (E,A_0)$ admit holomorphic sections $s,\
s_0$, with transverse zero sets $(s)^{ }_0,\ (s_0)^{ }_0$.  Then $CS$
defined by (\ref{CS}) may also be described as follows.  As the zero
sets are homologous, write $(s)^{ }_0-(s_0)^{ }_0=\partial\Delta$ for
some singular 3-chain $\Delta$. Then, modulo periods,
$CS(A)-CS(A_0)=\int_\Delta\theta$.
\end{Prop}

The second formula is the complex analogue of the classical formula
for the Chern-Simons functional of a connection $A$ on a bundle $E\to
M$ on a real 3-manifold: bound $M$ by a 4-manifold $N$ and extend
$(E,A)$ to a bundle and connection $(\mathbb E,\mathbb A)$ on $N$. Then
$$
CS(A)-CS(A_0)=\int_Np_1(\mathbb A)-p_1(\mathbb A_{\,0}),
$$
where $p_1(\mathbb A)=(1/4\pi^2)\tr F_{\mathbb A}\wedge F_{\mathbb A}$
is the Chern-Weil differential form representing the first Pontryagin
class of $\mathbb E$.  This can be rephrased in a way that will make
the complex analogy more apparent in terms of the long exact homology
and cohomology sequences of the pair $(N,M)$, by the commutative
pairings
$$
\begin{array}{ccccc}
& ^{[CS(A)]} & ^{\shortmid\!\to} & ^{[p_1(\mathbb A)-p_1(\mathbb
A_{\,0})]} \\ \to & H^3(M) & \Rt{d} & H^4(N,M) & \to\,0\,
\\ & \otimes && \otimes \\ \ot & H_3(M) & \stackrel{\partial}{\Ot} &
H_4(N,M) & \ot\,0. \\ & _{[M]} & _{\otspam} & _{[N]} \end{array}
$$
That is, the fundamental class of $M$ is in the image of the lower
map, coming from the fundamental class of $N$, so to find
$\int_MCS(A)$ we can map the class $[CS(A)]$ into $H^4(N,M)$ and
evaluate on $[N]$ to give the result.

The holomorphic analogue replaces the exact sequence of the pair
$(N,M)$ by the sheaf cohomology sequence of the pair $(Y,X)$,
\begin{equation} \label{fano}
0\to K_Y\Rt{s}\OO_Y\to\OO_X\to0,
\end{equation}
in the case that $X$ is an anticanonical divisor in a 4-fold $Y$. Here
then we think of $X$ as being bounded by the ``Fano'' $Y$; we use the
term Fano loosely to mean a variety $Y$ with a section $s$ of its
anticanonical bundle $K_Y^{-1}$ with a smooth zero set $X$, which is
its ``boundary'' -- it is Calabi-Yau by the adjunction formula. In
fact about a point of $X\subset Y$, choosing a local coordinate $z$
whose zero locus is $X$, to leading order $s^{-1}=\frac1{2\pi i}
\theta\frac{dz}z$ uniquely defines a holomorphic 3-form $\theta$ on
$X$ (\cite{GH} p 147). Then we obtain

\begin{Theorem}
Suppose that the Calabi-Yau 3-fold $X$ is a smooth effective
anticanonical divisor in a 4-fold $Y$ defined by $s\in H^0(K_Y^{-1})$.
If $E\to X$ is a bundle that extends to a bundle $\mathbb E\to Y$,
then for a $\db$-operator $A$ on $E$, let $\mathbb A$ be any
$\db$-operator on $\mathbb E$ extending $A$. We have, modulo periods,
$$
CS(A)=\frac1{4\pi^2}\int_Y\tr F^{0,2}_{\mathbb A}\wedge
F^{0,2}_{\mathbb A}\wedge s^{-1}.
$$
\end{Theorem}

\begin{Proof}
Notice that $H^4(\OO_Y)=H^0(K_Y)^*=0$: if $t\in H^0(K_Y)$ then $s.t$
is a holomorphic function on $Y$ vanishing on $X$, thus $t=0$. So the
sequence (\ref{fano}) gives us the following commutative diagram of
Serre duality pairings,
\begin{equation} \label{1}
\begin{array}{ccccll}
& ^{[CS(A)]} & ^{\shortmid\!\to} & \hspace{-7mm}^{[p_1(\mathbb A)-} &
\hspace{-9mm}^{p_1(\mathbb A_{\,0})]\wedge s^{-1}} \\ \to & H^3(\OO_X)
& \To & H^4(K_Y) & \To & 0\, \\ & \otimes && \otimes \\ \ot &
H^0(\OO_X) & \Ot & H^0(\OO_Y) & \Ot & 0, \\ & _{[\,1\,]} & _{\otspam}
& _{[\,1\,]} \end{array}
\end{equation}
(the first pairing is by integrating against $\theta$) since the upper
map takes a holomorphic $(0,3)$-form on $X$, extends it to a
$C^\infty$ form on $Y$, and takes $\db(\,\cdot\,)\wedge s^{-1}$ of the
result. Setting $CS(A_0)=\int_Yp_1(\mathbb A_{\,0})$ to fix constants
gives the result.
\end{Proof}

Just as the real case $CS(A)=\int_Np_1(\mathbb A)$ can be proved
directly by Stokes' theorem, the above amounts to an application of
Stokes' theorem \emph{and} the Cauchy residue theorem, hence reducing
dimensions by two, as also observed in \cite{Kh}. If $\nu_{\delta}(X)$
denotes a small tubular neighbourhood
of $X\subset Y$ then, by Stokes' theorem,
$$
\int_Yp_1(\mathbb A)\wedge s^{-1}=\int_Yd\,(CS(\mathbb A)\wedge
s^{-1}) =\lim_{\delta\to0}\int_{\partial\nu_{\delta}}CS(\mathbb
A)\wedge s^{-1},
$$
which can be integrated first over the fibres of the circle bundle
$\partial\nu_{\delta}\to X$, and then along $X$, by Fubini's theorem.
As $s^{-1}\sim\frac1{2\pi i}\theta\frac{dz}z$, integration over the
fibres gives, by the Cauchy reside formula, $\int_XCS(A)\wedge\theta$
in the limit of $\delta\to0$.

Hence, just as $d$ is adjoint to the boundary operator $\partial$ in
real geometry, $\db(\,\cdot\,)\wedge s^{-1}$ is adjoint to this
complex boundary operation of taking the anticanonical divisor $(s)^{
}_0$ with its induced complex volume form. An application of this to
holomorphic linking is given in \cite{KR} -- the $\db$-Green's
function for the current represented by a complex curve (weighted by a
holomorphic one-form) is represented, as a current, by any Fano
surface containing it as an anticanonical divisor. Thus integrating
this over another curve (against a one-form), to give Atiyah's
holomorphic linking number, is the same as intersecting this second
curve with the complex surface, and weighting intersection numbers by
ratios of the holomorphic volume forms at intersection points. \\

Finally we mention briefly the holomorphic analogue of Casson's
original approach to counting flat connections by splitting a
3-manifold $M$ (in fact a homology sphere) across a Riemann surface
$\Sigma$,
$$
M=M_1\cup_\Sigma M_2.
$$
The orientation of $M$ induces a symplectic structure on $\Sigma$, and
so one on (the smooth locus of) its moduli space of flat connections
$\M_\Sigma$. Then those connections on $\Sigma$ that extend to flat
connections on $M_1$ form a Lagrangian submanifold in $\M_\Sigma$, the
image of the restriction map $\M_{M_1}\to\M_\Sigma$. Similarly for
$\M_{M_2}$. These are both of half dimension so we expect them to
intersect in a finite number of points -- the flat connections on
$\Sigma$ that extend to both $M_1$ and $M_2$, i.e. the flat
connections on $M$. Casson overcomes the technical difficulties to
define just such an invariant, counting (one half of) all the flat
connections except the trivial one.

Although not yet rigorous there is a holomorphic analogue of this
\cite{DT}, \cite{Kh} following work of Tyurin \cite{Ty1}. We review it
briefly because it motivates some examples and is verified in all of
them. Our complex analogue of gluing across a boundary is to take two
Fano 3-folds $X_i$ with a common anticanonical divisor $S$, and form
the normal crossings space
$$
X=X_1\cup_SX_2,
$$
which is a singular Calabi-Yau. The (singular) holomorphic volume form
on $X_i$, with poles along $S$, induces a complex symplectic structure
on the surface $S$ (this is just the adjunction formula) and so on its
moduli space of (stable) bundles $\M_S$ \cite{Mu1}. Then those bundles
on $S$ that extend to holomorphic bundles on $X_1$ form a complex
Lagrangian submanifold given by the restriction map $\M_{X_1}\to\M_S$
(at least where this is defined, i.e. where stability is preserved on
$S$); similarly for $X_2$. As before intersecting these
\begin{equation} \label{tyc}
\M_{X_1}\cap\M_{X_2}
\end{equation}
in $\M_S$ we expect to get a finite number of holomorphic bundles that
extend to both $X_i$, i.e. bundles on $X$. In the examples we consider
$X$ will be smoothable and the number of bundles will be preserved on
smoothing to give the holomorphic Casson invariant of the smooth
3-fold.

\section{Virtual moduli cycles}

To count (stable) holomorphic bundles on a Calabi-Yau 3-fold, there
are two things we require of the moduli problem -- compactness and
transversality. In gauge theory such results are easier in lower
dimensions. In two real dimensions moduli of stable bundles are both
compact and of the right dimension. In three dimensions this can be
achieved after a perturbation \cite{Ta} (leaving aside problems with
reducibles) for flat connections. In dimension four we need both
perturbations to achieve transversality, and a compactification to
take account of the non-compactness caused by conformal invariance
\cite{DK}.  Just recently we now have the results of Tian \cite{Ti} in
higher dimensions, proving just about everything that one would like
to be true, giving a natural analogue of the Donaldson-Uhlenbeck
compactification of four dimensions. For a K\"ahler 3-fold this
involves ideal instanton singularities along holomorphic curves in the
3-fold, but also some codimension 3 singularities that are harder to
deal with.

What is missing, however, is a transversality result. Staying within
the confines of algebraic Calabi-Yau manifolds we cannot hope to get a
moduli space of the correct dimension; there is no result along the
lines of Donaldson's generic smoothness result for moduli spaces on
algebraic surfaces. One reason is that, in the rank two case for
instance, relating the deformation theory of a bundle to that of a
zero set of a section via the Serre construction, points are
unobstructed on a surface but curves can be obstructed on a 3-fold
(there is a fuller discussion of this and other such issues in
\cite{T1}, and an example in \cite{T2}).

So we would like more perturbations. There is an elliptic perturbation
of the Hermitian-Yang-Mills equation valid on any almost-complex
symplectic manifold, which I learnt from Donaldson:
\begin{eqnarray} \nonumber
F^{0,2}_A &\!=\!& \db^*u, \\ \iL &\!=\!& \lambda I. \nonumber
\end{eqnarray}
(The problem is that the Hermitian-Yang-Mills equations appear over
determined, but are not because of the Bianchi identity $\db
F^{0,2}_A=0$ on a K\"ahler manifold.  This can be formalised by
introducing the (0,3)-form $u$ as above, and then $\db^*u$
vanishes. It need not be zero in the almost complex case.)

It seems that Tian's work \cite{Ti} should also apply to these
equations (with the singularities now along pseudo-holomorphic curves
and at points) so long as we can get a bound on $\|F_A\|_{L^2}^2$
similar to that given by characteristic class formulae in the
integrable case. A Weitzenb\"ock formula shows that this is the case
if, for instance, the scalar curvature is everywhere positive. Thus,
at present, this would work best for Fano 3-folds and the Calabi-Yau
case is borderline. Either way a transversality result, proving that
for generic almost complex structures the moduli space of solutions is
of the correct dimension, seems a long way off, as does understanding
the codimension three singularities. \\

So we turn to algebraic geometry where we have the now standard
compactification of the moduli space of stable bundles by semistable
sheaves, due to Gieseker, Maruyama and Simpson. This moduli space will
invariably be singular and of too high a dimension, but it is often
clear what the contribution of a particular component to the ``number
of bundles'' should be -- the Euler number of its cotangent bundle in
the case it is smooth, two if it is a scheme-theoretic double point,
etc. In the general case there is the machinery of \cite{LT},
\cite{BF} to produce a ``virtual moduli cycle'' of the correct
dimension (zero, for us) inside any moduli space satisfying certain
conditions; we briefly outline the picture.

Suppose a variety $M$ (which will eventually be our moduli space $\M$)
sits inside a smooth ambient $n$-fold $Z$, cut out by a section $s$ of
a rank $r$ vector bundle $E\to Z$. Then the ``virtual dimension'' of
$M$ is $(n-r)$ -- the dimension it would be were $s$ transverse.  If
it is not but, for instance, $s$ lies in a subbundle $E'\subset E$ and
is a transverse section of $E'$, then it is clear the ``correct''
$(n-r)$-cycle we should take is the Euler class of the cokernel bundle
$E/E'$ over $M$ -- this is homologous to the zero set of a transverse
perturbation of $s$ if one exists. In the general case dealt with by
Fulton-MacPherson intersection theory \cite{Fu}, $s$ induces a cone in
$E\res M$, which can be thought of as ``$s$ made vertical'', i.e.  the
limit of the images of $\lambda s$ as $\lambda\to\infty$. We may then
intersect this with the zero set $M$ inside $E$ to get a cycle in $\M$
(whose image in $Z$ represents the top chern class of $E$, as required).

The point here is that we worked entirely on $M$ and not in the
ambient space $Z$, and so we might hope the method is applicable to
moduli problems where the ambient space $Z$ does not exist. Instead
the deformation theory of the moduli problem often gives us the
infinitesimal version of $(Z,\,E,\,s)$ on $M$, namely the derivative
of $s$, yielding the exact sequence
\begin{equation} \label{sm}
0\to TM\to TZ\res M\Rt{\!ds}E\res M\to\mathrm{ob}\to0,
\end{equation}
for some cokernel ob which in the moduli problem becomes the
obstruction sheaf. In the general case we require a global version of
this, namely a two term locally free resolution
\begin{equation} \label{pre}
0\to\T_1\to E_1\to E_2\to\T_2\to0,
\end{equation}
of the tangent-obstruction functors to be introduced presently
(\ref{t}, \ref{o}).

Here $E_1$ and $E_2$ play the roles of $TZ\res M$ and $E\res M$ in the
above motivation (these last two have the same fibre rank at each
point of $\M$, and hence are locally free) and have difference in
ranks equal to the virtual dimension of the moduli problem. It is
shown in \cite{LT}, \cite{BF} that such data on $\M$ is in fact
sufficient to obtain a cone in the vector bundle $E_1$ which can be
intersected with the zero set $\M$ to give a virtual moduli cycle with
the correct properties. The precise statement is given below. \\

First then, we need to develop the necessary sheaf deformation theory.
An earlier version of this paper used the cotangent complex approach
of \cite{BF}, and Lehn's description \cite{L} of equations cutting
out the Quot scheme to calculate it; the deformation theory here is
more classical and natural, even if it is a little longer, dealing
with higher order deformations.

An excellent reference for the sheaf theory we use is \cite{HL}; we
shall assume familiarity with Gieseker stability (here always referred
to just as stability), slope stability, the fact that
$$
\mathrm{slope\ stable\ }\Rightarrow\mathrm{\ stable\ }\Rightarrow
\mathrm{\ semistable\ }\Rightarrow\mathrm{\ slope\ semistable,}
$$
and that for rank and degree coprime the circle is completed by
slope semistability implying slope stability. Also refer to \cite{HL}
for the fact that for either form of stability
$$
\mathrm{(slope)\ stable\ }\Rightarrow\mathrm{\ simple,}
$$
i.e. the only endomorphisms of the sheaf are the scalars $\C\,.\,$id.
Finally, stable sheaves are pure, i.e. torsion-free on restriction to
their support.

\subsection{Sheaf deformation theory}

In this section all schemes will be complex and quasi-projective, and
all sheaves coherent. We begin by recalling the definition of the
trivial thickening
\begin{equation} \label{thick}
S*\n
\end{equation}
of a scheme $S$ by a coherent sheaf of $\OO_S$-modules $\n$. Namely
make the $\OO_S$-module $\OO_S\oplus\n$ into a sheaf of rings by
stipulating that $\n^2=0$, giving the trivial ring extension
$\OO_S*\n$; the associated scheme is $S*\n$. Then sheaf deformation
theory is built on the following standard result.

\begin{Lemma} \label{defs}
Let $\n$ be a sheaf on a scheme $S$, and let $S*\n$ be the trivial
extension of (\ref{thick}) above. Then deformations of a sheaf $\E$ on
$X\x S$, flat over $S$, to a sheaf on $X\x(S*\n)$, flat over
$S*\n$, are in 1-1 correspondence with $\mathrm{Ext}^1_{X\x
S}(\E,\E\otimes\n)$.
\end{Lemma}

\begin{Proof}
Given such a deformation $\F$, tensoring with the sequence
\begin{equation} \label{i}
0\to\n\to\OO_{S*\n}\to\OO_S\to0,
\end{equation}
gives a sequence, exact by flatness,
\begin{equation} \label{j}
0\to\E\otimes\n\to\F\to\E\to0
\end{equation}
on $X\x(S*\n)$, using the fact that the left and right hand terms
of (\ref{i}) are $\OO_S$-modules. The sequence (\ref{j}) defines the
class in Ext$^1_{X\x S}(\E,\E\otimes\n)$.

Conversely such a class gives a sequence of $\OO_{X\x S}$-modules
(\ref{j}), which defines an $\OO_{X\x(S*\n)}$-module $\F$ since
there is an obvious action of $\OO_S*\n$ on $\F$: the $\n$-action is
given by projecting $\F\to\E$ and tensoring with $\n$, mapping to
$\E\otimes\n\subset\F$. Flatness over $X\x(S*\n)$ also follows
from the sequence (\ref{j}) and the following Lemma.
\end{Proof}

\begin{Lemma} \label{flat}
Let $S\subset Y$ be a subscheme with ideal $\n\subset\OO_Y$ such that
$\n^2=0$. Then an $\OO_Y$-module $\F$ is flat over $Y$ if and only if
$\F\otimes\n\to\F$ is injective and $\F\res S=\F\otimes\OO_S$ is flat
over $S$.
\end{Lemma}

\begin{Proof}
We must show (\cite{Ha1} III 9.1A) that $\F\otimes J\to\F$ is
injective for any ideal $J\subset\OO_Y$. Given such a $J$, we have an
exact sequence
$$
0\to J\cap\n\to J\to J'\to0,
$$
where $J'$ is annihilated by $\n$ so is naturally an ideal in $\OO_S$.
Thus the diagram
\spreaddiagramrows{-0.8pc}
\spreaddiagramcolumns{-0.8pc}
$$ \diagram
&&& 0 \dto \\ & \F\otimes(J\cap\n) \rto\dto & \F\otimes J
\rto\dto & \F\res S\otimes J' \rto\dto & 0 \\ 0 \rto & \F\otimes\n
\rto & \F \rto & \F\res S \rto & 0
\enddiagram
$$
is exact by our hypotheses. To show the central vertical map is
injective, then, it is sufficient to show that the first vertical map
is injective. But this map is $\F\res S\otimes(J\cap\n)\to\F\res
S\otimes\n$, through which $\F\res S\otimes(J\cap\n)\to\F\res S$
factors, and $\F\res S$ is assumed flat so this last map is an
injection.
\end{Proof}

Next we consider obstructions to deformations. We will repeatedly use
the following set-up.
\begin{description}
\item[\ \ $\bullet$\ \ ] Let $S\subset Y\subset Y_1$ all be schemes
\emph{over $S$}, and denote the ideals of $S\subset Y,\ S\subset Y_1,
\ Y\subset Y_1$ by $\n,\,\m$ and $\I$ respectively. Assume also that
$\m\,.\,\I=0$, giving an exact sequence of $\OO_Y$-modules
\begin{equation} \label{ideals}
0\to\I\to\m\to\n\to0.
\end{equation}
\end{description}
Notice we do \emph{not} assume that $Y=S*\n$ this time, or even
that $\n^2=0$.

Given a sheaf $\E$ over $X\x Y$ (flat over $Y$ and restricting to
$\E_0$ over $S$) we get an exact sequence
\begin{equation} \label{En}
0\to\E\otimes\n\to\E\to\E_0\to0,
\end{equation}
giving a class $e\in$\,Ext$^1_{X\x S}(\E_0,\E\otimes\n)$ (we will consider
all terms as $\OO_S$-modules using the projections $Y\to S,\ Y_1\to S$,
so all Exts will be over $X\x S$ from now on).

We would like to lift this to an $\F$ on $Y_1$ to give a sequence
\begin{equation} \label{Em}
0\to\E\otimes\m\to\F\to\E_0\to0,
\end{equation}
(here $\F\otimes\m\cong\E\otimes\m$ since $\m\,.\,\I=0$) defining
$f\in$\,Ext$^1(\E_0,\E\otimes\m)$. $f$ is a lift of $e$ in
the sequence
\begin{equation} \label{dagger}
\mathrm{Ext}^1(\E_0,\E\otimes\m)\to\mathrm{Ext}^1(\E_0,\E\otimes\n)
\Rt{\partial}\mathrm{Ext}^2(\E_0,\E_0\otimes\I),
\end{equation}
obtained by applying Hom\,$(\E_0,\ .\ )$ to the sequence
$\E\otimes\,$(\ref{ideals}):
\begin{equation} \label{Eideals}
0\to\E_0\otimes\I\to\E\otimes\m\to\E\otimes\n\to0.
\end{equation}
(Exactness follows from the flatness of $\E$ over $Y$.)

\begin{Prop} \label{obsclass}
$\E$ over $X\x Y$ as above extends to a sheaf $\F$ over $X\x
Y_1$, flat over $Y_1$, if and only if there is a lift $f\in
\mathrm{Ext}^1 (\E_0,\E\otimes\m)$ \emph{(\ref{Em})} of $e\in
\mathrm{Ext}^1(\E_0,\E\otimes\n)$ \emph{(\ref{En})}, i.e. if and only
if $\partial e\in\mathrm{Ext}^2(\E_0, \E_0\otimes\I)$ in the above
sequence \emph{(\ref{dagger})} is zero.
\end{Prop}

\begin{Proof}
We are left with showing that the existence of a lift $f$ of $e$ gives
such an $\F$. $f$ gives a sequence (\ref{En}) \emph{of
$\OO_S$-modules} lifting (\ref{Em}), and so a diagram
$$
\spreaddiagramrows{-0.8pc}
\spreaddiagramcolumns{-0.8pc}
\diagram
0 \rto & \E\otimes\m \rto^(.6)\iota\dto & \F \rto\dto^\pi & \E_0
\rto\ar@{=}[d] & 0 \\
0 \rto & \E\otimes\n \rto\dto & \E \rto\dto & \E_0 \rto & 0 \\ &
0 & \ \,0\,.
\enddiagram
$$
Here $\F$ is an $\OO_S$-module, and we would like to make it an
$\OO_{Y_1}$-module, where $\OO_{Y_1}\cong\OO_S\oplus\m$ via the maps
$S\rightleftarrows Y_1$. So we define the action $\F\otimes\m\to \F$
by $\iota\circ(\pi\otimes\id)$ in the above diagram.

Flatness of $\F$ over $Y_1$ then follows from Lemma \ref{flat} on
noting that $\F\res{Y}=\E$ is flat over $Y$, and $\I^2=0$ since
$\I\subset\m$ and $\m\,.\,\I=0$.
\end{Proof}

We note in passing that $\partial$ is cup product with the element
$e'\in$\,Ext$^1(\E\otimes\n,\E_0\otimes\I)$ defining the extension
(\ref{Eideals}), so $\E$ extends to $\F$ if and only if $e'\cup e
\in$\,Ext$^2(\E_0,\E_0\otimes\I)$ is zero.

\subsection*{The trace map}

We have now pretty much found the tangent-obstruction (in the sense of
\cite{LT}) complex of the moduli problem for stable sheaves, as we
shall see below.  Unfortunately we are more interested in the moduli
problem for sheaves of fixed determinant, for which we need the
machinery of the Mukai-Artamkin trace map (in a little more generality
than \cite{Mu1}, \cite{Ar}, for higher order deformations).

Given a coherent sheaf $\F$ on a quasi-projective scheme $X$ and an
affine open cover $\mathcal U=\{U_i:\,i=1,\ldots,n\}$ of $X$, denote
the \v Cech complex by
$$
\check C^p(\F)=\prod_{i_0<\ldots<i_p}\Gamma(\F\res{U_{i_0}\cap\cdots
\cap U_{i_p}}),
$$
with the usual \v Cech differential $\delta:\,\check C^p\to\check
C^{p+1}$ (\cite{Ha1} III 4.1). This computes the sheaf cohomology of
$\F$, and as the construction is functorial in $\F$ it can also be
applied to a complex of sheaves $\F\udot$ to give a double complex
whose associated total complex has cohomology the hypercohomology
$\mathbb H^{\,*}(\F\udot)$.

If $X$ is smooth then any sheaf $\E$ has a finite locally free
resolution $E\udot\to\E\to0$, and the trace map is defined as follows
(by standard arguments it will be independent, up to
quasi-isomorphism, of the choices $\mathcal U$ and $E\udot$).

Given any sheaf $\II$, form the complex
$\Homd(E\udot,E\udot\otimes\II)$ with
$\Hom^i(E\udot,E\udot\otimes\II)=\oplus_j\Hom(E^j,E^{i+j}\otimes\II)$
and differential $d\phi=d^{i+j}\circ\phi-(-1)^i\phi\circ d^j$ for
$\phi \in\Hom(E^j,E^{i+j}\otimes\II)$. This admits cochain maps
\begin{equation} \label{trid}
\Homd(E\udot,E\udot\otimes\II)\,\rightleftarrows\,\II,
\end{equation}
with the upper map given by
$$
\tr=\sum_i(-1)^i\,\tr^i\otimes\id_\II
$$
(where $\tr^i:\,\Hom(E^i,E^i)\to\OO$ is the usual trace map on locally
free sheaves), and the lower map
$$
\id=\sum_i\id_{E^i}\otimes\id_\II.
$$
That these are cochain maps follows from the easy computations
$\tr\circ\,d=0$ and $d\circ\id=0$.

Notice that $\tr\,\circ\,\id=\sum(-1)^i\,\tr^i\,\circ\,\id_{E^i}=
\sum(-1)^i\,\mathrm{rk}\,(E^i)=\mathrm{rk}\,(\E)$, so, for
$\mathrm{rk}\,(\E)>0$, we have a splitting
$$
\Homd(E\udot,E\udot\otimes\II)=\Homd_0(E\udot,E\udot\otimes\II)
\,\oplus\,\II,
$$
where $\Homd_0$ is the kernel of tr.

Thus (\ref{trid}) induces cochain maps between \v Cech complexes
$$
\check C\udot(\Homd(E\udot,E\udot\otimes\II))\,\rightleftarrows\,
\check C\udot(\II),
$$
inducing maps tr and id on cohomology
$$
\mathrm{Ext}^i(\E,\E\otimes\II)\,\rightleftarrows\, H^i(\II),
$$
such that $\tr\,\circ\,\id=\mathrm{rk}\,\E$. Thus for
$\mathrm{rk}\,(\E)>0$, there is a splitting
$\mathrm{Ext}^i(\E,\E\otimes\II)=\mathrm{Ext}^i_0
(\E,\E\otimes\II)\,\oplus\,H^i(\II)$ with $\mathrm{Ext}_0$ the kernel
of tr.

To work towards showing that taking the trace of deformations and
obstructions of a sheaf gives the deformations and obstructions of the
determinant of the sheaf, we first need this (rather elementary) fact
for locally free sheaves. Of course, phrasing the deformation theory
of holomorphic vector bundles in terms of connections or transition
functions this is simple; the work is then in showing the deformation
theory coincides with the abstract sheaf deformation theory of the
last section. This involves simple but very large computations with
transition functions as in \cite{Fr}. Here we prefer to work directly
with our definitions above; this then makes the proof below a little
long, but the reader could take it on trust.

\begin{Prop} \label{LF}
Take $S\subset Y\subset Y_1$ to be as in (\ref{ideals}).
Suppose we have a rank $r$ \emph{locally free}
sheaf $E$ on $X\x Y$, a flat deformation of $E_0$ on
$X\x S$ giving an extension $e\in \mathrm{Ext}^1(E_0,E\otimes\n)$
(\ref{En}). Assuming
first that $\n^2=0$, then the extension defined similarly by the
determinant $\Lambda^rE$ is
$$
\tr(e)\in H^1(\n)=\mathrm{Ext}^1(\Lambda^rE_0,\Lambda^rE_0\otimes\n).
$$
Likewise, for any $\n$, given the obstruction $\partial
e\in\mathrm{Ext}^2(E_0,E_0\otimes\I)$ of Proposition
\ref{obsclass} to extending $E$ over $X\x Y$ to $F$ over $X\x
Y_1$ (flat over $Y_1$), the obstruction to extending $\Lambda^rE$ is
given by
$$
\tr(\partial e)\in H^2(\I)=\mathrm{Ext}^2(\Lambda^rE_0,
\Lambda^rE_0\otimes\I).
$$
\end{Prop}

\begin{Proof}
We begin by giving explicit descriptions of $e$ and the obstruction
$\partial e$.

Applying Hom$(E_0,\ .\ )$ to the extension (of $\OO_{X\x
S}$-modules)
\begin{equation} \label{E0I}
0\to E\otimes\n\to E\to E_0\to0
\end{equation}
gives a connecting homomorphism $\mathrm{Hom}(E_0,E_0)\to
\mathrm{Ext}^1(E_0,E\otimes\n)$ under which the image of $\id_{E_0}$
is the extension class $e\in\mathrm{Ext}^1(E_0,E\otimes\n)$.

So consider the exact sequence of complexes below given by applying
the exact functor $\Hom(E_0,\ .\ )$ (recall that $E_0$ is a locally
free $\OO_S$-module) to the sequence (\ref{E0I}) and taking \v Cech
complexes:
\begin{equation} \label{cc}
0\to\check C\udot(\Hom(E_0,E\otimes\n))\Rt{\iota} \check
C\udot(\Hom(E_0,E))\to\check C\udot(\Hom(E_0,E_0))\to0.
\end{equation}
$\id_{E_0}$ gives a closed element of $\check C^0(\Hom(E_0,E_0))$
(i.e.  a global section). Lift this to some $a\in\check
C^0(\Hom(E_0,E))$, giving, for each $U,\,V$ in some affine open cover
$\mathcal U$ of $X$, $a_U,\,a_V$ such that over $U\cap V$, $a_V-a_U=:
\iota(e_{U\cap V})$ defines the (closed) element
$$
e\in\check C^1(\Hom(E_0,E\otimes\n))
$$
that represents $e\in\mathrm{Ext}^1(E_0,E\otimes\n)$.

To identify $\partial e$ we use the similar exact sequence of \v Cech
complexes
$$
0\to\check C\udot(\Hom(E_0,E_0\otimes\I))\Rt{\iota}\check
C\udot(\Hom(E_0,E\otimes\m))\to\check C\udot(\Hom(E_0,E\otimes\n))\to0,
$$
associated to the sequence (\ref{dagger}).

That is, choose lifts $b_{UV}\in\mathrm{Hom}_{\,UV}(E_0,E\otimes\m)$
of $e_{UV}$. Then over $U\cap V\cap W$, $b_{UV}+b_{VW}+b_{WU}=
\iota((\partial e)_{UVW})$ defines $(\partial e)_{UVW}\in
\mathrm{Hom}_{\,UVW}(E_0,E_0\otimes\I)$ giving the obstruction class
$\partial e\in\mathrm{Ext}^2(E_0,E_0\otimes\I)$. \\

We now use these explicit calculations to repeat the above exercise on
the induced deformation of determinants
\begin{equation} \label{dets}
0\to\Lambda^rE\otimes\n\to\Lambda^rE\to\Lambda^rE_0\to0.
\end{equation}
Here for $\Lambda^r E$ we have taken the $r$th exterior power of $E$
\emph{as a sheaf of $\OO_{X\x Y}$-modules} (\emph{not} as a sheaf
of $\OO_{X\x S}$-modules) but then we consider the result and the
above sequence as $\OO_{X\x S}$-modules. Thus we have the analogue
$$
0\!\to\!\check C\udot(\Hom(\Lambda^rE_0,\Lambda^rE\otimes\n))
\rt{\iota}\check C\udot(\Hom(\Lambda^rE_0,\Lambda^rE))
\!\to\!\check C\udot(\Hom(\Lambda^rE_0,\Lambda^rE_0))\!\to\!0
$$
of (\ref{cc}).

Use $\Lambda^ra_U\in\mathrm{Hom}_{\,U}(\Lambda^rE_0,\Lambda^rE)$ to
lift $\id_{\Lambda^rE_0}=\Lambda^r\id_{E_0}\in\mathrm{Hom}_{\,U}
(\Lambda^rE_0,\Lambda^rE_0)$.  Over $U\cap V$,
$\Lambda^ra_U-\Lambda^ra_V=\iota(\epsilon_{UV})$ defines
$\epsilon_{UV}\in\mathrm{Hom}_{\,UV}(\Lambda^rE_0,\Lambda^rE
\otimes\n)$, the extension class $\epsilon\in\,$Ext$^1(\Lambda^rE_0,
\Lambda^rE\otimes\n)=H^1(\n)$ of (\ref{dets}).

Over $U\cap V$ there is a splitting of (\ref{E0I}) induced by $a_U$:
\begin{equation} \label{split}
E\cong(E\otimes\n)\oplus E_0,
\end{equation}
with respect to which $a_U=0\oplus\id_{E_0}$, $a_V=e_{U\cap
V}\oplus\id_{E_0}$, and (\ref{dets}) splits as
$\Lambda^rE\cong(\Lambda^rE\otimes\n)\oplus\Lambda^rE_0$ (we are
omitting some $\res{U\cap V}$\,s for clarity).

So over $U\cap V$, $\epsilon_{UV}=\Lambda^r
(e_{UV}\oplus\id)-\id,\ \epsilon_{UW}=\Lambda^r \,(e_{UW}\oplus\id)
-\id$ with respect to this splitting, while $\epsilon_{VW}=
\Lambda^r(e_{UW}\oplus\id)-\Lambda^r(e_{UV}\oplus\id)$.

For the first part of the Proposition we want to show that, for
$\n^2=0$, these $\epsilon$\,s are (the images under $\iota$ of) the
traces of the corresponding $e$\,s (i.e. $\epsilon_{UV}=\tr(e_{UV})$,
etc.). But this is clear from the expansion
$$
\Lambda^r(\id_{E_0}\oplus e)=\id_{\Lambda^rE_0}+\tr_U(e)+2\tr_U
(\Lambda^2e)+\ldots,
$$
for $e:\,E_0\to E\otimes\n$, where $\tr_U(\Lambda^ke)$ is a map of the
right kind, i.e. an element of
$(\Lambda^rE_0)^*\otimes\Lambda^rE\otimes\n$, on defining $\tr_U$ by
the composition
$$
\Lambda^kE_0^*\otimes\Lambda^kE\otimes\n\to\Lambda^kE_0^*\otimes
\Lambda^kE_0\otimes\n\to\n\to\Lambda^rE_0^*\otimes\Lambda^rE_0\otimes
\n\to\Lambda^rE_0^*\otimes\Lambda^rE\otimes\n
$$
of the projection, trace, identity and $\Lambda^ra_U$ maps
respectively. This is just a glorified version of the expansion of the
determinant in terms of the elementary symmetric polynomials of $e$,
but over the ring $\OO_Y$. For $\n^2=0$ it is just the usual trace,
independent of $a_U$, since any two $a$\,s differ by something in the
ideal $\n$. In this case all the higher order terms in the above
expansion become zero anyway leaving just $\tr(e)$, giving the
required result.

To identify the obstruction we pick the obvious lifts of the
$\epsilon$\,s. That is, over $U\cap V$ in the splitting (\ref{split}),
$\epsilon_{UV}=\tr_U(e_{UV})+2\,\tr_U(\Lambda^2e_{UV})+\ldots$, so
we set $\beta_{UV}=\tr_U(b_{UV})+2\,\tr_U(\Lambda^2b_{UV})+\ldots\,$.
($\tr_U$ is defined by the same formula as before but with $\m$
replacing $\n$. Then $\beta_{UV}$ is actually skew symmetric with
respect to $U$ and $V$, so is well defined, though it takes a
calculation relating $\tr_U$ and $\tr_V$ to check it.)

The class of the obstruction we are seeking is given, on $U\cap
V\cap W$, by
\begin{equation} \label{trobs}
\beta_{UV}+\beta_{VW}+\beta_{WU}=\tr_U(b_{UV})+\tr_V(b_{VW})+
\tr_W(b_{WU})+\ldots.
\end{equation}
We claim that only the terms with a linear dependence on the $b$\,s
are non-zero up to coboundaries; this can be checked
by a large local calculation or the following cheat. We know that
$b_{VW}= -b_{UV}-b_{WU}+i$, where $i$ lies in the ideal $\I$ ($i$
is of course $\iota(\partial e)_{UVW})$.) If we had chosen the
different lift $-b_{UV}-b_{WU}$ of $e_{VW}$ then (\ref{trobs}) would
of course give zero up to coboundaries. Thus (\ref{trobs}),
considered as a polynomial in $i$, has zero constant term, and any
term of order $\ge2$ (in the $b$\,s or $i$) involves $i$ multiplying
something in $\m$, which vanishes since $\m\,.\,\I=0$.  Thus the
obstruction class is the just the linear part of (\ref{trobs}), which
is $\tr((\partial e)_{UVW})$ since
$b_{UV}+b_{VW}+b_{WU}=\iota((\partial e)_{UVW})$.
\end{Proof}

Similar explicit calculations over patches of a cover $\mathcal U$
also give the standard results that, for line bundles $L_i$ with
deformation and obstruction classes $e(L_i)$ and $\partial(e(L_i))$ as
above, we have
\begin{equation} \label{tensor}
e(\otimes_iL_i)=\sum_ie(L_i)\in H^1(\n), \qquad\mathrm{and}\qquad
\partial(e(\otimes_iL_i))=\sum_i\partial(e(L_i))\in H^2(\n).
\end{equation}

Also recall that for any sheaf $\E$ with a finite locally free
resolution $E\udot$, the determinant of $\E$ is defined to be the line
bundle
\begin{equation} \label{det}
\det\E=\bigotimes_i(\det E^i)^{(-1)^i},
\end{equation}
which is independent of the resolution.

So as before let $X$ be a smooth quasi-projective variety and let
$S\subset Y\subset Y_1$ be as in (\ref{ideals}). Fix a sheaf $\E$ on
$X\x Y$ that is flat over $Y$ and restricts to $\E_0$ on $X\x S$. Then
since $X$ is smooth, $\E$ has such a finite locally free resolution
$E\udot$ that restricts on $X\x S$, by flatness, to a finite locally
free resolution $E_0\udot$ of $\E_0$.

\begin{Theorem} \label{tr}
In the above situation, denote by
$e\in\mathrm{Ext}^1(\E_0,\E\otimes\n)$ and $\partial
e\in\mathrm{Ext}^2(\E_0,\E_0\otimes\I)$ the deformation and
obstruction classes of $\E$ (\ref{En},\,\ref{obsclass}).  Then
the obstruction class of $\det\E$ is $\tr(\partial e) \in
H^2(\I)=\mathrm{Ext}^2(\det\E_0,\det\E_0
\otimes\I)$, and, if $\n^2=0$, the deformation class of $\det\,
\E$ is $\tr(e)\in H^1(\n)=\mathrm{Ext}^1(\det\E_0,
\det\E_0\otimes\n)$.
\end{Theorem}

\begin{Proof}
Again we explicitly chase \v Cech cocycles representing the extension
and deformation classes. The exact sequence $0\to\n \to\OO_{X\x
Y}\to\OO_{X\x S}\to0$ gives the exact diagram of resolutions
\spreaddiagramrows{-0.6pc}
\spreaddiagramcolumns{-0.4pc}
$$ \diagram 0 \rto & E\udot\otimes\n \rto\dto & E\udot \rto\dto &
E_0\udot \rto\dto & 0 \\ 0 \rto & \E\otimes\n \rto\dto & \E \rto\dto &
\E_0 \rto\dto & 0 \\ & 0 & 0 & \ \,0\,.  \enddiagram
$$
In turn this gives an exact sequence of \v Cech complexes
$$
0\to\check
C\udot(\Hom\udot(E_0\udot,E\udot\otimes\n))\rt{\iota}
\check C\udot(\Hom\udot(E_0\udot,E\udot))\to\check C\udot(\Hom\udot
(E_0\udot,E_0\udot))\to0,
$$
where the first complex computes Ext$^1(\E_0,\E\otimes\n)$, etc.
$\id\in\check C^0(\Hom\udot(E_0\udot,E_0\udot))$ is $\oplus_i\id_i$
with $\id_i\in\check C^0(\Hom^0(E_0^i,E_0^i))$.  Lift this to
$a\in\check C^0(\Hom\udot(E_0\udot,E\udot))$ as $\oplus_ia_i$ with
$a_i\in\check C^0(\Hom^0(E_0^i,E^i))$.

We now apply the total differential $d+\delta$ to $a$, with $d$ the
differential on $\Hom\udot$ and $\delta$ the \v Cech differential,
keeping track of degrees. We get
$$
(\delta+d)\,a=\iota\,(\oplus_ie_i+\oplus_if_i),
$$
where $e_i\in\check C^1(\Hom^0(E_0^i,E^i\otimes\n))$ and $f_i\in
\check C^0(\Hom(E_0^i,E^{i-1}\otimes\n))\oplus\check C^0
(\Hom(E_0^{i+1},\newline E^i\otimes\n))$.  Thus
$\iota(e_i)=\delta(a_i)$ (and $\iota(f_i)=d(a_i)$).

So $\oplus_ie_i+\oplus_if_i$ is closed under $d+\delta$ and
represents the class $e\in\mathrm{Ext}^1(\E_0,\E\otimes\n)$ of the
deformation $\E$, whereas $e_i$ is closed under $\delta$ and
represents the class in $\mathrm{Ext}^1(E_0^i,E^i\otimes\n)$ of the
deformation $E^i$. Suppose that $\n^2=0$ so that $E\otimes\n=E_0
\otimes\n$, etc. and $\tr$ is defined. Then since $\tr$ only acts on
$\Hom^0$ components of the complex (\ref{trid}), we see that
$$
\tr(e)=\tr(\oplus_ie_i+\oplus_if_i)=\tr(\oplus_ie_i)=
\sum_i(-1)^i\tr(e(E^i))=\sum(-1)^i\,e(\det E^i),
$$
by Proposition \ref{LF}. By (\ref{tensor}) this is
$e(\det\E)$, as required.

To deal with obstructions use $0\to\I\to\m\to\n\to0$ to give the exact
diagram
$$ \diagram 0 \rto & E_0\udot\otimes\I \rto\dto & E\udot\otimes\m
\rto\dto & E\udot\otimes\n \rto\dto & 0 \\ 0 \rto & \E_0\otimes\I
\rto\dto & \E\otimes\m \rto\dto & \E\otimes\n \rto\dto & 0 \\ & 0 & 0
& \ \,0\,, \enddiagram
$$
giving the exact sequence of \v Cech complexes
$$
0\to\check
C\udot(\Hom\udot(E_0\udot,E_0\udot\otimes\I))\rt{\iota}
\check C\udot(\Hom\udot(E_0\udot,E\udot\otimes\m))\to\check
C\udot(\Hom\udot(E_0\udot,E\udot\otimes\n))\to0.
$$
Lift our class $e=\oplus_ie_i+\oplus_if_i$ to
$\oplus_ib_i+\oplus_ic_i$ in $\check
C\udot(\Hom\udot(E_0\udot,E\udot\otimes\n))$, and apply $\delta+d$ to
give, in different degrees,
$$
\oplus_i\delta(b_i)\ +\ \oplus_i(\delta(c_i)+d(b_i))\ +\
\oplus_id(c_i).
$$
This is $\iota(\partial e)$, by definition, with the component of
$\partial e$ in $\check C^2(\Hom^0)$ being a sum $\oplus_io_i$ of
terms such that $\iota(o_i)=\delta(b_i)$. Thus $o_i$ is, by the
definition of $b_i$, the obstruction $\partial\,e(E^i)$ to the
extension of $E^i$, and since $\tr$ only acts on $\Hom^0$ components
of the complex (\ref{trid}), i.e. on only the $\oplus_i\partial
e(E^i)$ parts of $\partial e$, we have
$$
\tr(\partial e)=\tr(\oplus_i\partial\,e(E^i))=\sum_i(-1)^i\tr (\partial
\,e(E^i))=\sum_i(-1)^i\partial\,e(\det E^i),
$$
by Proposition \ref{LF}. But by (\ref{tensor}) this is $\partial\,
e(\det\E)$, as required.
\end{Proof}

We are finally in a position to find the ``tangent-obstruction
complex'' of our moduli problem, as defined in \cite{LT} (though our
moduli problem is contravariant, not covariant).

Fix a smooth quasi-projective scheme $X$ and Chern classes $c_i\in
H^{2i}(X)$, and consider the moduli functor $\M$ that assigns to any
scheme $S$ the set of isomorphism classes of sheaves $\E$ on $X\x S$,
flat over $S$, whose restriction to each fibre is \emph{stable and has
Chern classes $c_i$}. Here two sheaves are considered isomorphic if
they differ by tensoring with a line bundle on $S$. This moduli
functor has a coarse moduli space which we also denote by
$\M=\M(X,c_i)$ \cite{HL}, such that any sheaf $\E$ as above
induces a morphism $f:\,S\to\M$. Similarly for the sub-moduli problem
$\M_L$ for sheaves of fixed determinant $L$ (\ref{det}), with its
moduli space $\M_L\subset\M$. We shall call this standard data:
\begin{description}
\item[$\bullet$] $S$ an affine scheme,
\item[$\bullet$] $\E_0$ on $X\x S$, stable on each fibre of
$p:\,X\x S\to S$ and flat over $S$,
\item[$\bullet$] Chern classes $c_i(\E_0)\in H^{2i}(X)$, a rank $r(\E)
\in H^0(X)$, and a line bundle $L$ on $X$ with $c_1(L)=c_1$,
\item[$\bullet$] The corresponding classifying
morphism $f:\,S\to\M=\M(X,c_i)$,
\item[$\bullet$] An $\OO_S$-module $\II$.
\end{description} \vspace{-1cm}
\begin{equation} \label{std} \vspace{3mm}
\end{equation}
Here we have also included the arbitrary $\OO_S$-module $\II$, along
which we will consider deformations.

\begin{Definition} \emph{\cite{LT}}
Given standard data (\ref{std}) as above, the tangent functor of the
moduli functor assigns to each $\II$ an $\OO_S$-module
$\T_{\E_0}^1(\II)$ such that the set of sheaves on $X\x(S*\II)$
restricting to $\E_0$ on $X\x S$ is isomorphic to
$$
\Gamma_S(\T_{\E_0}^1(\II)).
$$
We also require that given $\E_0$ on $X\x S_2$, a morphism of
schemes $f:\,S_1\to S_2$ and a homomorphism $f^*\II_2\to\II_1$ induce a
canonical homomorphism
$$
f^*\T_{\E_0}^1(\II_2)\to\T^1_{f^*\E_0}(\II_1)
$$
compatible with base-change as in \emph{\cite{LT}}.
\end{Definition}

Denote by $\Ext_p$ the right derived functors of $\Hom_p=p_{*\,}\Hom$,
where $p:\,X\x S\to S$. The functoriality
and compatability with base-change properties will follow from those
of $\Ext_p$ in our case. For $S$ affine, so that $H^1$ of any
$\OO_S$-module vanishes, the Leray spectral sequence shows that
Ext$_{X\x S}^1(\E_0,\E_0\otimes p^*\II)=\Gamma_S
(\Ext_p^1(\E_0,\E_0\otimes p^*\II))$ and similarly, for its trace-free
counterpart,  Ext$_{X\x S}^1(\E_0,\E_0\otimes p^*\II)^{ }_0=\Gamma_S
(\Ext_p^1(\E_0,\E_0\otimes p^*\II))^{ }_0$. Thus we have, by Lemma
\ref{defs} and Theorem \ref{tr},

\begin{Prop} \label{t}
The tangent functor for the above sheaf moduli problem $\M$
assigns to any $\OO_S$-module $\II$ the $\OO_S$-module
$$
\T_{\E_0}^1(\II)=\Ext_p^1(\E_0,\E_0\otimes p^*\II).
$$
For the moduli functor $\M_L$ the same applies to the trace-free
$\OO_S$-module $(\T_{\E_0}^1)^{ }_0(\II)=\Ext_p^1(\E_0,\E_0\otimes
p^*\II)^{ }_0$.
\end{Prop}

\begin{Definition} \label{obdef} \emph{\cite{LT}}
Given standard data (\ref{std}) as above, an obstruction sheaf for the
moduli functor is an $\OO_S$-module $\T_{\E_0}^2$ satisfying the
following conditions.

Let $S\subset Y\subset Y_1$ be as in (\ref{ideals}), and take a
sheaf $\E$ over $X\x Y$ that is flat over $Y$, stable on the
$X$-fibres, and restricts to $\E_0$ on $S$. Then there is an
obstruction class
$$
\mathrm{ob}\,(\E,Y_0,Y_1)\in\Gamma_S(\T_{\E_0}^2\otimes\I)
$$
whose vanishing is necessary and sufficient for there to be an
$\OO_{Y_1}$-flat extension of $\E$ to $\E_1$ over $X\x Y_1$.  This
should be functorial and canonical under base-change (as in
\emph{\cite{LT} 1.2}, but with contravariance replacing covariance).
\end{Definition}

\begin{Theorem} \label{o}
Given any affine scheme $S$ and a sheaf $\E_0$ on $X\x S$ as above
such that $\mathrm{dim\,Ext}^i(\E_0\res{X_s},\E_0\res{X_s})$ is
constant for all $s\in S$ and $i\ge3$, the sheaves
$$
\T_{\E_0}^2=\Ext_p^2(\E_0,\E_0) \qquad\mathrm{and}\qquad
(\T_{\E_0}^2)^{ }_0=\Ext_p^2(\E_0,\E_0)^{ }_0
$$
are obstruction sheaves for both $\M$ and $\M_L$.
\end{Theorem}

\begin{Proof} Theorem \ref{obsclass} shows that
$$
\mathrm{ob}\,(\E,Y_0,Y_1):=\partial(e_\E)\in\mathrm{Ext}^2_{X\x S}
(\E_0,\E_0\otimes p^*\I)
$$
defines an obstruction class, and Theorem \ref{tr} shows that in fact
it lies in the subspace Ext$^2_{X\x S}(\E_0,\E_0\otimes p^*\I)^{ }_0$
since deformations of the line bundle $\det\E$ are unobstructed
(Pic\,$X$ is smooth since $X$ is smooth). By the Leray spectral
sequence $\mathrm{Ext}^2_{X\x S}(\E_0,\E_0\otimes p^*\I)=\Gamma_S
(\Ext_p^2(\E_0,\E_0\otimes\I))$ (and similarly for the trace-free
Exts, for which the rest of the proof is also similar and is omitted).

There is a spectral sequence Tor$_j(\Ext_p^i(\E_0,\E_0),\I)
\Longrightarrow\Ext_p^{i-j}(\E_0,\E_0\otimes p^*\I)$ (by standard nonsense
and the flatness of $\E_0$ so that Tor$_j(\E_0,p^*\I)=0,\ j>0$).  By
base-change (e.g. \cite{BPS}) the hypothesis of the Theorem implies
that $\Ext^i_p(\E_0,\E_0)$ is locally free for $i\ge3$, so its Tors
vanish (as do those of $\Hom_p(\E_0,\E_0)\cong\OO_S$). So the spectral
sequence for $\Ext_p^2(\E_0,\E_0\otimes p^*\I)$ degenerates to
$\Ext_p^2(\E_0,\E_0)\otimes\I$, and $\Ext_p^2(\E_0,\E_0)$ is an
obstruction sheaf as defined in (\ref{obdef}).

The necessary base-change for $\T_{\E_0}^2$ follows from base change
for $\Ext_p$.
\end{Proof}

The condition required to employ the machinery of \cite{LT}, \cite{BF}
to obtain a virtual cycle is then the following, as promised in
(\ref{pre}).

\begin{Definition} \emph{\cite{LT}}
$\M$ has a \emph{perfect} tangent-obstruction complex if there is a
complex $E_1\to E_2$ of locally free sheaves on $\M$ resolving the
tangent obstruction functors (\ref{t}, \ref{o}), in the following
sense. Given standard data (\ref{std}) we require that the cohomologies
of the complex
$$
f^*E_1\otimes\II\to f^*E_2\otimes\II
$$
are $\T_{\E_0}^1(\II)$ and $\T_{\E_0}^2\otimes\II$, in degrees 1
and 2, respectively. Similarly for $\M_L$ using the
$(\T^i_{\E_0})^{ }_0$\,s.
\end{Definition}

\begin{Remarks} The first condition is easily seen to be equivalent to
the exactness of $E_2^*\to E_1^*\to\Omega_\M\to0$, because of course
the tangent functor assigns to data (\ref{std}) the $\OO_S$-module
$$
\T_{\E_0}^1(\II)=\Hom_S(f^*\Omega_\M,\II),
$$
by the standard deformation theory of the morphism $f$. One can also
prove the equality
$\Hom_S(f^*\Omega_\M,\II)=\Ext^1_S(f^*\E,f^*\E\otimes\II)$ directly
using (\ref{t}) to identify $\Omega_\M$ with
$\Ext_p^{n-1}(\E,\E\otimes K_X)$ (recovering the result of \cite{L};
here $\E$ is a local universal bundle on
$X\x\M\rt{p}\M$ as will be described below) and using some
base-change and relative Serre duality \cite{Ha3}.

There is a weaker notion of perfect in (\cite{LT}; comments
following Corollary 3.6) which we will need later -- namely that $E_1$
need only exist locally on $\M$, with $E_2$ still a global vector
bundle on $\M$ surjecting onto $\T^2$.
\end{Remarks}

\begin{Theorem} \label{resol}
Let $X$ be a smooth, polarised, complex projective variety, and fix Chern
classes $c_i\in H^{2i}(X)$ and a line bundle $L$ on $X$ with $c_1(L)=c_1$.
Let $\M$ denote the corresponding moduli space of \emph{stable} sheaves,
and $\M_L$ the subscheme of those with determinant $L$ (\ref{det}).
If the numbers
$$
\mathrm{dim\,Ext}^i(\E,\E),\quad i\ge3,
$$
are constant over $\M\ni\E$ (e.g. if $\mathrm{Ext}_0^i(\E,\E)=0\ \,\forall
\E\in\M,\ \forall i\ge3$), then the tangent-obstruction complex of $\M$
given by $\T^1_{\E_0},\,\T^2_{\E_0}$ (\ref{t}, \ref{o}) is perfect.
Similarly, for rank $r>0$, the tangent-obstruction complex of $\M$
given by $\T^1_{\E_0},\,(\T^2_{\E_0})^{ }_0$ is also perfect, as is
$(\T^1_{\E_0})^{ }_0,\,(\T^2_{\E_0})^{ }_0$ for $\M_L$.
\end{Theorem}

\begin{Remarks}
The required 2-step resolution is given in \cite{LT} for $X$ a
surface, but the method does not generalise to higher dimensions. The
proof below can be seen (see \cite{T1}) as first working out the
tangent-obstruction theory of the Quot scheme (\ref{*}) and then
passing to the quotient $\M$ of the relevant subset of
Quot by the appropriate projective
linear group (\ref{simp}) -- but we shall work purely algebraically.

We work in the universal case of $S$ being (an open subset of) $\M$ and
resolve the sheaves $\Ext^i_p(\E,\E\otimes\II),\ i=1,\,2$, where $\II$
is an $\OO_\M$-module. Pulling back via maps $f:\,S\to\M$ gives all that
we need; the only problem is that the universal sheaf $\E$ on
$X\x\M\rt{p}\M$ may of course not exist. However we shall ignore
this irritation since the usual methods (see e.g. \cite{HL} 10.2)
circumvent it -- $\E$ exists locally on open subsets
$X\x S$ ($S\subset\M$) and choices differ by line bundles pulled
up from $S$. But $\Ext^*_p(\E,\E\otimes\II)$ is invariant under
twisting by such line bundles and so exists uniquely and globally on
$X\x\M$. Thus all sequences in the proof below are really local,
except those involving $\Ext_p$\,s, which patch together globally.
\end{Remarks}

\begin{Proof}
We denote by $\OO(1)$ the pull-back to $X\x\M$ of the polarisation
on $X$, by $P(n)=\chi(\E(n))$ the Hilbert polynomial associated to the
Chern classes $c_i$, and by $\E$ the (local) universal bundle (see the
Remarks above). Choose $n_1\gg0$ such that $\E(n_1)$ is generated by its
fibrewise sections and has no other cohomology, i.e. we have a sequence
\begin{equation} \label{quot}
0\to\K\to p^*\left(p_{*\,}\E(n_1)\right)(-n_1)\to\E\to0,
\end{equation}
for some kernel $\K$, and $R^ip_*(\E(n_1))=0\ \ \forall i\ge1$.
$p_*(\E(n_1))$ is locally free as it has fibres of fixed dimension
$P(n_1)$, and $\K$ is flat over $\M$ because the other two terms are.

Now twist with $n_2$ sufficiently large such that $\K(n_2)$ and
$\E(n_2)$ are generated by fibrewise sections with no $R^ip_*$, and
take cohomology:
\begin{equation} \label{gr}
0\to p_*(\K(n_2))\to p_{*\,}\E(n_1)\otimes\V\to p_*(\E(n_2))\to0.
\end{equation}
$\V$ denotes $p_*(\OO(n_1-n_2))$, and to pull the sheaf
$p_{*\,}\E(n_1)$ through $p_*p^*$ we have used its local freeness.

Define an $\OO_\M$-flat sheaf $\KK$ by
$$
0\to\KK\to p^*p_{*\,}\K(n_2)\to\K(n_2)\to0.
$$
Then applying $\Hom_p(\ \cdot\ ,\E(n_2)\otimes p^*\II)$, for any
$\OO_\M$-module $\II$, yields
\vspace{2mm}
\begin{equation} \label{ex}
\hspace{-4cm} 0\to\Hom_p(\K,\E\otimes p^*\II)\to(p_{*\,}\K(n_2))^*\otimes
p_{*\,}\E(n_2)\otimes\II
\end{equation} $$ \hspace{3cm}
\to\Hom_p(\KK,\E(n_2)\otimes p^*\II)\to\Ext_p^1(\K,\E\otimes p^*\II)\to0,
$$
where the final zero comes from the choice of $n_2\gg0$, and the second
term is produced by the projection formula for the $\OO_\M$-flat sheaves
$p^*p_{*\,}\K(n_2)$ and $\E(n_2)$.

The higher terms in the above long exact sequence (\ref{ex}), with $\II$
trivial, give
$$
0\to\Ext^i_p(\KK,\E(n_2))\to\Ext^{i+1}_p(\K,\E)\to0,\quad i\ge1,
$$
while the long exact $\Hom_p(\ \cdot\ ,\E)$ sequence of (\ref{quot})
yields
$$
0\to\Ext_p^j(\K,\E)\to\Ext_p^{j+1}(\E,\E)\to0, \quad j\ge1,
$$
by the choice of $n_1\gg0$. Thus
\begin{equation} \label{3}
\Ext_p^i(\KK,\E(n_2))\cong\Ext_p^{i+1}(\K,\E)\cong\Ext_p^{i+2}(\E,\E),
\end{equation}
for all $i\ge1$.

The last term is locally free, by base-change and the constancy of
dim\,Exts in the hypothesis. Therefore so is the first term for
$i\ge1$, so its $i=0$ counterpart
$$
E_2:=\Hom_p(\KK,\E(n_2))
$$
is also locally free by base-change \cite{BPS}. Here we
have used the facts that $X$ is smooth and that the sheaves concerned
are flat over $\M$, which also implies that Tor$_i(\E(n_2),p^*\II)=0$
for $i>0$ giving a spectral sequence Tor$_i(\Ext_p^j(\KK,\E(n_2)),\II)
\Longrightarrow\Ext_p^{i-j}(\KK,\E(n_2)\otimes p^*\II)$ as in the proof
of Theorem \ref{o}. By (\ref{3}) this vanishes for $i\ge1,\ j\ge1$ since
$\Ext_p^{j+2}(\E,\E)$ is locally free for $j\ge1$, so the spectral
sequence degenerates to give $\Hom_p(\KK,\E(n_2)\otimes p^*\II)\cong
E_2\otimes\II$. Therefore (\ref{ex}) has become
\begin{equation} \label{*}
0\to\Hom_p(\K,\E\otimes p^*\II)\to E_1'\otimes\II\to E_2\otimes
\II\to\Ext_p^2(\E,\E\otimes p^*\II)\to0,
\end{equation}
with $E_2$ a vector bundle on $\M$; here we have defined
$E_1':=(p_{*\,}\K(n_2))^*\otimes p_{*\,}\E(n_2)$.

Applying $\Hom_p(\ \cdot\ ,\E\otimes p^*\II)$ to (\ref{quot}) yields
$$ \hspace{-3cm}
0\to\Hom_p(\E,\E\otimes p^*\II)\to(p_{*\,}\E(n_1))^*\otimes
p_{*\,}\E(n_1)\otimes\II \vspace{-2mm}
$$ $$ \hspace{3cm} \to\Hom_p(\K,\E\otimes p^*\II)\to
\Ext_p^1(\E,\E\otimes p^*\II)\to0.
$$
The first two terms are $\II$ (by base-change since $\E$ is stable,
and so simple, on the fibres) and $\E\hspace{-1pt}nd_\M(p_{*\,}\E(n_1))
\otimes\II$ with the identity map between them, giving
$$
0\to\E\hspace{-1pt}nd_0(p_{*\,}\E(n_1))\otimes\II\to\Hom_p
(\K,\E\otimes p^*\II)\to\Ext_p^1(\E,\E\otimes p^*\II)\to0.
$$

Fit this into the sequence (\ref{*})
$$
\begin{array}{rccl}
& 0 \\ & \downarrow \\ & \E\hspace{-1pt}nd_0(p_{*\,}\E(n_1))\otimes\II \\ &
\downarrow & \searrow \\ 0 \to\AND \Hom_p(\K,\E\otimes p^*\II) \AND\to\AND
E_1'\otimes\II\to E_2\otimes\II\to\Ext_p^2(\E,\E\otimes p^*\II)\to0 \\ &
\downarrow \\ & \Ext_p^1(\E,\E\otimes p^*\II) \\ & \downarrow \\ & 0
\end{array} \vspace{-14mm}
$$
\begin{equation} \label{simp} \vspace{6mm}
\end{equation}
and divide by the two injections of $\E\hspace{-1pt}nd_0(p_{*\,}\E(n_1))$
to give the required sequence
\begin{equation} \label{atlast}
0\to\Ext^1_p(\E,\E\otimes p^*\II)\to E_1\otimes\II\to
E_2\otimes\II\to\Ext_p^2(\E,\E)\otimes\II\to0,
\end{equation}
where I claim that $E_1$ is locally free. To prove this it is enough
to show the above map of vector bundles $\E\hspace{-1pt}nd_0
(p_{*\,}\E(n_1))\to(p_{*\,}\K(n_2))^*\otimes
p_{*\,}\E(n_2)=E_1'$ is nowhere zero. But this follows easily by
repeating all of the above analysis at a single point of $\M$, instead
of relative to $\M$, and using base-change: the same exact sequences
show that the map on fibres End$_0(H^0(\E(n_1))\to
H^0(\K(n_2))^*\otimes H^0(\E(n_2))$ is an injection.

Finally, take the cokernel of the map $R^1p_*\OO\to E_1$, and/or the
kernel of the map $E_2\to R^2p_*\OO$, in the following diagram
$$
\spreaddiagramrows{-1pc}
\spreaddiagramcolumns{-0.8pc}
\diagram
0 \rto & \Ext^1_p(\E,\E\otimes p^*\II) \rto &
E_1\otimes\II \rto & E_2\otimes\II \rto\drto &
\Ext_p^2(\E,\E\otimes p^*\II) \dto^\tr\rto & 0 \\
& R^1p_*\OO\otimes\II \uto^\id\urto &&& R^2p_*\OO\otimes\II \dto \\
& 0 \uto &&& \,0.
\enddiagram
$$
\vspace{-15mm} \begin{equation} \label{atlast0} \end{equation}
For rank $r>0$ these maps are injective and surjective respectively,
and so give locally free resolutions of the trace-free
tangent-obstruction complexes.
\end{Proof}

\begin{Corollary} \label{virt}
Let $X$ be a smooth projective 3-fold with trivial or anti-effective
canonical bundle, let $\M$ denote the projective moduli space of
semistable (with respect to the projective polarisation) sheaves
of some fixed rank $r$ and Chern classes, and $\M_L$ those sheaves of fixed
determinant $L$. Suppose that all such sheaves are stable
(e.g. if rank and degree are coprime). Then for $r>0$ there is a virtual
cycle $Z_0\subset\M_L$, defined by the tangent-obstruction functors
$(\T^1)^{ }_0,\,(\T^2)^{ }_0$, of dimension the virtual dimension
$$
vd=\sum_{i=0}^3(-1)^{i+1}\mathrm{dim\,Ext}_0^i(\E,\E).
$$
Its class in the Chow group $A_{vd}(X)$ is independent of the resolution
(\ref{resol}). For any $r$ there are similar classes $Z\subset\M$ using the
tangent-obstruction functors $\T^1,\,(\T^2)^{ }_0$ (of dimension
$vd+h^{0,1}(X)$) and $\T^1,\,\T^2$ (dimension $vd+h^{0,1}(X)-h^{0,2}(X)$).
If $\M$ is smooth then the appropriate obstruction sheaf is locally free,
and the virtual cycle is its top Chern class.
\end{Corollary}

\begin{Proof}
Hom$_0(\E,\E)=0$ for any stable sheaf $\E$ in the moduli space, so by Serre
duality Ext$^3_0(\E,\E)\cong$\,Hom$_0(\E,\E\otimes K_X)^*$ and the
assumptions on the canonical bundle $K_X$, Ext$^3_0(\E,\E)$ vanishes
also. Thus by Theorem \ref{resol} we may apply (\cite{LT} 3.7) to give
the required virtual cycle.
\end{Proof}

For $r>0$ we consider the trace-free obstruction class $(\T^2)^{ }_0$,
since $\T^2$ contains a trivial $R^2p_*\OO$ factor whose Chern class vanishes
making the class of the virtual moduli cycle zero. For a Calabi-Yau manifold
the above cycle $Z_0$ has dimension zero and its length will give our
definition of the holomorphic Casson invariant. For this to be a sensible
definition, however, we would like it to be deformation invariant.
So we need to work out the deformation and obstruction theory
of a sheaf $\E$ on $X$ as the complex structure of $X$ varies over an
affine curve. We will do this in the next section; meanwhile since
$\Pee^3$ has no complex deformations we can use the following Corollary
to give an application of the virtual cycle.

\begin{Corollary} \label{hilb}
Hilbert schemes of curves in a 3-fold $X$ with trivial or negative
canonical bundle (such as $\Pee^3$) have a virtual moduli cycle.
\end{Corollary}

\begin{Proof}
We use the particular moduli space of rank 1 sheaves of trivial
determinant that contains the ideal sheaves of the curves. All we need
to show is that any sheaf $\E$ of the same Hilbert polynomial and
determinant is also an ideal sheaf (these are then automatically stable).
But since $\E$, being stable, is torsion free, it is
contained in its double dual, which is rank one, reflexive, and so a
line bundle (\cite{OSS} pp 154--156).  So it equals its own
determinant $\OO$, so $\E\subset\OO$ is a sheaf of ideals.
\end{Proof}

Corollary (\ref{virt}) allows us to define Donaldson-like invariants
for such a 3-fold $X$ by doing intersection theory on moduli spaces of
sheaves using characteristic classes of universal sheaves, in the
usual way. For instance we may now define (deformation invariant,
by the results of the next section)
Gromov-Witten-like invariants (with \emph{integral coefficients})
using Corollary (\ref{hilb}): choose $n$ homology classes $\alpha_i$
in $X$ of total codimension equal to the virtual dimension of a fixed
Hilbert scheme $H$ of curves in $X$, plus two. There is a universal
curve $C$ over $H$ with a \emph{flat} morphism to $H$, and similarly
the projection of the $n$th fibre product to $H$,
$$
C\x_H\ldots\x_HC\to H,
$$
is flat. Thus (\cite{Fu} 1.7) we may pull back the virtual cycle in
$H$ to a cycle $Z$ in the total space and push it forward via the
universal evaluation map
$$
ev:\,C\to X \mathrm{\qquad inducing \qquad}
ev^n:\,C\x_H\ldots\x_HC \to X^n,
$$
to $X^n$. Then the invariant is the integer intersection number
(in the \emph{smooth} variety $X^n$)
\begin{equation} \label{GW}
ev^n_*(Z)\,.\left(\x_{i=1}^n\alpha_i\right).
\end{equation}
Similarly for a Calabi-Yau 3-fold $X$ we may simply count the number
of points in the virtual moduli cycle (i.e. take its length as a
scheme) to get a count of curves in $X$.

These invariants differ from the Gromov-Witten invariants since
Hilbert schemes contain many nasty components representing things
other than curves. An example that will also be relevant later
(explained to me by Jun Li) is given by two disjoint $\Pee^1$\,s in
3-space coming together at a single point in a flat
family. Consideration of the Euler characteristic of the structure
sheaf, or just looking at the equations defining the subscheme, shows
that the limiting curve must have a fat point at the intersection
point (pointing in the direction in which the curves came
together). In a separate flat family this point can break off to give
a $\Pee^1$ and a distinct point in the same Hilbert scheme.

\begin{Remark} A naive approach to creating a virtual fundamental class
would simply be to take the top Chern class of the obstruction
sheaf. This is correct if $\M$ is smooth and so we have an
obstruction \emph{bundle}.  In general Pidstrigatch \cite{Pi} and
Siebert \cite{Si} have shown the correct formula is the natural
generalisation of this given (\ref{sm},\,\ref{pre}), namely the
$vd$-dimensional part of
$$
Z=c\,(E_1-E_0)\cap[c_F(\M)],
$$
where $c_F(\M)$ is Fulton's total Chern class of the scheme $\M$
(\cite{Fu} 4.2.6).
\end{Remark}

\subsection*{Deformation invariance}

Fix a quasi-projective scheme
$\X$ with a flat map to a smooth affine curve $C$, with projective
fibres $\iota:\,X_t\into\X$ over $t\in C$.
Given a stable (in particular, simple)
sheaf $\E$ on $X_0$ we study deformations of its (stable and simple)
pushforward $\iota_*\E$ to $\X$, thus allowing it to move onto other
fibres $X_t$ (it is easy to see, using stability or simplicity,
that this is all it can do: its support must remain over a finite
number of points in $C$ and stability, which is an open property,
prevents it from splitting over more than one fibre; the deformation
theory below will show this for instance).

We will need here and later a technical result, part of whose proof
was worked out with the generous help of Brian Conrad.

\begin{Lemma} \label{exty}
Suppose $\iota:\,D\subset Z$ is a Cartier divisor in a
quasi-projective scheme $Z$, with normal bundle $\nu=\OO_D(D)$. Then
for coherent sheaves $\E$ and $\F$ on $D$ there is a long exact sequence
$$
\to\mathrm{Ext}^i_D(\E,\F)\to\mathrm{Ext}^i_Z(\iota_*\E,\iota_*\F)\to
\mathrm{Ext}^{i-1}_D(\E,\F\otimes\nu)\rt{\delta}\mathrm{Ext}^{i+1}_D
(\E,\F)\to\ .
$$
\end{Lemma}

\begin{Proof}
The sheaf sequence $0\to\OO(-D)\to\OO\to\iota_*\OO_D\to0$ yields
$$
\Hom\,(\OO,\iota_*\F)\to\Hom\,(\OO(-D),\iota_*\F)\to\Ext^1(\iota_*\OO_D,
\iota_*\F)\to0.
$$
The first map, multiplication by the section of $\OO(D)$ defining $D$,
is zero since $\iota_*\F$ is supported on $D$. Thus we see that
$\Ext^1(\iota_*\OO_D,\iota_*\F)\cong\iota_*(\F\otimes\nu)$, so
$\Ext^1(\iota_*E,\iota_*\F)\cong\iota_*(\Hom\,(E,\F)\otimes\nu)$ for
any locally free sheaf $E$ on $D$. Thus, for a locally free sheaf $E$
that is sufficiently negative (so that Ext$^1=H^0(\Ext^1)$ in what
follows) we have
\begin{equation} \label{hom}
\mathrm{Ext}^1(\iota_*E,\iota_*\F)\cong\mathrm{Hom\,}(E,\F\otimes\nu),
\end{equation}
with higher Exts zero.

Now take a (not necessarily finite) locally free resolution of
$E\udot\to\E\to0$ of $\E$ on $D$ with the $E^i\,$s ($i=0,1,\ldots$)
sufficiently negative with respect to $\iota_*\F$ as above (this is
possible since $D$ is projective and we may take
$E^{i+1}=H^0(E^i(N))\otimes\OO_D(-N)$ for large $N$). Take also an
injective resolution $0\to\iota_*\F\to I\udot$ on $Z$
(${\scriptstyle\bullet}=0,1,\ldots$). Then the
cohomology of the complex Hom$\udot(\iota_*\E,I\udot)$ computes
Ext$\udot(\iota_*\E,\iota_*\F)$. But $I\udot$ is a complex of
injectives bounded from below so respects quasi-isomorphisms like
$\iota_*E\udot\to\iota_*\E\to0$ (since $\iota_*$ is exact). Thus the
double complex Hom$\udot(\iota_*E\udot,I\udot)$ also computes
Ext$\udot(\iota_*\E,\iota_*\F)$.

Now, by (\ref{hom}) above, the associated single complex Hom$\udot
(\iota_*E^i,I\udot)$ (for \emph{fixed} $i$, and with differential
$\delta_i^j:\,$Hom$\,(\iota_*E^i,I^j)\to\,$Hom$\,(\iota_*E^i,I^{j+1})$)
has cohomology only in degrees 0 and 1. Thus, truncating all of these
complexes (as $i$ varies) simultaneously by setting terms in degree
${\scriptstyle\bullet}\ge2$ to zero, and replacing the degree one term
by ker$\,\delta_i^1$, we get a quasi-isomorphic complex $B^{i,\bull}$,
sitting in an exact sequence of complexes
$$
0\to\mathrm{ker\,}\delta_i^0\to B^{i,\bull}\to\mathrm{coker\,}
\delta_i^0\,[-1]\to0.
$$
(Here $[-1]$ means shift the complex one place to the right.) So by
(\ref{hom}) this sequence is just
\begin{equation} \label{cxseq}
0\to\mathrm{Hom}\,(E^i,\F)\to B^{i,\bull}\to\mathrm{Hom}\,(E^i,
\F\otimes\nu)\to0.
\end{equation}
The above complexes compute Ext$^i(\E,\F)$, Ext$^i(\iota_*\E,\iota_*\F)$
and Ext$^{i-1}(\E,\F\otimes\nu)$ respectively, so taking the long exact
sequence in cohomology of the exact sequence of total complexes gives
the required result.
\end{Proof}

\begin{Remark} The usual arguments show that the sequence is independent
of choice of resolutions, and that the maps are the natural ones. The
map Ext$^1_D(\E,\F)\to$ Ext$^1_Z(\iota_*\E,\iota_*\F)$ pushes forward
by $\iota_*$ (which is exact) an extension of $\E$ by $\F$ on $D$ to
the corresponding extension of $\iota_*\E$ by $\iota_*\F$ on $Z$.  The
map Ext$^1_Z(\iota_*\E,\iota_*\F)\to\,$Hom$_D(\E,\F\otimes\nu)$ is a
little harder to describe -- taking an extension class on $Z$ and a
(local) section of $\E$ on $D$, we must produce a section of
$\F\otimes\nu$. But pulling back the extension (of $\iota_*\E$ by
$\iota_*\F$) by the section of $\E$ gives a section of
$\Ext^1(\iota_*\OO_D,\iota_*\F)$ which is shown in the proof to be
canonically isomorphic to $\iota_*(\F\otimes\nu)$.
\end{Remark}

So in our usual set up (\ref{ideals}) of $S\subset Y\subset Y_1$,
consider the obstructions and deformations of a sheaf $\F:=\iota_*\E$
on $\X\x Y$, where $\E$ is a stable sheaf on $X_0\x Y$ (flat
over $Y$), and $\iota$ denotes the inclusions $X_0\subset\X$ and
$X_0\x Y\subset\X\x Y$ of the fibre of $\X$ over $0\in C$.
Restricting to $S\subset Y$ gives $\F_0=\iota_*\E_0$.

The sequence (\ref{dagger}) for the deformations of $\F=\iota_*\E$,
and the corresponding sequence for $\E$ in one lower degree, are
linked by the sequence of Lemma \ref{exty} (applied to the sheaves
$\E$ and $\E\otimes\I$ on the Cartier divisor $X_0\x Y\subset\X\x
Y$ and written vertically below) to form the diagram:
\spreaddiagramrows{-0.6pc}
\spreaddiagramcolumns{-0.6pc}
$$
\diagram
& \mathrm{Ext}^1_{X_0\x S}(\E_0,\E\otimes\n) \dto^{\iota_*}
\rto^{\partial_0\ \,} &
\mathrm{Ext}^2_{X_0\x S}(\E_0,\E_0\otimes\I) \dto^{\iota_*} \\
\mathrm{Ext}^1_{\X\x S}(\F_0,\F\otimes\m) \rto\dto &
\mathrm{Ext}^1_{\X\x S} (\F_0,\F\otimes\n) \dto\rto^{\partial\ \,} &
\mathrm{Ext}^2_{\X\x S}(\F_0,\F_0\otimes\I) \dto^\phi \\
\mathrm{Hom}_{X_0\x S}(\E_0,\E\otimes\m) \rto &
\mathrm{Hom}_{X_0\x S}(\E_0,\E\otimes\n) \rto &
\mathrm{Ext}^1_{X_0\x S}(\E_0,\E_0\otimes\I). \enddiagram
$$
Here we have used the fact that the normal bundle to $X_0\x S$ in
$\X\x S$ is of course trivial. Since our extension class $e_\F\in$
Ext$^1(\F_0,\F\otimes\n)$ is in the image of $\iota_*$ in
the above diagram, the obstruction $\partial(e_\F)=\iota_*
\partial_0(e_\E)$ is in the kernel of $\phi$. Extending the right hand
vertical sequence upwards therefore shows the following.

\begin{Theorem} \label{fobs}
The obstruction map $\partial$ of (\ref{obsclass}) takes values in
$\mathrm{coker}\,(\delta)$, with $\delta$ the last map in the
following sequence (\ref{exty}) relating the first order deformations
of $\F_0=\iota_*\E_0$ to those of $\E_0$:
\begin{equation} \label{pf}
\vspace{1mm} \hspace{-4cm}
0\to\mathrm{Ext}^1_{X_0\x S}(\E_0,\E_0\otimes\I)\to
\mathrm{Ext}^1_{\X\x S}(\iota_*\E_0,\iota_*\E_0\otimes\I)
\end{equation}
\vspace{-6mm} $$ \hspace{4cm}
\to\mathrm{Hom}_{X_0\x S}(\E_0,\E_0\otimes\I)\Rt{\delta}
\mathrm{Ext}^2_{X_0\x S}(\E_0,\E_0\otimes\I).
$$
The first map pushes deformations on $X_0\x S$ forward to $\X\x
S$. For $\E_0$ simple (e.g. stable) the penultimate term is just
$H^0(\I)$, and $\delta$ is the obstruction to first order
deformations of $\E_0$ off the fibre $X_0\x S\subset\X\x S$.
\end{Theorem}

We want to repeat this result for trace-free determinants. Of course
there can be an obstruction to deforming the determinant of a sheaf
to a nearby fibre, so we need to assume this vanishes by fixing a
determinant that extends to all of $\X$. So choose
line bundle $L$ on $\X$, and study stable sheaves on fibres
$X_t$ of $\X$ whose determinant \emph{on $X_t$} is $L\res{X_t}$.
We will also now insist that rank$\,\E_0=r>0$.

Then we showed above (\ref{fobs}) that the obstruction to extending
$\F$ is the image in Ext$^2(\F_0,\F_0\otimes\I)$ of the obstruction
$\partial(e_\E)\in$\,Ext$^2(\E_0,\E_0\otimes\I)$ to extending
$\E$ from $X_0\x Y$ to $X_0\x Y_1$. By Theorem (\ref{tr}) this
last obstruction in fact lies in the trace-free part
Ext$^2_0(\E_0,\E_0\otimes\I)$. So to get the trace-free analogue of
the sequence (\ref{fobs}) we want to show that the map
$$
\mathrm{Hom}_{X_0\x S}(\E_0,\E_0\otimes\I)\Rt{\delta}
\mathrm{Ext}^2_{X_0\x S}(\E_0,\E_0\otimes\I)
$$
of (\ref{fobs}) in fact factors through Ext$^2_0$. But this map fits
into the diagram
$$
\diagram
\mathrm{Hom}_{X_0\x S}(\E_0,\E_0\otimes\I) \rto^{\delta}\dto^{\tr} &
\mathrm{Ext}^2_{X_0\x S}(\E_0,\E_0\otimes\I) \dto^{\tr} \\
H^0_{X_0\x S}(\I) \rto^{\delta} & H^2_{X_0\x S}(\I).
\enddiagram
$$
where the bottom row is the corresponding sequence for det$\,\E$,
so the map $\delta$ is the obstruction to extending det$\,\E$, which
vanishes by our assumption that det$\,\E$ is the restriction of
a global line bundle $L$ on $\X$. Now this diagram commutes by the
following observation of Brian Conrad, for which I am very grateful.
Namely, the map Hom$\,(\E,\F)\rt{\delta}$\,Ext$^2(\E,\F)$ of Lemma
\ref{exty}, is covariant in $\F$ and contravariant in $\E$, making
it covariant in $\R\Hom(\E,\F)=\Homd(E\udot,\F)$, where $E\udot$ is
a finite locally free resolution of $\E$. Therefore applying the trace map
to $\R\Hom(\E,\E)$ gives the commutativity of the above diagram.

So we can improve Theorem \ref{fobs} to

\begin{Theorem} \label{fobs0}
Taking $\E_0$ to have rank $r>0$ and determinant $L\res{X_0}$, the
restriction of a global line bundle $L$ on $\X$, Theorem \ref{fobs}
holds with $\mathrm{Ext}^2_{X_0\x S}(\E_0,\E_0\otimes\I)$ replaced by
its trace-free part $\mathrm{Ext}^2_{X_0\x S}(\E_0,\E_0\otimes\I)_0$.
\end{Theorem}

\begin{Corollary} \label{TO}
Use the projections $p:\,X_0\x S\to S$ and $q:\,\X\x S\to S$, and
let $\E$ be a
stable sheaf on $X_0\x Y$ (flat over $Y$) with restriction $\E_0$
to $X_0\x S$ such that $\mathrm{dim\,Ext}^i_0(\E_0\res{X_s},
\E_0\res{X_s})=0$ for all $s\in S$ and $i\ge3$. Then
$$
\T^2_{\F_0}:=\mathrm{coker}\,\{\OO_S\Rt{\delta}
\Ext^2_p(\E_0,\E_0)\}=\mathrm{image}\,\{\Ext^2_p(\E_0,\E_0)\to
\Ext^2_q(\F_0,\F_0)\}
$$
is an obstruction sheaf for $\F=\iota_*\E$. If $r>0$ and 
$\mathrm{det}\,\E_0$ is the restriction to $X_0\x S$ of a global
line bundle $L$ on $\X\x S$ (pulled back from $\X$), then
$$
(\T^2_{\F_0})^{ }_0:=\mathrm{coker}\,\{\OO_S\Rt{\delta}
\Ext^2_p(\E_0,\E_0)^{ }_0\}=\mathrm{image}\,\{\Ext^2_p(\E_0,\E_0)^{ }_0
\to\Ext^2_q(\F_0,\F_0)\}
$$
is also an obstruction sheaf for $\F=\iota_*\E$. Also,
$$
(\T^1_{\F_0})^{ }_0(\II):=\mathrm{coker\,}\{R^1p_*\OO\otimes\II\to
\Ext^1_q(\F_0,\F_0\otimes\I)\}
$$
(where the map is the identity into $\Ext^1_p(\E_0,\E_0\otimes\I)$
followed by the inclusion into $\Ext^1_q(\F_0,\F_0\otimes\I)$)
is a tangent functor for deformations of $\F_0$ with fixed determinant
$L$ on the $X_t$ fibres.

Thus tangent-obstruction functors of the moduli problems for $\E$ on $X_0$
and $\F$ on $\X$ fit into the exact sequences of $\OO_S$-modules
\begin{equation} \label{vary}
0\to\T^1_{\E_0}(\II)\to\T^1_{\F_0}(\II)\to\II\to\T^2_{\E_0}\otimes\II
\to\T^2_{\F_0}\otimes\II\to0,
\end{equation}
and
\begin{equation} \label{vary0}
0\to(\T^1_{\E_0})^{ }_0(\II)\to(\T^1_{\F_0})^{ }_0(\II)\to\II\to
(\T^2_{\E_0})^{ }_0\otimes\II
\to(\T^2_{\F_0})^{ }_0\otimes\II\to0,
\end{equation}
for any $\OO_S$-module $\II$.
\end{Corollary}

\begin{Proof} Use the results of (\ref{o}, \ref{t}) for $\E_0$
and $\F_0=\iota_*\E_0$, using the condition that the higher Ext groups
of $\E_0$ vanish to show the same for $\F_0$ (this follows from the exact
sequence (\ref{exty}) as usual.) Then Theorems \ref{fobs}, \ref{fobs0}
give the required exact sequences.
\end{Proof}

\begin{Theorem} \label{compat}
The tangent-obstruction functors of $\E$ on $X_0$ (\ref{t}, \ref{o})
and of $\F=\iota_*\E$ on $\X$ (\ref{TO}) are compatible in the sense
of Li-Tian (\emph{\cite{LT}} Definition 3.8).
\end{Theorem}

\begin{Proof}
The compatibility of Li-Tian says roughly that for an obstructed sheaf
$\E$ on $X_0$ to extend without obstruction inside $\X$, the obstruction
must cancel the obstruction to extending $\iota_*\E_0$ in the direction
of the base $C$.

The precise statement is that there should be an exact sequence
(\ref{vary}) (or (\ref{vary0}) in the trace-free case) such that the
following holds. Let $S\subset Y\subset Y_1$ be as in (\ref{ideals}),
and take a sheaf $\E$ on $X_0\x Y$,
giving a corresponding class $e\in$\,Ext$^1(\E_0,\E\otimes\n)$.
Suppose that the obstruction class $\partial_0(e)\in$\,Ext$^2(\E_0,\E_0
\otimes\I)$ to extending it to $\X_0\x Y_1$, has vanishing image in
Ext$^2(\F_0,\F_0\otimes\I)$. Thus $\F=\iota_*(\E)$ extends to
a sheaf on $\X\x Y_1$, giving a class $f\in\,$Ext$^1(\F_0,\F\otimes\m)$
and a diagram \vspace{-6mm}
$$
\spreaddiagramcolumns{-2mm}
\spreaddiagramrows{0mm}
\diagram
&&& \hspace{-8mm}e\in\,\mathrm{Ext}^1(\E_0,\E\otimes\n) \dto^{\iota_*} \\
&& \hspace{-8mm}f\in\,\mathrm{Ext}^1(\F_0,
\F\otimes\m) \rto\dto^\phi & \mathrm{Ext}^1(\F_0,\F\otimes\n) \dto^\psi \\
\hspace{24mm}0 \rto & \mathrm{Hom}\,(\E_0,\E_0\otimes\I) \rto\dto^\delta &
\mathrm{Hom}\,(\E_0,\E\otimes\m) \rto & \mathrm{Hom}\,(\E_0,\E\otimes\n) \\
e\in\,\mathrm{Ext}^1(\E_0,\E\otimes\n) \rto^{\,\partial_0} &
\mathrm{Ext}^2(\E_0,\E_0\otimes\I),
\enddiagram
$$
where all horizontal maps come from sequences of the form (\ref{dagger}),
and vertical maps from (\ref{exty}).
Since $\iota_*(e)$ is mapped to zero under $\psi$, $\phi(f)$ is in the
image of a unique class $c\in$\,Hom$\,(\E_0,\E_0\otimes\I)$.
Then Definition 3.8 of \cite{LT} requires that $\delta(c)=-\partial_0(e)$.

In our situation this holds by abstract homological algebra. We first
need to lift the big commutative diagram of Ext groups that the above
diagram is a part of to the level of short exact sequences of complexes.
Given an $\OO_S$-module $J$, we use the complex $B_J\udot$, whose
cohomology is Ext$\udot(\F_0,\F\otimes J)$, from the proof of Lemma
\ref{exty}. This is functorial in $J$, and we denote
the other two complexes in (\ref{cxseq}) by $A_J\udot$ and $C_J\udot$
(computing the cohomology of Ext$\udot(\E_0,\E\otimes J)$ and
Ext$^{\bull-1}(\E_0,\E\otimes\nu\otimes J)$ respectively).

Thus we get a (horizontal) sequence of sequences of the form (\ref{cxseq})
(written vertically) by setting $J$ to be the different modules in the
exact sequence $0\to\I\rt{j}\m\rt{\pi}\n\to0$, yielding the exact commutative
diagram of complexes
$$
\diagram
& 0 \dto & 0 \dto & 0 \dto \\
0 \rto & A\udot_\I \dto^{\iota_\I}\rto^{j_A} & A\udot_\m \dto^{\iota_\m}
\rto^{\pi_A} & A\udot_\n \dto^{\iota_\n} \rto & 0 \\
0 \rto & B\udot_\I \dto^{p_\I}\rto^{j_B} & B\udot_\m \dto^{p_\m}
\rto^{\pi_B} & B\udot_\n \dto^{p_\n} \rto & 0 \\
0 \rto & C\udot_\I \rto^{j_C}\dto & C\udot_\m \rto^{\pi_C}\dto &
C\udot_\n \rto\dto & 0 \\
& 0 & 0 & \ \,0\,.
\enddiagram
$$
We now follow the previous diagram chase around this diagram at the level
of complexes. So we start with
$e\in A^1_\n$ with $\delta e=0$ (we will denote all coboundary
operators by $\delta$) which we want to lift to $f\in B^1_\m$ with
$\delta f=0$.

Lift $e$ to $\beta\in A^1_\m$ (using the fact that $\pi_A$ is onto). This
is then assumed not coclosed; its coboundary is by definition
$\delta(\beta)=j_A(\partial_0(e))$ where $\partial_0(e)\in A^2_\I$
represents the obstruction class $\partial_0(e)\in$\,Ext$^2(\E_0,
\E_0\otimes\I)$ to lifting $e$. Pushing forward $\tilde f:=\iota_\m(\beta)$
to a class in $B^1_\m$, we see that a coclosed class $f\in B^1_\m$ lifting
$\iota_\n(e)\in B^1_\n$ exists if and only if there is an $\alpha\in
B^1_\I$ with $\delta(j_B(\alpha))=-\delta(\tilde f)$; the required $f$ is
then $f=\tilde f+j_B(\alpha)$. We are assuming $f$, and so $\alpha$, exist.

Thus, since $j_B(\delta(\alpha))=-\delta(\tilde f)=-\delta(\iota_\m(\beta))
=-\iota_\m j_A(\partial_0(e))=-j_B\,\iota_\I(\partial_0(e))$, and since $j_B$
is an injection, we have
\begin{equation} \label{cob}
\delta(\alpha)=-\iota_\I(\partial_0(e)).
\end{equation}
In particular then, $p_\I(\alpha)$ is coclosed, defining a class
$c\in$\,Hom$\,(\E_0,\E_0\otimes\I)$ whose coboundary (by which we mean
now the connecting homomorphism $\delta$ in the cohomology sequence of
the left hand column that arises in the previous diagram) $\delta(c)$ is
$-\partial_0(e)$ by (\ref{cob}).

Thus we are left with showing that $c$ is the same $c$ as defined above,
i.e. that $j_C(c)$ is $p_\m(f)$. But $j_C(c)=j_Cp_\I(\alpha)=p_\m j_B
(\alpha)=p_\m(f-\tilde f)=p_\m(f)$ as required, since $\tilde f=
\iota_\m(\beta)$ is in the kernel of $p_\m$.
\end{Proof}

\begin{Corollary} \label{!}
The virtual moduli cycles of Corollary \ref{virt} are deformation
invariant in the following sense. Given a family $\X\to C$ of smooth
projective varieties $X_t,\,t\in C$, consider the family of moduli
spaces $\M_t$ of stable sheaves $\E$ as in (\ref{virt}) on the fibre
$X_t$ (containing the virtual moduli cycle $Z_t$). These form the
moduli space $\M\to C$ of sheaves $(\iota_t)_*\E$. Then under
the same conditions as Corollary \ref{virt}, $\M$ has a virtual
moduli cycle $\mathcal Z$ of dimension $\mathrm{dim}\,Z_t+1$,
and as elements of the Chow group they satisfy $\iota_t^!\mathcal Z
=Z_t$. (Here $\iota_t^!$ is the Gysin homomorphism (\emph{\cite{Fu}} 6.2)
of the inclusion $\iota_t:\,\{t\}\to C$.) There is also the corresponding
result for sheaves of fixed determinant $L_t$, with $L$ a fixed
line bundle on all of $\X$, for rank $r>0$.
\end{Corollary}

\begin{Proof}
We first need to show that the tangent-obstruction functors of $\M\to C$
are perfect in the weaker sense of the Remark preceding Theorem
\ref{resol}. The $E_2$ of that Theorem exists on all of $\M$ and
surjects onto $\T^2_{\E_0}$, and this in turn surjects onto $\T^2_{\F_0}$.
So all we need is a local vector bundle $E_1$ giving a resolution
of $\T^1_{\F_0},\,\T^2_{\F_0}$ over an open set of $\M$.

The map $\OO_\M\to\T^2_{\E_0}$ of (\ref{vary}) locally factors through
$E_2\to\T^2_{\E_0}$ via a local lift $\OO\to E_2$. Combining this lift
with the map $E_1\to E_2$ of Theorem \ref{resol} into a map $E_1\oplus
\OO\to E_2$ gives an exact sequence, for any $\OO_S$-module $\II$,
$$
0\to\T^1_{\F_0}(\II)\to(E_1\oplus\OO)\otimes\II\to E_2\otimes\II
\to\T^2_{\F_0}\otimes\II\to0,
$$
by combining (\ref{atlast}) with (\ref{vary}). Thus we can use \cite{LT}
to produce the virtual moduli cycle $Z$ (and similarly for the trace-free
versions, using (\ref{atlast0}) and (\ref{vary0}) instead).

$\iota_t^!\mathcal Z=Z_t$ follows from (\cite{LT} 3.9), given the
compatibility of Theorem \ref{compat}.
\end{Proof}

\subsection*{The holomorphic Casson invariant}

\begin{Definition} Fix a smooth projective Calabi-Yau 3-fold $X$, and a
rank $r$ and Chern classes $c_i$ such that the moduli space $\M$ of
semistable sheaves with this data (and fixed determinant $L$ if rank$\,>0$)
contains only stable sheaves (for instance if the rank and
degree are coprime). Then we define the holomorphic Casson invariant
$\lambda_{\{c_i\}}(X)$ to be the length as a scheme of the
zero dimensional projective virtual moduli cycle $Z_0\subset\M_L$ of
Corollary \ref{virt}. It is invariant under deformations of $X$ in any
projective family to which $L$ extends (e.g. if $h^{0,2}(X)=0$ this
is immediate).
\end{Definition}

The deformation invariance comes from the relation $\iota_t^!Z=Z_t$ of
(\ref{!}) and the resulting ``conservation of number'' (\cite{Fu}
10.2).  This invariant is clearly similar in nature to the
Gromov-Witten invariants of $X$. In fact we might expect to recover GW
invariants either by considering moduli of ideal sheaves of curves as
in (\ref{hilb}, \ref{GW}), or by relating rank two
bundles to curves via zero sets of
their sections (and vice-versa by the Serre construction). However, we
have already remarked that the first case gives something slightly
different to GW invariants.  As for rank two bundles, under the Serre
construction spheres and tori tend to correspond to unstable bundles,
and for a higher genus curve to correspond to a bundle its tangent
bundle must extend to a line bundle on $X$, cutting down the space of
admissible curves. Since the space of curves has expected dimension
zero anyway we tend to find (for instance in the examples below) that
rank two bundles correspond to non-generic high-dimensional families
of curves. In this case deforming the curve corresponds to deforming
the section, not the bundle, so GW invariants do not arise.

We might also like to count fewer singular sheaves, i.e. to only count
bundles. We would expect from Donaldson theory to have to include
sheaves with codimension two singularities (limits of stable bundles,
where ``bubbling'' occurs), but in the rank two case one might hope to
be able to ignore sheaves with codimension three
singularities. However again the example above of a flat family of
curves producing a distinct point shows that bundles can degenerate to
sheaves with codimension three singularities, and so such sheaves can
lie in the same connected component as bundles.

In all of the examples we consider, however, we will be able to show
there are no such singular sheaves, and in fact GW invariants will
arise in a completely different way in the last section. \\

\noindent \textbf{Examples\ } Simple examples of the invariant are
given by considering ideal sheaves of 1, 2 or 3 points in a Calabi-Yau
3-fold $X$. The moduli space is then $X$, Hilb$^2X$ or Hilb$^3X$, and
so smooth, with the invariant giving the Euler number of the cotangent
bundle (since this is the obstruction bundle), i.e. $-\chi(X),\
\chi($Hilb$^2X)$, and $-\chi($Hilb$^3X)$ respectively. \\

A deeper example is motivated by Donaldson's reinterpretation
\cite{DT} of work of Mukai (\cite{Mu1} 0.9) as an example of the
Tyurin-style Casson invariant of (\ref{tyc}).

Fix a smooth quadric $Q_0$ in $\Pee^5$, in a fixed $\Pee^2$-family of
quadrics spanned by $Q_0,\ Q_1$ and $Q_2$, say. The singular quadrics
in the family lie on the sextic curve
$$
C=\{[\lambda_0;\lambda_1;\lambda_2]\in\Pee^2:\det 
(\lambda_0Q_0+\lambda_1Q_1+\lambda_2Q_2)=0\}\subset\Pee^2
$$
where the quadratic form defining the quadric becomes singular.

For every point of $\Pee^2\take C$ we get two tautological rank 2
bundles $A$ and $B$ over the corresponding quadric (thinking of it as
a Grassmannian of 2-planes in $\C^4$, $A$ and $B$ are defined by the
tautological sequence $0\to A^*\to\C^4\to B\to0$), and so also over
the $K3$ surface
$$
S=Q_0\cap Q_1\cap Q_2.
$$
In fact these $A$ and $B$ bundles give a double cover $\M$ of $\Pee^2$
branched along the sextic curve $C$ (the $A$ and $B$ bundles coincide
on the singular quadrics) as the moduli space of bundles of the same
topological type over $S$. $\M$ is a (complex symplectic) $K3$
surface, notice.

Similarly the Fano $X_1=Q_0\cap Q_1$ lies in the pencil spanned by
$Q_0$ and $Q_1$, a line $\Pee^1$ in our $\Pee^2$-family. The 2-fold
branched cover of this line induced by $\M\to\Pee^2$, i.e. the set of
$A$ and $B$ bundles on the quadrics in this pencil, is the moduli
space for $X_1$. Similarly for $X_2=Q_0\cap Q_2$ and the cover of the
line $\langle X_0,X_2\rangle\subset \Pee^2$.

So we have an example of the Tyurin-Casson invariant (\ref{tyc}), with
two complex Lagrangians (the curves covering the $\Pee^1$s) as the
moduli spaces of bundles on the two Fanos, injecting into the complex
symplectic moduli space of bundles on the common anticanonical divisor
$S$. Their intersection, namely the double cover of the intersection
point $\{Q_0\}$ of the lines in $\Pee^2$, corresponds to the two
stable bundles $A_{Q_0}$ and $B_{Q_0}$ on the singular Calabi-Yau that
is the union of $X_1$ and $X_2$.

Deforming this singular quartic in $Q_0$ to a smooth Calabi-Yau we
would like, then, to prove the following.

\begin{Theorem}
Let $Q_0$ be a smooth quadric in $\Pee^5$, and let $X$ be a smooth
quartic hypersurface in $Q_0$. Then the bundles $A$ and $B$ on $Q_0$
restrict to stable, isolated bundles of the same topological type on
$X$, and they are the only semistable sheaves in the moduli space.
Thus the corresponding holomorphic Casson invariant is 2.
\end{Theorem}

\begin{Proof}
Standard exact sequences and geometry on the Grassmannian
Gr$\,(2,4)\cong Q_0$ show that the bundles are stable and isolated
(\cite{T1} 2.3.1). The more difficult part is to show
that any semistable sheaf $\E$ of the same Chern classes is either $A$
or $B$. We will do this by controlling $\E$'s cohomology by studying
it on hyperplane sections, and then using Riemann-Roch to produce
sections (c.f. \cite{Ha2}).

Let $\dd$ denote the double dual (or ``reflexive hull'') of $\E$. Let
$H=\Pee^4\cap X$ be a smooth hyperplane section, with
$\Pee^4\subset\Pee^5$ sufficiently generic that $\dd\res H$ is the
double dual of $\E\res H$ and so a bundle $F$, say. Then the
Riemann-Roch formula for $F$ is
\begin{equation} \label{RR}
2h^0(F)-h^1(F)=12-c_2(F),
\end{equation}
by Serre duality, $K_H=\OO_H(1)$, and an unpleasant computation. Thus
$F$ has at least 4 sections, as $c_2(F)=c_2(\dd)\,.\,\omega\le
c_2(\E)\,.\,\omega=4$ (recall that passing to double duals lowers $c_2$
because of the exact sequence embedding a sheaf inside its double
dual).

Also, as $H$ is generic, we may assume its only line bundles are the
$\OO(n)$ bundles, by Noether-Lefschetz theory (see e.g. \cite{GH2}).
$\dd$ is slope semistable and so slope semistable on restriction to
the generic hyperplane (see e.g. \cite{Ha2} 3.2). Thus $F(-1)$ cannot
have any sections, so the four sections of $F$ must vanish only on
points, giving us a Koszul resolution (\cite{GH} p 688) of the form
\begin{equation} \label{Kos}
0\to\OO\to F\to\I_c(1)\to0,
\end{equation}
where $\I_c$ is the ideal sheaf of functions vanishing at
$c=c_2(F)\le4$ points. Taking sections and using the Riemann-Roch
formula (\ref{RR}) gives
$$
h^0(\I_c(1))=h^0(F)-1\ge5-c/2.
$$
However $h^0(\OO_H(1))=5$, and the $c$ points impose at least
min$(2,c)$ conditions on the sections of $\OO_H(1)$ as they are the
restriction of the sections of $\OO_{\Pee^4}(1)$ on the $\Pee^4$
hyperplane in $\Pee^5$. Thus $5-c/2\le h^0(\I_c(1))\le5-$\,min\,$(2,
c)$, i.e. $c\ge2$\,min\,$(2,c)$, whose only integral solutions for
$0\le c\le4$ are $c=0$ and $c=4$. We can rule out $c=0$ by stability
(either by the Bogomolov inequality or the fact that (\ref{Kos}) would
split); thus $c=c_2(F)=4$. Therefore $\dd$ has the same second Chern
class as $\E$ on $X$ (since $H^4(X;\C)$ is generated by a class
non-zero on $H$), and can only differ from it in codimension three.

Notice also that as we must have $h^0(\I_4(1))=3$, the points lie
on a web of hyperplanes in the $\Pee^4$ hyperplane, i.e. \emph{on a
line in $\Pee^4$}.

To control the cohomology of $F$ we use (\ref{Kos}) to give, for
$t\ge2$,
$$
\begin{array}{ccccc}
0\to & H^1(F(t-1)) \AND \to H^1(\I_4(t))\to \AND H^2(\OO_H(t-1)) \AND
\to0, \\ & \downarrow{\wr} && \downarrow{\wr} \\ \AND H^1(F(1-t))^*
\AND\AND H^0(\OO_H(2-t))^*
\end{array}
$$
by stability. But the sequence $0\to\I_4(t)\to\OO_H(t)\to\OO_4(t)\to0$
shows that $h^1(\I_4(t))=0\ \ \forall t\ge3$ since we can find a
polynomial of any degree$\,\ge3$ taking any prescribed values at 4
points on a line (in $\Pee^4$). Similarly $h^1(\I_4(2))=1$, and we
get $H^1(F(-n))=0\ \ \forall n\ge1$.

Thus we can now pass up to $X$ using the sequence
$0\to\dd(-1)\to\dd\to \dd\res H\to0$, which is exact because $\dd$ is
torsion free, giving
$$
H^1(\dd(-n-1))\to H^1(\dd(-n))\to H^1(F(-n)).
$$
Since $H^1(\dd(-n))$ vanishes for large $n$ it therefore vanishes for
all $n\ge1$; in particular $H^2(\dd)=0$. This and stability simplify
the Riemann-Roch formula for $\dd$ to
$$
h^0(\dd)-h^1(\dd)=4+c_3(\dd)/2.
$$
The third Chern class of a rank two reflexive sheaf on a smooth 3-fold
is always nonnegative and vanishes if and only if the sheaf is locally
free (\cite{Ha2} 2.6), so $\dd$ has at least 4 sections which do not
vanish on divisors, by stability. Thus we have a presentation
(\cite{Ha2} 4.1)
$$
0\to\OO\to\E\to\I_C(1)\to0,
$$
for some degree four curve $C$. Computing the third Chern class of
such an extension to be zero shows that $\dd$ is locally free with the
same Chern classes as $\E$, i.e. $\E\cong\dd$ is locally free. (The
point here is that the extension is locally free, and not just
reflexive, because the determinant of $\dd$, restricted to $C$, is
isomorphic to the determinant of the normal bundle to $C$, with the
isomorphism set up by the determinant of the section. Thus the section
vanishes transversally along $C$, the above sequence becomes the
Koszul resolution (\ref{Kos}) for this section, and $\E$ is locally
free. In fact the Serre construction (\cite{OSS} p 93)
uniquely constructs a sheaf $\E$ (which is then locally free) from the
above resolution when the
appropriate determinants are equal; it is only when they differ
that reflexive sheaves aride from Hartshorne's generalisation of the
Serre construction \cite{Ha2}.)

But now $h^0(\I_C(1))\ge3$, so $C$ lies in a web of hyperplanes in
$X\subset Q_0\subset\Pee^5$. Thus it lies on a linear $\Pee^2$ plane
$P\subset\Pee^5$. Since $C$ lies in the quadric $Q_0$, the plane $P$
must do too, otherwise the quadric would intersect $P$ in a conic
curve containing the degree four curve $C$, which is impossible.

But the planes in $Q_0$ are precisely the standard planes in
Gr$\,(2,4)$ -- zero sets of sections of $A$ and $B$. $P$ uniquely
defines either the $A$ or the $B$ bundle on $Q_0$ via the Serre
construction (see e.g. \cite{OSS} p 93)
$$
0\to\OO\to A/B\to\I_P(1)\to0,
$$
by the extension data Ext$^1(\I_P(1),\OO_{Q_0})\cong H^0(\OO_P)\ni1$.
This extension restricts, on the quartic $X$, to $1\in H^0(\OO_C)
\cong\mathrm{Ext}^1(\I_C(1),\OO_X)$, defining our bundle $\E$
$$
0\to\OO\to\E\to\I_C(1)\to0
$$
by uniqueness, since $H^0(\OO_C)\cong\C$ ($C$ is a curve in $\Pee^2$,
so is connected). Therefore $\E$ is one of $A$ or $B$ restricted to
$X$.
\end{Proof}

\section{$K3$ fibrations}

We now turn to calculating the invariants on $K3$-fibred Calabi-Yau
manifolds with no reducible or multiple fibres, using the nice
properties of moduli of bundles and sheaves on a $K3$ due to Mukai
(\cite{HL} Chapter 6). These do not quite generalise to the singular
$K3$ fibres. Although these have a (trivial) dualising sheaf as the
fibres are complete intersections (\cite{Ha1} III 7.11) and so the
usual Serre duality holds (\cite{Ha1} III 7.6), the Serre duality we
have been using \cite{Mu1},
\begin{equation} \label{sd}
\mathrm{Ext}^i(\E,\E)^*\cong\mathrm{Ext}^{2-i}(\E,\E),
\end{equation}
only holds for sheaves $\E$ with a \emph{finite} locally free
resolution.  On the singular fibres there will be sheaves with
unbounded locally free resolutions for which the result does not
hold.

We will consider stability on $X$ using a \emph{suitable} polarisation
$\omega$, in the sense of Friedman. This means that the fibres are
small so that semistable sheaves restrict, on the generic fibre, to
semistable sheaves. The basic idea is to add a large multiple of the
fibre divisor $f$ to any fixed polarisation, so that the sign of the
degree $\omega\,.\,\omega\,.\,c_1$ of a possibly destabilizing
subsheaf is the same as the sign of the degree $f.\,\omega\,.\,c_1$ on
a generic fibre. Since such subsheaves that we need to consider form a
bounded family we can do this:

\begin{Prop} \label{suit}
Let $X$ be a surface-fibred projective 3-fold. Choose rank and Chern
classes such that slope semistability implies slope stability on
smooth fibres (e.g. if rank and degree are coprime). Then we may add a
sufficient number of fibre classes to the polarisation such that
sheaves of the given rank and Chern classes are stable if and only if
they are stable on the generic fibre, and there are no strictly
semistable sheaves. Such a polarisation is called suitable.
\end{Prop}

\begin{Proof}
This should really be proved directly, but a cheat goes as follows.
Choose $N\gg0$ such that any slope semistable sheaf $\E$ (of the given
Chern classes) restricts to a slope semistable sheaf on a generic
hyperplane $H$ in the linear system $|\OO(N)|$ (\cite{HL} 7.2.1). Now
there is an $M$ such that $(\omega+Mf)\res H$ is a suitable
polarisation on $H\to C$ by (\cite{HL} 5.3) (here $\omega$ denotes the
K\"ahler form on $X$, and $f$ the class of a fibre). Thus $\E$ is
slope semistable on the generic fibre $H_t$ of $H\to C$, and so
slope semistable on the generic fibre $X_t$ of $X\to C$ (we are using
the easy fact that the (slope/semi) stability of a sheaf on a
hyperplane section implies the same on the whole space).  By
assumption, then, $\E$ is slope stable on $X_t$.

Reversing the argument, by increasing $N$ if necessary, if $\E$ is
slope stable on $X_t$, then by a theorem of Bogomolov (\cite{HL}
7.3.5) $\E$ is stable on $H_t$. Since $(\omega+Mf)\res H$ is a
suitable polarisation, this implies that $\E$ is slope stable on $H$,
and so on $X$.
\end{Proof}

\begin{Definition} \label{data}
Let $X$ be a smooth polarised 3-fold, $K3$-fibred over a smooth curve
$C$, let $c_i\in H^{i,i}_\Z(X_t),\ i=1,2$ be Chern classes in the
cohomology of a generic fibre, and let $r\in\Z_{\ge0}$ be a rank. Then
we say that $(X,r,c_i)$ is \emph{admissible} if and only if
\begin{description}
\item[\ \ $\bullet$\ \ ] $X\to C$ has no reducible or non-reduced fibres,
\item[\ \ $\bullet$\ \ ] the polarisation on $X$ is chosen to be
suitable for $(r,c_i)$ using Proposition \ref{suit},
\item[\ \ $\bullet$\ \ ] $c_1$ is the restriction of the first Chern
class of a global line bundle on $X$, i.e. is in the image of $H^{1,1}_\Z
(X)\to H^{1,1}_\Z(X_t)$,
\item[\ \ $\bullet$\ \ ] gcd\,$(r,c_1\,.\,\omega,\frac12c_1^2-c_2)=1$
(where $\omega$ is the K\"ahler form of the induced polarisation on the
fibre), and
\item[\ \ $\bullet$\ \ ] on any fibre $X_t$, slope semistability of
sheaves with Chern classes $(r,c_i)$ implies slope stability (e.g. if
rank $r$ and degree $c_1\,.\,\omega\res{X_t}$ are coprime).
\end{description}
We then set $d=2rc_2-(r-1)c_1^2-2(r^2-1)$, the dimension $\mathrm{dim\,
Ext}^1(\E,\E)=\sum_i(-1)^{i+1}\,\mathrm{dim\,Ext}^i_0(\E,\E)$ (by Serre
duality and stability of $\E$) of any nonempty moduli space of
stable sheaves $\E$ of Chern classes $(r,c_i)$ on a smooth fibre.
$d$ is even.
\end{Definition}

The penultimate condition ensures that a universal sheaf exists on
the product of any fibre and its moduli space (\cite{HL} 4.6.6),
and also that semistable sheaves on the fibres are in fact
stable (\cite{HL} 4.6.8 -- here is where we use the assumption
that the fibres are reduced and irreducible to ensure the ranks of
possible destabilising subsheaves are integers.) Thus the fibre moduli
spaces are \emph{fine}. But in general this may not be enough to
ensure slope stability, so we require the last condition. Under these
conditions we will often just talk about stability. For some results
we will restrict to fibrations $X\to C$ whose singular fibres have only
ordinary/rational double point (ODP/RDP) singularities, since we then
understand something about
reflexive sheaves on such $K3$ fibres (e.g. \cite{Is}, \cite{La}).

\begin{Definition} \label{data2}
We say the above data (\ref{data}) is \emph{very admissible} if each
singular fibre of $X\to\Pee^1$ contains only a single ODP, and $r$ and
$c_1\,.\,\omega$ are coprime.
\end{Definition}

\begin{Theorem} \label{univ}
Choose an admissible triple $(X,r,c_i)$ as in Definition \ref{data}.
Then there is a projective scheme $\M\to C$ with fibres $\M_t$ that are
the moduli spaces of stable sheaves of Chern classes $(r,c_i)$ on $X_t$.
Suppose that for some smooth fibre $X_t,\ \M_t\ne\emptyset$. Then
$\M\to C$ is surjective with generic fibres smooth of dimension $d$,
and there is a universal sheaf $T$ on $X\x_C\M$ (unique up to twisting
by the pull-back of a line bundle on $C$).

Suppose now that the singular fibres of $X\to C$ have only RDPs.
If $d=2$ and $r\ge2$ then all elements of $\M$ are reflexive
on their supporting fibre (and so locally free on smooth fibres).
If $d=0$ and $r\ge1$, $T$ is locally free, each $\M_t$ is a single
reduced point, $\M\cong C$, and so $T$ is a \emph{bundle on} $X$.
\end{Theorem}

\begin{Proof}
Simpson's construction of a projective proper $\M\to C$ is now
standard \cite{HL}, given the assumption that there are no
semistable sheaves. If we have a stable sheaf $\E$ in the moduli
space of a smooth fibre $X_t$ then its deformations are unobstructed
since by Serre duality Ext$^2_0(\E,\E)\cong\Endo(\E)^*=0$
and det$\,E$ is unobstructed by the assumptions on $c_1$. Thus it may
deformed off $X_t$ (\ref{fobs}), making $\M\to C$ onto over an open
set of $C$, and so onto all of $C$ by properness. Since
Ext$^2_0(\E,\E)=0,\ \M_t$ is smooth of the correct dimension $d$.

Consider sheaves $\E_t\in\M_t$ to be torsion
sheaves on $X$, by pushing them forward to $(\iota_t)_*\E_t$, where
$\iota_t:\,X_t\to X$ is the inclusion. Then $\M\to C$ is part of the
moduli space of sheaves on $X$ of the same Hilbert polynomial (in fact
deformation theory and stability show it an entire component of the
moduli space). There is a universal sheaf on $X\x\M$ as the numerical
conditions of (\cite{HL} 4.6.6) are satisfied by $(\iota_t)_*\E_t$:
they are satisfied by $\E_t$ on $X_t$, as gcd\,$(r,c_1\,.\,\omega,
\frac12c_1^2-c_2)=1$, and they only depend on the Hilbert polynomial,
which is the same for $\E_t$ and $(\iota_t)_*\E_t$.

The universal sheaf on $X\x\M$ is supported on the image of the diagonal
map $X\x_C\M\to X\x\M$, and so defines $T$ on $X\x_C\M$.

The statements about reflexivity of sheaves of rank $r\ge1$ are standard
arguments of Mukai (\cite{HL} 6.1.6, 6.1.9) for smooth fibres -- taking the
double dual of a sheaf does not affect $c_1$ but decreases $c_2$ (look
at the exact sequence of the sheaf injecting its double dual,
which it does if $r\ge1$ since it must then be torsion-free
by stability), which decreases $d$ (\ref{data}) by $2r$ times as much.
Since the double dual is also stable it sits inside a moduli space
of the correct dimension, thus $d(\dd)$ must be greater than or equal to
zero. Thus if $d=0$, or $d=2$ and $r\ge2$, the sheaf must be its own
double dual and so reflexive.

For a sheaf $\E$ on a singular fibre $X_t$ with only RDPs we again take its
double dual $\dd$, and then pull this up to the minimal desingularisation
$\pi:\,\widetilde X_t\to X_t$ of the fibre (another $K3$), and divide by
torsion. This gives a vector bundle $\widetilde\E$ on $\widetilde X_t$
which we show is stable with respect to the (degenerate) polarisation
$\OO(1)$ pulled up from $X_t$. Any subsheaf $\F\into\widetilde\E$ can be
pushed down to a subsheaf of $\pi_*\widetilde\E$, and $\pi_*\widetilde\E
\cong\dd,\ R^1\pi_*\widetilde\E=0$ (see e.g. \cite{Is}). So for $n\gg0$,
$$
\chi(\F(n))=\chi(\pi_*\F(n))-\chi(R^1\pi_*\F)\le\chi(\pi_*\F(n))
\le\chi(\dd(n))=\chi(\pi_*\widetilde\E(n))=\chi(\widetilde\E(n))
$$
by the stability of $\dd$, demonstrating the stability of $\widetilde\E$.
Thus it is also stable for nearby nondegenerate K\"ahler forms.

In particular $\widetilde\E$ is simple, so that its topological invariant
$d(\widetilde\E)$ (\ref{data}) gives the dimension Ext$^1(\widetilde\E,
\widetilde\E)$ of the moduli space it sits in. This must be greater than
or equal to zero, so the previous argument (in the smooth case) goes
through as before to show that $\E$ is locally free for $d(\E)=0,\ r\ge1$
and reflexive for $d(\E)=2,\ r\ge2$, if we can show the inequality
\begin{equation} \label{ddd}
d(\widetilde\E)\le d(\dd),
\end{equation}
with equality if and only if $\dd$ is locally free (notice
that we already know that $d(\dd)\le d(\E)$ with equality if and only if
$\E$ is reflexive).

So we want to compare the topological invariant $d$ on smooth
fibres with $d(\dd)$ on the resolution of the singular fibre $X_t$.
This is more-or-less contained in the work of Langer \cite{La}, and I am
grateful to him for explaining it to me. He defines a second Chern class
for reflexive sheaves such as $\dd$ on singular surfaces such as $X_t$,
which we will denote by $c_2^L(\dd)$. This is not the same as what we
will denote by $c_2$, namely the class which gives the right contribution
to the Riemann-Roch formula (the deformation invariant $c_2(\dd):=
c_2(\E)-l$, where $l$
is the length of the torsion sheaf cokenel of $\E\into\dd$, and $c_2(\E)$
is measured on a nearby smooth fibre).

Langer's definition (building on work of Wahl) is given in terms of the
sheaf $\widetilde\E=(\pi^*\dd/\,$torsion) upstairs, as (\cite{La} Section 3)
$$
c_2^L(\dd):=c_2(\widetilde\E)-\sum_yc_2(\widetilde\E,y),
$$
where $y$ runs over the RDPs of $X_t$, and $c_i(\widetilde\E,y)$ is the
local $i$th Chern class on its resolution (\cite{La} 2.2, 2.3).
He defines (\cite{La} 2.7)
$a_y(\dd)$ to be the local difference between $c_2^L$ and the $c_2$
we use that fits into the Riemann-Roch formula; that is
$$
c_2(\dd)=c_2^L(\dd)-\sum_ya_y(\dd).
$$
The last result of (\cite{La} Section 6) is that for an RDP $y$,
$$
a_y(\dd)=\frac12c^2_1(\widetilde\E,y)-c_2(\widetilde\E,y),
$$
giving
$$
c_2(\widetilde\E)=c_2(\dd)+\frac12\sum_yc^2_1(\widetilde\E,y).
$$
Also, the local first Chern class satisfies
$$
c_1^2(\widetilde\E)-c_1^2(\dd)=\sum_yc^2_1(\widetilde\E,y),
$$
so that putting these two formulae together gives
\begin{eqnarray} \nonumber \hspace{-6mm}
\big[2rc_2-(r-1)c_1^2-2(r^2-1)\big](\widetilde\E)-\big[2rc_2-
(r-1)c_1^2-2(r^2-1)\big](\dd) && \\
=\sum_yc^2_1(\widetilde\E,y). && \label{dim}
\end{eqnarray}
Thus $d(\widetilde\E)-d(\dd)=\sum_yc^2_1(\widetilde\E,y)$ is the sum of
squares of divisors like $-c_1(\widetilde\E,y)$, a positive multiple of
which is effective and supported entirely on the exceptional set. This
is therefore negative, and zero if and only if $c_1(\widetilde\E,y)=0$ 
for all $y$, if and only if $\dd$ is locally free at all singular
points $y$ \cite{Is}.

Finally, we want to show that if $d=0$ then $\M_t$ is a single
point. This is an argument of Mukai for smooth fibres (\cite{HL} 6.1.6)
which generalises to singular fibres with RDPs since we have just shown
that the sheaves are locally free, so the duality Ext$^i(\E,\E)^*\cong
\,$Ext$^{2-i}(\E,\E)$ (\ref{sd}) holds, which is all that is used.
\end{Proof}

Alternatively, in the $d=0$ case, we can give an easier proof of the
existence of $T$: using the same diagonal map $X\x_C\M\into X\x\M$
and Luna's \'etale slice theorem for the quotient map from the
relevant Quot scheme to $\M$, $T$ exists locally over $C$.
It is simple, so is patched together by
nonzero scalars on overlaps of open sets. Thus the obstruction to the
patchings satisfying the cocycle condition lies in $H^2(\OO^*_C)$,
which is zero. We start by analysing the holomorphic Casson invariant
in this $d=0$ case.

\begin{Theorem} \label{0d}
Let $(X,r,c_i)$ be as in Definition \ref{data}, $K3$-fibred over
$C=\Pee^1$, with $d=0,\ r\ge1$. If the fibres have only RDPs then the
bundle $T\to X$ constructed in Theorem \ref{univ} is slope stable,
isolated, and the only semistable sheaf with the same Chern classes.
\end{Theorem}

\begin{Remarks} Though the method of this proof and that of Theorem
\ref{2d} below could be applied to other such ``adiabatic limit''
problems for stable bundles in any dimension, the conditions we
require (that the fibre moduli spaces should be empty for stable
sheaves of the same rank and $c_1$, lower $c_2$, and any $c_i\ \,
i\ge3$) are so stringent that they are probably only effective on
$K3$-- (or, with some modifications to take account of the fundamental
group, $T^4$--) fibred 3-folds, and Fano-surface-fibred
3-folds. Again, allowing for the fundamental group (i.e. fixing
determinants), we can extend the above to any base curve $C$, but we
restrict to $\Pee^1$ for simplicity and because it is the case
relevant to Calabi-Yau 3-folds.

Notice we are not claiming that the moduli spaces are so simple for
\emph{any} classes on the total space satisfying
$2rc_2-(r-1)c_1^2-2(r^2-1)=0$ on the fibres. Even when a stable bundle
with these classes exists on a fibre, we could modify $T$ by giving it
ideal-sheaf singularities contained in fibres, or make it unstable
on some fibres, or add codimension 3 singularities, to produce
stable sheaves in different, more complicated moduli spaces. What is
remarkable about $T$'s moduli space is that no such singular
phenomena can occur.
\end{Remarks}

\begin{Corollary} Any such $K3$-fibred 3-fold admits isolated bundles
if the generic fibre does (with $c_1$ a class coming from the total
space), and with respect to some polarisation these bundles are stable
and the only point in their moduli space of sheaves. In particular,
the relevant holomorphic Casson invariants of $K3$-fibred Calabi-Yau
3-folds are one. The same is true for the dual Chern classes
$(-1)^ic_i\in H^{i,i}(X_t)$.
\end{Corollary}

\begin{Proof}
The point here is to show that the Casson invariant is defined, i.e.
that semistability implies stability. But this follows from the choice
of polarisation (\ref{suit}), implying that slope semistability
implies slope stability, which is stronger.

As for deformations, in a polarised family of Calabi-Yau manifolds the
fibration structure survives (the obstructions to the survival of a
generic fibre $X_t$ in a deformation lie in $H^1(\nu_{X_t|X})\cong
H^1(\OO_{K3})=0$), and the suitability of the polarisation is
preserved (the volume of the fibres is unchanged, for instance).

For the dual Chern classes we take the dual bundle on a smooth fibre and
repeat the construction, giving the unique bundle $T^*$ in the moduli
space.
\end{Proof}

\begin{Proof}\emph{of Theorem\ }
Denoting the projection map by $\pi:\,X\to\Pee^1$, I claim that
$\Ext_\pi^1(T,T)$ is zero by standard base-change arguments. Namely:
pick a finite, very negative, locally free resolution $E\udot\to T$,
so that Ext$^i(E^j\res{X_t},T\res{X_t})=0$ for $i>0$ and
$t\in\Pee^1$. Then $\Ext^*_\pi(T,T)$ is the cohomology of the complex
$F\udot=p_*\,\Hom(E\udot,T)$. $T$ is flat over $\Pee^1$ so we can
restrict to any fibre $X_t$ to give a resolution of $T\res{X_t}$; thus
$F\udot\res{X_t}=$\ Hom$_{X_t}(E\udot\res{X_t},T\res{X_t})$ computes
Ext$_{X_t}^*(T\res{X_t},T\res{X_t})$. Since the higher Exts vanish by
construction, and since $T$ is flat over $\Pee^1$, the dimension of
each Hom$_{X_t}(E\udot\res{X_t},T\res{X_t})$ is constant in $t$. Thus
$F\udot$ is a complex of \emph{locally frees}. By stability, its
zeroth cohomology on each fibre is canonically
End$\,(T\res{X_t})=\C\,.\,$id, so we get a nowhere vanishing map
$\OO_{\Pee^1}\Rt{\sim}$\
ker\,$(F^0\rt{d^0}F^1)$, and the image
$d^0(F^0/\OO_{\Pee^1})\subset F^1$ is locally free.  Therefore
$\Ext_\pi^1(T,T)$ is the kernel of the map $\widehat{d^1}$ that $d^1$
induces on the locally free cokernel of $d^0$. Similarly
Ext$^1(T\res{X_t},T\res{X_t})$ is the kernel of
$\widehat{d^1}\res{X_t}$ on coker\,$d^0\res{X_t}$. But this is zero
for $X_t$ a smooth fibre, so that $\widehat{d^1}$ is generically
an injection of vector bundles. Thus its kernel, as a map of sheaves,
is zero, and $\Ext_\pi^1(T,T)$ vanishes.

So the Leray spectral sequence for $\pi_*\E\hspace{-1pt}nd_0(T)$ and
$\Ext_\pi^i(T,T)$ shows that $T$ is isolated. (Alternatively, the
proof below extended to $X\x\Spec\C[\epsilon]/(\epsilon^2)$
shows there is a unique bundle on this thickened space.) Stability
follows from fibrewise stability and the choice of polarisation.

What we want to prove is that any stable sheaf $\E$ of the same Chern
classes is isomorphic to $T$. The basic idea is that, firstly, taking
double duals decreases $c_2$ and the fibrewise moduli spaces of lower
$c_2$ are empty, and secondly, that if sheaves have fibres on which
they are unstable, this only ever increases $c_2$. Playing these two
phenomena off against each other ensures we have no non-reflexive
sheaves, and no unstable fibres.

We first replace the stable sheaf $\E$ by its double dual $\dd$, which
on the generic fibre is the double dual of the restriction of $\E$ to
that fibre, and hence slope stable and locally free there.

We have a sequence
$$
0\to\E\to\dd\to\F\to0,
$$
where $\F$ is a sheaf supported on a subvariety $Z$ of codimension two
or higher. Therefore
\begin{equation} \label{ineq}
c_2(T).\,\omega=c_2(\E).\,\omega=c_2(\dd).\,\omega+\int_Z\mathrm{rk}_Z
(\F)\omega\ge c_2(\dd).\,\omega,
\end{equation}
with a similar inequality $c_2(\dd\res{X_t})\le c_2(T\res{X_t})$ on a
generic fibre $X_t$. But the dimension $d$ of the moduli space on a
smooth fibre decreases by $2r\ge2$ for every decrease in $c_2$, and must
remain non-negative on the generic fibre where $\dd\res{X_t}$ is stable.
Since $d=0$ we see that $\dd\res{X_t}$ is reflexive and in the same
moduli space as $T\res{X_t}$, and so it \emph{is} $T\res{X_t}$, recalling
that the moduli space is a single point.

Therefore, using stability, on the generic fibre
Hom$\,(T\res{X_t},\dd\res{X_t})$ is a copy of $\C$, so
$$
\pi_*\,\Hom\,(T,\dd)
$$
is a rank one torsion free sheaf on $\Pee^1$, i.e. $\OO(-n)$ for some
$n$. This gives the exact sequence
\begin{equation} \label{star}
0\to T(-n)\Rt{\phi}\dd\to Q\to0,
\end{equation}
where $T(-n)$ denotes the twist of $T$ with the pullback of $\OO(-n)$
from $\Pee^1$, and $Q$ is a sheaf supported on some finite number $d$
of possibly singular fibres $\{X_{t_i}\}_{i=1}^d$.  (Some of these
fibres might be infinitely close; $d$ is the total number counted with
multiplicities: the length of the scheme in $\Pee^1$
over which the union of the fibres sit.) Write $Q=
\oplus_{i=1}^d\iota_*Q_i$, where $Q_i$ is a sheaf supported on
$X_{t_i}$.

Since $\dd$ is torsion free we have the sequence
$$
0\to\dd(-1)\to\dd\to\dd\res{X_{t_i}}\to0, \vspace{-3mm}
$$
inducing
$$
0\to\mathrm{Hom}(T(-n+1),\dd)\to\mathrm{Hom}(T(-n),\dd)\to\mathrm{Hom}
(T(-n),\dd\res{X_{t_i}}).
$$
Therefore if $\phi$ in (\ref{star}) is zero on $X_{t_i}$ it comes from
an element of the first group in the above sequence. Thus, by reducing
$n$ if necessary, we may assume that $\phi\res{X_{t_i}}\ne0$, so that
the rank of $Q_i$ is at most $r-1$ on its support.

Since $\dd$ is reflexive on a smooth variety, it has homological
dimension at most one, i.e. $\Ext^i(\dd,\F)=0\ \ \forall i\ge2$, for
all coherent sheaves $\F$. Also $T(-n)$ is locally free, so from
(\ref{star}) we see that $\Ext^i(Q,\F)=0\ \ \forall i\ge2$, and $Q$
has homological dimension at most one. Therefore it is supported
in exactly codimension one; i.e. the rank of each $Q_i$ is at least
1 (since there are no reducible fibres). That is $1\le r_i\le r-1$,
where $r_i=\mathrm{\,rk\,}Q_i$.

Now the slope stability of $T\res{X_{t_i}}$, and the sequence
$$
T\res{X_{t_i}}\to\dd\res{X_{t_i}}\to Q_i\to0,
$$
imply that the slope of $Q_i$ is less than that of $T\res{X_{t_i}}$.
That is,
\begin{equation} \label{slope}
c_1(Q_i)\,.\,\omega\res{X_{t_i}}<\frac{r_i}r
c_1(T)\,.\,\omega_1\,.\,\omega,
\end{equation}
where $\omega_1$ is the pullback via $\pi$ of the standard K\"ahler
form on $\Pee^1$.

By the Grothendieck-Riemann-Roch theorem (or more elementary
considerations), $ch(\iota_*Q_i)=\iota_*ch(Q_i)$, where $\iota:\,
X_{t_i}\to X$ is the inclusion of a fibre (in particular it has a
trivial normal bundle), and $ch$ is the Chern character. Therefore
$$
c_1(\iota_*Q_i)=r_i\,\omega_1 \mathrm{\qquad and \qquad}
c_2(\iota_*Q_i)=-\iota_*c_1(Q_i),
$$
since $c_1(\iota_*Q_i)^2=0$. (The pushforward on cohomology $\iota_*$
increases the (complex) degree by one in cohomology; it is Poincar\'e
dual to the inclusion on homology.)

Combining this with (\ref{star}), and denoting the total Chern class
by $c$, we have
$$
c(\dd)=c(T(-n))\prod_{i=1}^dc(\iota_*Q_i)= \hspace{8cm} \vspace{-4mm}
$$ $$
\hspace{1cm}
\bigg(1+c_1(T)-rn\omega_1+c_2(T)-(r-1)n\omega_1c_1(T)\bigg)
\,.\,\left(1+\sum_{i=1}^dr_i\omega_1-\sum_{i=1}^d\iota_*c_1(Q_i)\right)
$$
up to second degree in cohomology. Since $\E$ and $T$ have the same
Chern classes, this gives
$$
c_2(\dd)-c_2(\E)=\left(\sum_{i=1}^d(r_i)-rn\right)\omega_1+\left(\sum_
{i=1}^d(r_i)-(r-1)n\right)\!\omega_1c_1(T)-\sum_{i=1}^d\iota_*c_1(Q_i).
$$
The degree one piece thus gives
$$
\sum_{i=1}^dr_i=rn,
$$
which with (\ref{slope}) yields
\begin{equation} \label{slop}
\sum_{i=1}^dc_1(Q_i)\,.\,\omega\res{X_{t_i}}<n\,c_1(T)\,.\,\omega_1
\,.\,\omega,
\end{equation}
unless the number of fibres $d$ is zero, in which case both sides
vanish. Taking the cup product of the second order piece with
$\omega$, however, gives
$$
n\,c_1(T)\,.\,\omega_1\,.\,\omega-\sum_{i=1}^dc_1(Q_i)\,.\,\omega\res
{X_{t_i}}=(c_2(\dd)-c_2(\E)).\,\omega,$$ which, by (\ref{ineq}), is
nonpositive. Therefore (\ref{slop}) cannot hold, and so $d=0=n$,
$Q=0$, and (\ref{star}) becomes
$$
0\to T\to\dd\to0.
$$
Thus the Chern classes of $\E$ and $\dd$ are the same (they are both
equal to that of $T$) and $\E\cong\dd\cong T$.
\end{Proof} \vspace{-5mm}

\subsubsection*{Examples}

\noindent $\ \bullet\ \,$ Consider the Calabi-Yau 3-fold that is a
smooth (2,2,3) divisor in a product of projective spaces,
$$
X_{2,2,3}\subset\Pee^1\x\Pee^1\x\Pee^2.
$$
The projection $\pi_1$ to the first $\Pee^1$ exhibits $X$ as a $K3$
fibration, with fibre a (2,3) divisor in $\Pee^1\x\Pee^2$.  In
turn this fibre is a double cover of $\Pee^2$ branched over a sextic
-- the zero locus of the discriminant of the quadratic on $\Pee^1$
(which has coefficients that are cubics on $\Pee^2$) defining the $K3$
fibre. It is well known that the pullback of the tangent bundle of
$\Pee^2$ is an isolated slope stable bundle on any such $K3$ with
respect to the pullback of the polarisation on $\Pee^2$ (the proof in
\cite{DK} 9.1.8 works even for the singular covers).  So taking the
polarisation $\pi_1^*\OO(N)\otimes\pi_3^*\OO(1),\ \, N\gg0$, in the
obvious notation, we see that
$$
\pi_3^*\ T\,\Pee^2\to X
$$
is slope stable and unique in its moduli space of sheaves, giving a
holomorphic Casson invariant of one. The same is true of the pullback
of $T\Pee^2$ to a more general double cover of $\Pee^1\x \Pee^2$
branched over a smooth (4,6)-divisor. \\

\noindent $\ \bullet\ \,$ On $K3\x T^2$, fixing determinants to be
trivial in the $T^2$ direction, the method of the proof of the theorem
shows that, given an isolated slope stable bundle on $K3$ with
gcd\,$(r,c_1\,.\,\omega,\frac 12c_1^2-c_2)=1$, we get a holomorphic
Casson invariant of $r^2$ by pulling the bundle back from $K3$ and
twisting by line bundles on $T^2$ whose $r$th power is
trivial. (Notice that here we are taking a polarisation with the $T^2$
fibres large, the opposite of \cite{FMW}.) It would be especially
interesting to study the higher dimensional moduli spaces. The parts
of the moduli space arising as pullbacks from $K3$ are smooth and the
corresponding invariant would be the Euler characteristic of the
moduli space, to within a sign. Generating functions of these give
modular forms, by the work of Vafa and Witten \cite{VW}, so it would
be interesting to compute the corrections given by other parts of the
moduli space. For instance, the next section will show there are no
corrections for the 2-dimensional moduli spaces and we get a Casson
invariant
$$
\lambda(K3\x T^2)=r^2\chi(\M_{K3}).
$$

\noindent $\ \bullet\ \,$ It would be fashionable, and by now almost
expected, to find modular forms arising from the Casson invariants on
more general Calabi-Yau manifolds. An easy, cheating, way to produce
them is to take the moduli spaces of stable bundles on Fano surfaces in
\cite{VW}, and push these forward to give torsion sheaves on any
Calabi-Yau 3-fold that is modeled locally on the total space of the
canonical bundle over the Fano. Deformation theory (\ref{exty}) shows
that all sheaves in
the same component of the moduli space are of the same form, so the
holomorphic Casson invariant is the Euler characteristic of the moduli
space if it is smooth, and some appropriate modification of it if not.
(In fact this gives a way of rigorously defining the Euler
characteristics of moduli spaces in \cite{VW} so long as there are no
semistable sheaves.) Thus generating functions of the invariants (for
Chern classes on the Calabi-Yau that are pushforwards from the Fano
surface) do indeed give modular forms in the cases studied in \cite{VW}. \\

\noindent $\ \bullet\ \,$ We can also consider the two examples of
Mark Gross \cite{Gr}, \cite{Ru}.  Let $E_1$ be the trivial rank four
bundle $\OO^{\oplus4}$ over $\Pee^1$,
$E_2=\OO(-1)\oplus\OO\oplus\OO\oplus\OO(1)$, and $P_i=\Pee(E_i)$.  Let
the $X_i$ be anticanonical divisors in $P_i$ (Gross shows they can be
chosen smooth). They are clearly both fibred $\pi:\,X_i\to\Pee^1$ by
$K3\,$s that are quartics in $\Pee^3$. Their second cohomologies are
generated by $t:=c_1(\OO_{\Pee(E_i)}(1))$ and the fibre class
$f:=\pi^*\omega_1$. We use a polarisation $t+Nf,\ N\gg0$.

The $X_i$ are in fact diffeomorphic via a diffeomorphism taking $t$ and
$f$ on $X_1$ to $t$ and $f$ on $X_2$. They are not, however, deformation
equivalent as projective or even symplectic manifolds (by work of
P.\,M.\,H. Wilson in the first case, and Gross and Ruan in the second
\cite{Ru}). We might hope to be able to prove the projective statement
(or even both if Tian's work mentioned in Section 3 can be developed
to give symplectic invariants) using the holomorphic Casson
invariant. The most obvious bundle to take is, in the first case,
$$
\pi^*\,T\,\Pee^3\to X_1.
$$

\begin{Prop} \label{S}
The restriction of $T\,\Pee^3$ to the generic $X_1\to\Pee^1$ is, on
restriction to each fibre, slope stable (with respect to the
Fubini-Study K\"ahler form), isolated, and unique in its moduli space.
The holomorphic Casson invariant of any $X_1$, for sheaves of total
Chern class $1+4t+6t^2+4t^3$, is 1.
\end{Prop}

\begin{Proof}
While it is well known that $T\,\Pee^3$ is stable on the generic
quartic this is not enough for us -- it is a priori possible that the
unique stable bundle on each quartic is generically $T\,\Pee^3$ but
something else on a closed subset of quartics. We need to show
that it is stable on each quartic $S$ in a generic one-parameter
family $X_1$. The other statements follow
by the usual methods if we can show that $T\,\Pee^3$ is stable, since it
is locally free and so satisfies Serre duality even on the singular
quartics.

So fix a (possibly singular) quartic $S\subset\Pee^3$. Suppose firstly
that $T\,\Pee^3\res S$ is (slope) destabilised by a rank two
subsheaf. Then the quotient map gives us a sequence
$$
T\,\Pee^3(-1)\res S\to\Ll\to0
$$
on $S$, with $\Ll$ a torsion free (without loss of generality) rank
one sheaf. $\Ll$ is generated by its sections since $T\,\Pee^3(-1)$
is, and of degree less than or equal to 1. If $\Ll$ were trivial then
the dual of the above sequence would give $\Omega_{\Pee^3}(1)\res S$ a
section which the restriction sequence from $\Pee^3$ shows it does
not have. Therefore $\Ll$ must have degree 1 and at least 2 sections
to be generated by them.

Sections of $\Ll$ vanish on a degree one curve in $S\subset\Pee^3$,
which must therefore be a line. Taking two sections of $\Ll$ we get
two homologous distinct $\Pee^1$\,s in $S$ which is a contradiction
if $S$ is smooth since they must have self-intersection $-2$
by adjunction.

As Adrian Langer pointed out to me this makes no sense if $S$ is
singular; nevertheless the argument goes through as there are no
1-parameter families of $\Pee^1$\,s in $S$ \emph{for $S$ a fibre of
a generic $X_1$} (for instance we can therefore assume that $S$ has
only ODP singularities).

Since we can assume $S$ irreducible and reduced it is now sufficient
to consider $T\,\Pee^3(-1)$ being destabilised by a rank one subsheaf.
Dualising, we get a sequence
$$
\Omega_{\Pee^3}(2)\res S\to\Ll\to0,
$$
for some other torsion free (without loss of generality) rank one
sheaf $\Ll$. Again this shows that $\Ll$ is generated by its sections,
and is of degree $\le2$ this time. But degree 2 curves in $\Pee^3$ are
curves of arithmetic genus 0, and a similar argument applies.

Therefore, by Theorem (\ref{0d}), $\pi^*\,T\,\Pee^3\to X_1$ is the
unique stable sheaf in its moduli space with total Chern class
$1+4t+6t^2+4t^3$. Thus we have a holomorphic Casson invariant of one,
for the generic $X_1$, and \emph{therefore for all} $X_1$ by
deformation invariance of the holomorphic Casson invariant.
\end{Proof}

Unfortunately this does not distinguish the $X_i$s since for the
corresponding
(under the diffeomorphism) Chern class $1+4t+6t^2+4t^3$ on $X_2$ the
invariant is also one, given by the relative tangent bundle down the
fibres of $\Pee(E_2)$:
$$
T_\pi\,\Pee(E_2)\res{X_2}\to X_2.
$$
So to distinguish the $X_i$\,s we have to consider two dimensional
moduli spaces on the $K3$ fibres.

\subsection*{The fibre dimension two case}

To study bundles with $d=2$ dimensional moduli spaces on the fibres we
follow the same method as before, considering sections of the
fibration of fibrewise moduli spaces -- the ``Mukai-dual'' 3-fold
$\M\to\Pee^1$ constructed by applying Mukai duality (replacing a $K3$
by a 2-dimensional moduli space of sheaves on it -- another $K3$)
fibrewise as in Theorem \ref{univ}. We will show that $\M$ is Calabi-Yau.

In this two dimensional situation a new feature can occur, namely the
sections can change homology class. However, just as introducing
codimension two singularities and unstable fibres drives $c_2$ higher,
we shall find that changing the section to one of higher degree does
the same (where degree is measured appropriately). Thus for sections
of the smallest possible degree this cannot happen in a fixed moduli
space.

Relating the deformation theories of the bundles and the sections will
show the Casson invariant of $X$ is equal to the appropriate
Gromov-Witten invariant of $\M$. Thus we
see Gromov-Witten invariants arising, but not in the way we might have
expected (with rank 2 bundles corresponding to curves via zero sets of
their sections), and in fact for arbitrary rank.

Throughout this section we will work with $K3$ fibrations whose singular
fibres have only single ODPs. Many of the results are true for arbitrary
reduced irreducible fibres (and in fact all of them should be,
but I cannot prove it). We start with some technicalities improving
on Theorem \ref{univ}.

\begin{Theorem} \label{ODPs}
Fix $(X\!\to\!\Pee^1,r,c_i)$ as in Definition \ref{data2}, with $d=2$
and $r\ge2$. Then the moduli space $\M\!\to\!\Pee^1$ is
smooth with fibres $\M_t$ having only ODPs as singularities. The
singular points represent the reflexive non locally free sheaves on
$X_t$; all other points correspond to vector bundles.
\end{Theorem}

\begin{Proof}
These results were conjectured in an earlier version of this paper based
on comparing monodromy in $\M\to\Pee^1$ around singular fibres with the
corresponding monodromy in $X\to\Pee^1$ (using Mukai's natural
isomorphisms between the cohomologies of the fibres $X_t$ and $\M_t$ as
in ([HL] 6.1.14); see Theorem \ref{muk} below).
Since then, however, a remarkable proof of the smoothness of $\M$
using Fourier-Mukai transforms has appeared \cite{BrM} (this applies
directly in our case, even without assuming RDP fibres: the assumptions
on Chern classes (\ref{data}) make the fibrewise moduli spaces of
``Fourier-Mukai type'' in the terminology of \cite{BrM} -- the Mukai
vector of the sheaves is primitive of square zero). The remaining
results can be deduced from standard theory of reflexive sheaves on
ODPs, and I would like to thank Akira Ishii for explaining this theory
to me.

From Theorem \ref{univ} all sheaves in $\M$ are reflexive, and locally
free if supported on a smooth fibre. So consider a reflexive, non locally
free sheaf $\E$ on a singular fibre $X_t$. As in the proof of Theorem
\ref{univ}, consider $\widetilde\E=\pi^*\E/\,$torsion on the minimal
resolution $\pi:\,\widetilde X_t\to X_t$. We proved that $\widetilde\E$
was simple and had $d(\widetilde\E)<d(\E)=2$, so that it sits in a moduli
space of dimension zero:
$$
\mathrm{Ext}^1(\widetilde\E,\widetilde\E)=0.
$$
In fact from equation (\ref{dim}) we see that the local first Chern class
satisfies $c_1^2(\widetilde\E,y)=-2$, where $y$ is the ODP (for instance,
dimension considerations force it to be between $-1$ and $-3$, and it
must be even since this is the only ODP, and the intersection form on
$K3$ is even).

Thus standard theory of reflexive modules (see e.g. \cite{Is}
Example 3.2) implies that $\widetilde\E$, restricted to the
exceptional $-2$-sphere $Z\subset\widetilde X_t$, is
$$
\OO_{\Pee^1}^{r-2}\oplus\OO_{\Pee^1}(1)\oplus\OO_{\Pee^1}(1).
$$
The local deformations of the corresponding reflexive module over the ODP
are then given by (\cite{Is} Theorem 4.3 (2))
pushing down bundles whose restriction to $Z$ is any deformation of
\begin{equation} \label{or}
\OO_{\Pee^1}^{r-2}\oplus\OO_{\Pee^1}(1)\oplus\OO_{\Pee^1}(-1)
\end{equation}
(or $\OO_{\Pee^1}^{r-2}\oplus\OO_{\Pee^1}(-1)\oplus\OO_{\Pee^1}(-1)$,
but this gives the same local reflexive module). Pushdowns of deformations
of the above bundle give the 2 dimensional local moduli space
of deformations of the sheaf. This has an ODP (the original reflexive
sheaf) away from which the sheaves are all locally free
(\cite{Is} Theorem 4.6), corresponding to the deformation
$\OO_{\Pee^1}^{\oplus r}$ of (\ref{or}) on the resolution.

Now $H^2(\E\hspace{-1pt}nd\,\E)\cong\,$Hom\,$(\E\hspace{-1pt}nd\,\E,
K_{X_t})^*$ by Serre duality, and this is$\ \C$ ($K_{X_t}$ is trivial and
$\E\hspace{-1pt}nd(\E)^*\cong\E\hspace{-1pt}nd\,\E$ since it is reflexive
-- it is the pushforward of the restriction of itself to the smooth
locus of a normal surface). (Thanks to Akira Ishii for this argument.)

The local-to-global Ext spectral sequence yields
$$
0\to H^1(\E\hspace{-1pt}nd\,\E)\to\mathrm{Ext}^1(\E,\E)\to
H^0(\Ext^1(\E,\E))\to H^2(\E\hspace{-1pt}nd\,\E),
$$
where the last term is all trace $H^2(\OO_{X_t})\cong\C$ that survives in the
spectral sequence for Ext$^2(\E,\E)$. Thus the penultimate map is onto, and
global deformations of $\E$ map \emph{onto} local deformations
$\Ext^1(\E,\E)$ of the module over the ODP. But both have dimension 3
(global deformations of $\E$ are the deformations of a singular point
in a surface $\M_t$ in a smooth 3-fold $\M$, so have Zariski tangent
space of dimension 3), so we get a local
isomorphism from the global deformations to the local model described
above. In particular we see that
there are locally free deformations of $\E$, the moduli space is
2 dimensional, and the non locally free sheaf is isolated in its
moduli space and represents an ODP.
\end{Proof}

We need to fix an appropriate polarisation on $\M\to\Pee^1$. In
fact, all that will concern us will be its restriction to the fibres,
so we consider a polarised $K3$ surface $S$ with at worst an ODP
singularity $y\in S$, and $\M_S$ a moduli space of semistable
sheaves of fixed determinant such that rank and
degree are coprime. Let $H$ be a generic smooth hyperplane section
of $S$ missing its singularity $y$,
$p:\,S\x\M_S\to\M_S$ and $q:\,S\x\M_S\to S$ be the projections, and
let $\E$ be a (local) universal sheaf. Then there is an obvious class
$[\E\hspace{-1pt}nd_0(\E\res H)]$ in the $K$-group $K^0(S\x\M_S)$ which
is just $\E\hspace{-1pt}nd_0\E$ if $\E$ is locally free, and more
generally the kernel of
the trace map on the restriction to $H$ of the tensor product of a
finite locally free resolution of $\E\res{S\backslash\{y\}}$ with
the dual complex. This is uniquely defined even though $\E$ is only
defined locally up to tensoring with a line bundle pulled back from $\M$.

Likewise the push down
$$
p_{\,!}:\,K^0(S\x\M_S)\to K^0(\M_S)
$$
(\cite{HL} 2.1.11) is $p_{\,!}=p_*-R^1p_*+R^2p_*$, or $p_*$ on a finite
resolution by sheaves with no higher cohomology (this always exists
\cite{Ha1} III 2.7). Then the result we need is the following.

\begin{Lemma} \label{above}
In the above set-up and notation, the determinant line bundle
$$
(\det p_{\,!}\,\E\hspace{-1pt}nd_0\E\res H)^*
$$
is ample on $\M_S$.
\end{Lemma}

\begin{Proof}
Fix a point $x\in S\backslash\{x\}$ and let
$\res{\{x\}\x\M_S}$ denote $\otimes\,q^*\OO_x$ in the K-group
(resolve by locally frees away from $y$ and restrict to $x$).
Let $\chi=\chi(\E\res H)=\chi(\E)-\chi(\E(-1))$.  Then it is a result
of Jun Li (\cite{Li}, \cite{HL} Section 8.2) that
$$
\det p_{\,!}\,(\chi\E\res{\{x\}\x\M_S}-r\E\res H)
$$
is ample on the moduli space of slope stable sheaves (which for us is
all of $\M_S$) \emph{for $S$ smooth}. But the proof goes through for
surfaces with isolated singularities for rank and degree coprime:
for a fixed \emph{bundle} $\E$ the
restriction to the generic curve in $|\OO(a)|,\ a\gg0$ is semistable
(Flenner's theorem (\cite{HL} 7.1) applies to normal varieties) and so
stable if $a$ is coprime to the rank, and the generic curve misses the
singularities and is smooth.  Now proceed as in (\cite{HL} 8.2) to
deduce ampleness of Jun Li's line bundle on the locus of locally free
sheaves. As the non locally free sheaves are isolated in the moduli
space (\ref{ODPs}) this is sufficient.

Thus it is sufficient to show that the two line bundles
have the same first Chern class.

But by Grothendieck-Riemann-Roch we have
\begin{eqnarray} \nonumber
c_1\,p_{\,!\,}(\chi\E\res{\{x\}\x\M_S}-r\E\res H) \AND=\AND
p_*\!\left[ch\,(\chi\E\res{\{x\}\x\M_S}-r\E\res H)\,\Td(H)
\right]_2 \\ \nonumber \AND=\AND p_*\big(\!-\!r\,ch_2(E\res
H)\big)+c_1(\E\res{\{x\}\x
\M_S})\left(\chi-\frac12r\,(2-2g)\right),
\end{eqnarray}
where $g$ is the genus of $H\subset S$. Therefore
$\chi-\frac12r\,(2-2g)$ is the degree of $\E\res H$,
i.e. $p_*\,c_1(\E)$, so we obtain
\begin{eqnarray} \nonumber
c_1\,p_{\,!\,}(\chi\E\res{\{x\}\x\M_S} \AND-\AND r\E\res H)=
p_*\big(\!-\!r\,ch_2(\E\res H)+\frac12c_1(\E\res H)^2\big) \\
\nonumber \AND=\AND p_*\big(r\,c_2(\E\res H)-\frac12(r-1)\,c_1(\E\res
H)^2\big) =-\frac12\,p_*\,ch_2(\E\hspace{-1pt}nd_0\E\res H).
\end{eqnarray}
But since $ch_1(\E\hspace{-1pt}nd_0\E\res H)=0$,
Grothendieck-Riemann-Roch equates the last term with
$\frac12c_1\big((p_{\,!\,}(\E\hspace{-1pt}nd_0\E\res H))^*\big)$.
\end{Proof}

So consider again the relative moduli space $\M\to\Pee^1$, with a
(local) universal sheaf $T$. Let $p$ be the projection $p:\,X\xp
\M\to\M$, and pick a smooth hyperplane $H$ of $X$, flat over
$\Pee^1$. Then the first Chern class $\Omega$ of the line bundle
$(\det p_{\,!\,}(\E\hspace{-1pt}nd_0 T\res H))^*$ is, by the above
lemma, ample on the fibres $\M_t$.

So denoting by $\omega_1$ the pull-back to $\M$ of the Fubini-Study
form on $\Pee^1$, $N\omega_1+\Omega$ is ample on $\M$ for $N$
sufficiently large. $N$ will be unimportant for us (we are concerned
with sections of $\M\to\Pee^1$ which all have fixed degree measured
against $\omega_1$), so we fix one such $N$ and measure the degree of
sections against the above polarisation.

In the following, for any section $t:\,\Pee^1\to\M$ we will denote by
$\tilde t$ the induced section $\tilde t:\,X\to X\xp\M$.

\begin{Theorem} \label{2d}
Let $(X,r,c_i)$ satisfy Definition \ref{data2}. Denote by $T$ a
universal sheaf on $X\xp\M$ as produced by (\ref{univ}), and choose a
section $\sigma$ of $\M\to\Pee^1$ of strictly minimal degree with respect
to the polarisation $N\omega_1+\Omega$ described above. Then all stable
sheaves with the same Chern classes as $\tilde\sigma^*\,T$ are
pull-backs of$\ T$ by sections in the same homology class, and these are
locally free.
\end{Theorem}

\begin{Remarks}
The theorem should be true in greater generality, without the assumptions
of Definition \ref{data2}. For fibrations with worse singular fibres
one can extend the proof, getting slightly weaker results. In fact one
can show, using \cite{BrM} to study the obvious Fourier-Mukai
transform on a singular fibre $X_t\times\M_t$, that smooth points
of the moduli space are in one-one correspondence with locally free
sheaves. Understanding when the singularities of $\M_t$ are isolated
is more difficult, and stops one proving the ampleness (\ref{above}) of
Jun Li's line bundle. But in examples this is obvious, so the results
still apply.

As mentioned earlier the theorem also applies to 3-folds fibred by
surfaces with negative canonical bundle, and arbitrary base curves if
we take account of the fundamental group. The $K3$-fibred case is most
relevant to us, however.
\end{Remarks}

\begin{Proof}
We want to study an arbitrary stable sheaf $\E$ with the same Chern
classes as $\tilde\sigma^*\,T$. We follow the proof of Theorem
\ref{0d}, adapting the argument where necessary.

As in (\ref{0d}) the restriction of the double dual $\dd$ to the
generic fibre must lie in the same moduli space as $T$, giving us a
rational (and so regular) section $s:\,\Pee^1\to\M$. Again, as before,
this gives us a sequence
$$
0\to \tilde s^*\,T(-n)\Rt{\phi}\dd\to Q\to0,
$$
for some $n$, with $Q=\oplus_{i=1}^d\iota_*Q_i$ supported on a finite
number of fibres $\{X_{t_i}\}_{i=1}^d$, rank\,($Q_i)=r_i,\ 1\le r_i\le r$,
and $\phi\res{X_{t_i}}\ne0$.

Computing Chern classes as in (\ref{0d}) gives, for $c_1$,
\begin{equation} \label{om1}
c_1(\tilde\sigma^*\,T)-c_1(\tilde s^*\,T)=\left(\sum_{i=1}^d(r_i)-rn
\right)\omega_1,
\end{equation}
which, with the stability inequality (\ref{slope}), yields
\begin{equation} \label{mess}
\sum_{i=1}^dc_1(Q_i)\,.\,\omega\res{X_{t_i}}\le\left[n\,\omega_1+
\frac1r\big(c_1(\tilde\sigma^*\,T)-c_1(\tilde s^*\,T)\big)\right].\,
c_1(\tilde s^*\,T)\,.\,\omega,
\end{equation}
with the inequality strict if $d\ne0$. For $c_2$ we obtain
$$
c_2(\dd)=-\sum_{i=1}^d\iota_*c_1(Q_i)+\left(\sum_{i=1}^d(r_i)-(r-1)n\right)
\omega_1\,.\,c_1(\tilde s^*\,T)\,+\,c_2(\tilde s^*\,T).
$$
Proceeding as before, with (\ref{mess}) this yields \vspace{2mm}
$$
\hspace{-92mm} 0\,\ge\,(c_2(\dd)-c_2(\E))\,.\,\omega\,\ge
$$$$
\hspace{22mm}
\left(1-\frac1r\right)\big[c_1(\tilde\sigma^*\,T)-c_1(\tilde
s^*\,T)\big].\,c_1(\tilde s^*\,T)\,.\,\omega\,+\,\big[c_2(\tilde
s^*\,T)- c_2(\tilde\sigma^*\,T)\big]\,\omega,
$$
with a strict inequality for $d\ne0$.

Now $\big[c_1(\tilde\sigma^*\,T)-c_1(\tilde s^*\,T)\big]^2=0$ from
(\ref{om1}) (this also shows that the inequality is independent of
twisting $T$ by a line bundle on $\Pee^1$, as it should be), so we
may rewrite the first term in the above inequality to arrive at
$$
0\,\ge\,\big(\tilde
s^*-\tilde\sigma^*\big)\,\omega\,.\left[c_2(T)-\frac12
\left(1-\frac1r\right)c_1(T)^2\right].
$$
The idea now is that the quantity in the square brackets is the
discriminant (\cite{HL} 3.4) of $T$ (divided by $2r$), which, by
stability, should be positive in some sense. In fact, letting $H$ be a
smooth hyperplane section dual to $\omega$, flat over $\Pee^1$ as
before, we have
$$
0\,\ge\int_X(\tilde\sigma^*-\tilde s^*)\ \omega\,.\ ch_2\,
(\E\hspace{-1pt}nd_0T)=\int_H (\tilde\sigma^*-\tilde s^*)\ ch\,
(\E\hspace{-1pt}nd_0 T\res H)\,\Td(T_p).
$$
Here $p$ denotes the projection $p:\,H\xp\M\to\M$, with relative
tangent bundle the class $T_p$ in K-theory, and the other terms in
$(ch\,.\,$Td) do not contribute as $ch_1=0$ and $ch_0\,\Td(T_p)$ is
the same under $\tilde s^*$ and $\tilde\sigma^*$.  Also $\res H$ is
meant in the sense of multiplication by $\OO_{H\x _{\Pee^1\,}\M}$
in K-theory, as before.

Thus, by Grothendieck-Riemann-Roch, we have
$$
0\,\ge\int_{\Pee^1}(\sigma^*-s^*)\,p_*\big(ch(\E\hspace{-1pt}nd_0
T\res H)\,\Td(T_p)\big)=\int_{\Pee^1}(\sigma^*-s^*)\ ch\,(p_{\,!}\,
\E\hspace{-1pt}nd_0 T\res H),
$$
where $p_{\,!}=p_*-R^1p_*$. Thus, finally,
$$
0\,\ge\int_{\Pee^1}(s^*-\sigma^*)\
c_1\big((\det p_{\,!}\,\E\hspace{-1pt}nd_0 T\res H)^*\big).
$$
Letting $\Omega=c_1\big((\det p_{\,!}\,\E\hspace{-1pt}nd_0 T\res
H)^*\big)$, which is ample on the fibres, $N\omega_1+\Omega$ is a
polarisation for $N\gg0$. $s$ and $\sigma$ are both sections of
$\M\to\Pee^1$ and so have the same degree with respect to $\omega_1$,
and so we arrive at
$$
0\,\ge\int_{\Pee^1}(s^*-\sigma^*)\,\Omega=\int_{\Pee^1}
(s^*-\sigma^*)\ (N\omega_1+\Omega).
$$
Thus if, as in the conditions of the theorem, $\sigma$ is of strictly
minimal degree then we must have $s$ and $\sigma$ in the same homology
class, and, as the inequality is now non-strict, $d=0=Q$. By
(\ref{om1}) $n=0$, and so $\tilde s^*\,T\cong\dd\cong\E$.

Finally, sections $s:\,\Pee^1\to\M$ miss the singular points of fibres
where the derivative of the projection to $\Pee^1$ vanishes. Thus
$s$ takes its image in the smooth points of $\M$, which correspond
to locally free sheaves on the fibres by (\ref{ODPs}).
Thus $\tilde s^*\,T$ is locally free.
\end{Proof}

\begin{Corollary} \label{GW2}
Under the conditions of Theorem \ref{2d},
the holomorphic Casson invariant of $X$ equals the algebraic
Gromov-Witten invariant $GW_\M([\sigma])$ of $\M$ (as defined in
\emph{\cite{LT}, \cite{Be}}) counting curves in the same homology
class as $\sigma$.
\end{Corollary}

\begin{Proof}
Let $\N$ denote the moduli space of curves in the homology
class (which is isomorphic to the space of stable maps as the curves are
sections -- they can have no bubbles in fibres by the assumptions on degree
-- and so regularly embedded). Use $s:\,\N\x\Pee^1\into\N\x\M$ to
denote the universal curve, with $\tilde s:\,\N\x X\into\N\x\M\xp X$
the induced map. Let the pull-back of the universal sheaf on
$X\xp\M$ to $\N\x X\xp\M$ be denoted by $\E$. Then the sheaf
$\tilde s^*\E$ on $\N\x X$ exhibits $\N$ as the relevant moduli space
of sheaves on $X$ as produced by the above Theorem, depicted by the diagram:
$$
\begin{array}{ccc}
\tilde s^*\E && \E \\ \downarrow && \downarrow \\
\N\x X\!\! & \stackrel{\tilde s\,}\into\!\! & \N\x\M\xp X \\
\,\ \ \downarrow p && \ \ \ \downarrow\pi \\
\N\x\Pee^1\!\! & \stackrel{s\,}\into\!\! & \N\x\M \\
\,\ \ \downarrow q && \,\ \ \downarrow\rho \\ \N && \ \N\x\Pee^1.
\end{array}
$$
We must show that the tangent-obstruction functor of the sheaves
is the same as that of the curves in the Gromov-Witten theory of \cite{LT}.
Since the curves are sections they miss the singularities of the
projection $\M\to\Pee^1$ and so, by the above discussion, lie
in the open set of $\M$ corresponding to locally free sheaves on fibres,
which therefore satisfy the Serre duality (\ref{sd}) on their
Ext$_0$\,s. Thus the only non-zero such Ext$_0$ is Ext$^1_{X_t}(\F_t,\F_t)_0$
(for $\F_t$ a sheaf on the fibre $X_t$), which is the fibrewise tangent
space to $\M\to\Pee^1$ at $\F_t$, i.e. the normal bundle to any section
through $\F_t$.

For smoothly embedded rational curves the deformation theory of \cite{LT}
reduces to $R^iq_*\nu_s\otimes\II,\ i=0,\,1$, where $\nu_s$ is the normal
bundle to $s$ and $\II$ is an arbitrary $\OO_\N$-module. This is then
$$
R^iq_*(s^*T_\rho)=R^iq_*(s^*\Ext_\pi^1(\E,\E\otimes\pi^*\II)_0)=
R^iq_*(\Ext^1_p(\tilde s^*\E,\tilde s^*\E\otimes p^*\II)),
$$
where for the last equality we have base-changed around the square in
the above diagram using
the flatness of $\E$ and $(\Ext^1_\pi)_0$ (since the other
$(\Ext^i_\pi)_0$\,s vanish). By the Leray spectral sequence this yields
$\Ext^{i+1}_{q\circ p}(\tilde s^*\E,\tilde s^*\E\otimes(q\circ p)^*\II)_0,
\ i=0,1$, which gives the deformation theory of the sheaf moduli
problem, as required.
\end{Proof}

\begin{Remark}
Degenerating the base $\Pee^1$ to a curve with one node degenerates $\M$
(which we show below is Calabi-Yau) to a normal crossings space of two
``Fanos'' (their anticanonical bundles are effective) joined across a
common anticanonical divisor $S$ -- the fibre of $\M$ over the
node. Then the Tyurin-Casson invariant picture sketched in (\ref{tyc})
becomes, in terms of the sections of $\M\to\Pee^1$, a formal picture
mentioned in \cite{T1} (where it was also motivated by relating bundles
to curves, but in that case via zero sets of rank 2 bundles). Namely, to
count curves in such a singular Calabi-Yau, one should count those in
each Fano component which meet in the ``boundary'' $S$. Analogously
to the Tyurin picture we find \cite{T1} that the image of the map
taking a curve in one Fano component to its intersection with $S$,
which generically lies in the complex symplectic space Hilb$^n(S)$, is
a complex Lagrangian. (Here $n$ is the intersection number of the curve
with $S$.) Intersecting the two Lagrangians should give
something like the GW invariant of a smoothing of the Calabi-Yau.  It
would be nice to rigorise this and the Tyurin-Casson invariant (of
course one is a special case of the other, by considering the case of
ideal sheaves of curves).
\end{Remark}

We next show that $\M$ is in fact a Calabi-Yau 3-fold. I have since
discovered this was conjectured to be true many years ago by Mark Gross.
The result has now also been given a new proof in \cite{BrM}, using
Fourier-Mukai transforms, without the need for assumptions on the
singular fibres of $X$. Notice that although (\ref{2d}) did not hold
for rank $r\le1$, this result does, provided that $\M$ is smooth
and sheaves supported on singular fibres without locally free resolutions
are isolated. Thus, in the case that the fibres $X_t$ are themselves
elliptically fibred, we may take as $\M_t$ the moduli of torsion sheaves
that are rank one, degree zero sheaves supported on elliptic fibres (so
that $\M_t$ is a compactified Jacobian of $X_t\to\Pee^1$). In a
hyperk\"ahler rotated complex structure (and for some appropriate choice
of ``B-field'') this should be the mirror $K3$ to $X_t$ (\cite{SYZ} 4.1)
and $\M$ is a obtained from $X$ by fibrewise mirror symmetry and
hyperk\"ahler rotation. However $\M$ is \emph{not} the mirror of
$X$ (its Hodge numbers are not flipped\, for instance; see below),
though it is ``T-dual'' to $X$.

\begin{Theorem} \label{CY}
Let $\M$ be a relative moduli space of slope stable sheaves on $X$
with data (\ref{data}), $d=2,\ r\ge2$, and single ODP singularities
in the fibres of $X\to\Pee^1$. Then $\M$ is a Calabi-Yau 3-fold.
\end{Theorem}

\begin{Proof}
As before consider $\M$ to be the moduli space of torsion sheaves
$(\iota_t)_*\E$, where $\E$ is a sheaf on $X_t\subset X$.  Let $T$
denote the universal sheaf on $X\xp\M$, and let $\iota:\,X\xp\M\to
X\x\M$ be the inclusion and $p:\,X\x\M\to\M$ the projection. Then, as
$\M$ is smooth (\ref{ODPs}), it has a tangent sheaf given by
$\Ext_p^1(\iota_*T,\iota_*T)$.

The holomorphic 3-form is the pairing on the tangent bundle given by
the cup product pairings
$$
\Ext_p^1(\iota_*T,\iota_*T)\otimes\Ext_p^1(\iota_*T,\iota_*T)\otimes
\Ext_p^1(\iota_*T,\iota_*T)\to\Ext_p^3(\iota_*T,\iota_*T)\cong\OO_\M,
$$
in the spirit of Mukai's symplectic structure for moduli spaces on
$K3$ surfaces \cite{Mu1}. Notice its dependence on the holomorphic
3-form on $X$ is through the final isomorphism. To show that this pairing
is non-degenerate everywhere and that $\M$ is Calabi-Yau it is enough
to show it does not vanish on any divisor in $\M$. But by (\ref{ODPs})
it is enough to show this at points of $\M$ corresponding to locally
free sheaves since their complement is isolated; these then satisfy
Serre duality on their Exts. As then all of the sheaves in the above
pairing are locally free and we can localise at a point to get, at the
level of tangent spaces, the cup product pairing
$$
\mathrm{Ext}^1(\F,\F)\otimes\mathrm{Ext}^1(\F,\F)\otimes\mathrm{Ext}^1(\F,\F)
\to\mathrm{Ext}^3(\F,\F)\cong\C
$$
at the point $\F=(\iota_t)_*\E\in\M$, where $\E$ is a sheaf on the
fibre $X_t$.

So now to show the pairing is non-degenerate it is enough to show, by
Serre duality, that the pairing
$$
\mathrm{Ext}^1(\F,\F)\otimes\mathrm{Ext}^1(\F,\F)\to\mathrm{Ext}^2(\F,\F)
$$
is onto. By Lemma \ref{exty} above the inclusion $X_t\subset X$
induces a long exact sequence
$$
\ldots\to\mathrm{Ext}^i_{X_t}(\E,\E)\to\mathrm{Ext}^i_X(\F,\F)\to
T_{t\,}\Pee^1\otimes\mathrm{Ext}^{i-1}_{X_t}(\E,\E)\to\ldots,
$$
which yields, for $X_t$ a $K3$ and $\E$ stable and satisfying Serre
duality,
$$
0\to\mathrm{Ext}^1_{X_t}(\E,\E)\to\mathrm{Ext}^1_X(\F,\F)\to
T_{t\,}\Pee^1\to0,
$$
and
$$
0\to\mathrm{Ext}^2_{X_t}(\E,\E)\to\mathrm{Ext}^2_X(\F,\F)\to
T_{t\,}\Pee^1 \otimes\mathrm{Ext}^1_{X_t}(\E,\E)\to0.
$$
The pairing respects these sequences in the obvious way. Firstly the
cup product of the image of two elements of Ext$^1_{X_t}(\E,\E)
\subset\ $Ext$^1_X(\F,\F)$ is just (the image of) the cup product on
$X_t$, with image in Ext$^2_{X_t}(\E,\E)\subset\ $Ext$^2_X(\F,\F)$,
and this is non-degenerate. Secondly pairing an element $x$ of
Ext$^1_X(\F,\F)$ with (the image of) an element $y$ of
Ext$^1_{X_t}(\E,\E)$ and projecting Ext$^2_X(\F,\F)\to
T_{t\,}\Pee^1\,\otimes\,$Ext$^1_{X_t}(\E,\E)$ gives the same as
projecting $x$ to $T_{t\,}\Pee^1$ and tensoring with $y$. Again this
is onto, so we are done.
\end{Proof}

We now study the cohomology of such an $\M$,
putting Mukai's isomorphism of Hodge structures
(\cite{HL} 6.1.14) between $X_t$ and $\M_t$ into a family $X\xp\M$. 
We avoid problems with the singularities of fibres of $X$ by
working on $X\x\M$.

Fix $X$ and a rank and Chern classes satisfying the conditions of
Theorem \ref{2d}, pick a universal sheaf $\E$ on $X\x\M$, and
let $\pi$ and $p$ be the projections to $X$ and $\M$ respectively.

Then define maps
$$
f:\,H^*(X;\C)\to H^*(\M;\C), \qquad f':\,H^*(\M;\C)\to H^*(X;\C),
\vspace{-2mm}
$$
by
$$
f(c)=p_*\!\left(ch\ve(\E)\sqrt{\mathrm{Td}(X\x\M)}\,.\,\pi^*(c)\right)\!,
\ \
f'(c)=\pi_*\!\left(\!-ch(\E)\sqrt{\mathrm{Td}(X\x\M)}\,.\,p^*(c)\right)\!,
$$
where $ch\ve=\sum(-1)^ich_i$.

\begin{Theorem} \label{muk}
$f\circ f'$ is the identity.
\end{Theorem}

\begin{Proof}
We simply mimic the exposition of Mukai's proof for $K3$ surfaces in
(\cite{HL} 6.1.13). Label two copies of $\M$ by $\M_1$ and $\M_2$, and
denote by $\E_i$ the pull back of the universal sheaf on $X\x\M_i$ to
$X\x\M_1\x\M_2$.  We must transfer a class $c$ from $\M_1$ to $X$ via
$f'$, then back to $\M_2$ via $f$. Pulling everything back to
$X\x\M_1\x\M_2$ via the commutativity of various push-down pull-back
squares, and pushing down $X$ before $\M_1$ (as in \cite{HL} 6.1.13),
we can reduce to the following diagram
$$ \begin{array}{c} X\x\M_1\x\M_2 \\ \ \downarrow \tilde p
\vspace{1mm} \\ \M_1\x\M_2 \\ q_2\swarrow \searrow q_1 \vspace{1mm} \\
\M_2 \quad \M_1\,, \end{array}
$$
arriving at
$$
ff'(c)=q_{2\,*\,}\tilde p_*\!\left[-ch\ve(\E_1)ch(\E_2)\sqrt{\Td
(\M_1\x\M_2)}\ \Td(X)\,.\,q_1^*(c)\right],
$$
where we have supressed some obvious pull-back maps for clarity.  By
multiplicativity of the Chern character this yields
\begin{equation} \label{pq}
q_{2\,*}\left[\tilde
p_*\big(\!-ch(\R\Hom(\E_1,\E_2)\,\Td(X)\big)\sqrt{
\Td(\M_1\x\M_2)}\,.\,q_1^*(c)\right].
\end{equation}

Now $\tilde p_*(ch(\R\Hom(\E_1,\E_2)\Td(X))=ch(\R\Hom_{\tilde p}(\E_1,
\E_2))$, by the Grothendieck-Riemann-Roch theorem. But we know that
\begin{equation} \label{van}
\mathrm{Ext}^i_X(\F_1,\F_2)=0, \quad \forall i,
\end{equation}
for sheaves $\F_1\ne\F_2$ in $\M$ \emph{unless both are non locally free
on their support} (where again we are considering $\M$
to be a moduli space of torsion sheaves supported on fibres $X_t$).
This is because if either is locally free then the Serre duality
(\ref{sd}) holds on the fibres and the usual arguments go through
(\cite{HL} 6.1.8, 6.1.10) to give the vanishing of all Ext$^i_{X_t}$s
on the fibre; thus by (\ref{exty}) all Ext$^i_X$s vanish too.

So the Ext$^i$s $(i\le3)$ vanish outside a subvariety of codimension
three (the union of the diagonal $X\x\Delta\subset X\x\M_1\x\M_2$ and
the codimension six
(\ref{ODPs}) locus of pairs of non locally free sheaves) for all
$i\le3$, from which it follows (\cite{Mu2} 2.26) that
$$
\Ext_{\tilde p}^i(\E_1,\E_2)\equiv0, \quad i<3,
$$
and $\Ext_{\tilde p}^3(\E_1,\E_2)$ is concentrated on $\Delta$. On
$\Delta$ we have Ext$^3(\F,\F)\cong H^3(\OO_X)$ via the trace map,
inducing an isomorphism
$$
\R\Hom_{\tilde p}(\E_1,\E_2)\,[3]\cong\Ext^3_{\tilde p\,}(\E,\E)
\Rt{\tr}R^3\tilde p_*\OO_{X\x\Delta}\cong\OO_\Delta.
$$
Again by Grothendieck-Riemann-Roch we have
$$
ch(\iota_*\OO_\Delta)\sqrt{\Td(\M_1\x\M_2)}=\iota_*(ch(\OO_\Delta)),
$$
so plugging everything into (\ref{pq}) (and noting that the shift by
[3] introduces a minus sign in $ch$) gives
$$
ff'(c)=q_{2\,*}(ch(\OO_\Delta)\,q_1^*(c))=c. \vspace{-8mm}
$$
\hspace{12cm}
\end{Proof}

\begin{Corollary}
Let $\M$ be a relative moduli space of slope stable sheaves on $X$
with data (\ref{data2}), and $d=2,\ r\ge2$. Then the groups $I^k=
\bigoplus_{i-j=k}H^{i,j}$ are isomorphic on $X$ and $\M$; thus their
Hodge numbers $h^{i,j}$ are the same.
\end{Corollary}

\begin{Proof}
Since the Fourier-Mukai machinery of \cite{BrM} now gives a quick
proof of this result we will only give a sketch of our previous argument.

Note that all we want to do is reverse the roles of $X$ and
$\M$, using the universal sheaf on $X\xp\M$, restricted to fibres
$\{x_t\}\times\M_t$, to exhibit $X$ as a fibrewise moduli space
of sheaves on $\M$ and deduce (\ref{muk}) that $f'\circ f=\id$.

Pick a smooth fibre $X_t$ and consider the universal bundle on $X_t
\times\M_t$. Firstly, by results of Mukai, (see e.g. \cite{Br2} Theorem
1.1), the corresponding Fourier-Mukai transform is an equivalence.
By using the corresponding transform for Hermitian-Yang-Mills
connections in differential geometry (see for example \cite{BBH}), one
can show that the transform maps slope stable
sheaves to slope stable sheaves, making $X_t$ a moduli space of
slope stable sheaves on $\M_t$ of ``Fourier-Mukai type'' (that is the
Mukai vector is primitive and of square zero; equivalently,
gcd\,$(r,c_1\,.\,\omega,\frac12c_1^2-c_2)=1$ for an appropriate choice
of polarisation).

Thus even on singular fibres semistable sheaves are stable and so
simple, so long as we can show that $X$ really does parametrise
sheaves on $\M$ over these fibres; i.e. we need to show that the
universal sheaf on $X\times\M$ is \emph{flat over} $X$.

But this follows by mimicking the argument of (\cite{Br1} Lemma 5.1)
one dimension up on $K3$ fibrations instead of elliptic fibrations.
(All that \cite{Br1} uses is the fact that the sheaves parametrised by
$\M$ are flat over $\M$ and have a locally free resolution of length 2
on $X$. In our case flatness over $\M$ is again immediate, and since
the sheaves are reflexive on the fibres by Theorem \ref{univ} they
have depth 2 and so homological dimension 1 by the Auslander-Buchsbaum
theorem (\cite{HL} 1.1).)

Using (Gieseker) stability we get the condition (\ref{van}) so that
the proof of Theorem \ref{muk} goes through with $X$ and $\M$
exchanged.

Since $ch$ and Td are of pure $(p,p)$ Hodge types, the
$I^k$s are preserved.  Calculating their dimensions in terms of
$h^{0,1},\ h^{1,1}$ and $h^{1,2}$, by symmetries of the Hodge diamond,
we see that these three numbers are preserved.
\end{Proof}

Theorem \ref{muk} is of course the fibrewise Fourier-Mukai transform
of \cite{BrM} at the level of K-theory. Theorem \ref{2d} can also be
interpreted in terms of this transform; the transform of
$\OO_{[\sigma]}$, where $[\sigma]$ is the image of a minimal degree
section of $\M\to\Pee^1$, is just $\tilde\sigma^*T$. Thus one might
try to use \cite{BrM} to study moduli of sheaves on $X$; the
difficulties are then displaced to understanding when the transform
gives stable sheaves, and not more exotic elements of the derived
category. Nonetheless, one might try to study moduli by
transforming the structure sheaves of connected curves in
a class that intersects the fibre class once. The singular nodal
sections with components lying in fibres will give rise to fibres on
which instabilities and singularities of the sheaves lie.

\subsubsection*{Examples}
We have already mentioned how the two dimensional results apply to
$K3\x T^2$. Now we again study Gross' examples, this time considering
two dimensional  moduli spaces, doing a family version of an example
in (\cite{HL} 5.3.7) and extending the proof of stability to \emph{all}
quartics in $\Pee^3$.

\begin{Prop}
Fix any irreducible, reduced quartic surface $S\subset\Pee^3$ which
contains no family of curves of arithmetic genus 0, and a
closed point $x\in S$. Then the rank two sheaf $E=E_x$ defined by the
sequence
$$
0\to E\to H^0(\I_x(1))\to\I_x(1)\to0,
$$
is slope stable with respect to the restriction of the Fubini-Study metric.
$E_x$ is a bundle for $x$ a smooth point of $S$, and the above gives an
isomorphism between the appropriate moduli space and $S$.
\end{Prop}

\begin{Proof} The dual $E^*$ of such a sheaf is generated by its
sections, and a destabilising rank one (torsion free, without
loss of generality) quotient sheaf $\Ll$ is therefore also generated
by its sections and of degree $\le2$. Therefore the proof of (\ref{S})
applies to obtain a contradiction.
\end{Proof}

So in this case, for the two Gross 3-folds $X_i$, the relative moduli
spaces $\M_i$ are in fact isomorphic to $X_i$ for the generic $X_i$
with fibres satisfying the above conditions. (We need only calculate
the holomorphic Casson invariant of the generic $X_i$, by deformation
invariance.) It is easy
to see that the polarisation on $\M_i$ is that on $X$ plus $N\omega_1$,
since this gives the right fibre polarisation on $\M_t$ for the generic
fibre $X_t$ of Picard number one. Unfortunately this example only
satisfies the conditions of Definition \ref{data},
not Definition \ref{data2}, so Theorem \ref{2d} does not really apply.
But the extra conditions of Definition \ref{data2} were only used
to show that Jun Li's line bundle (\ref{above}) was ample on the singular
fibres of $\M\to\Pee^1$, and that the smooth points of the fibres
corresponded to locally free sheaves, both of which, as mentioned
earlier, should be true in general anyway. I cannot prove it in the
general case, but it clearly holds in this example, so we may use
Theorem \ref{2d}.

Thus the holomorphic Casson invariants
are given by counting the number of lowest degree sections of
$$
\big(X_i\to\Pee^1\big)\subset\big(\Pee(E_i)\to\Pee^1\big)
$$
with respect to the form $t=c_1(\OO_{\Pee(E_i)}(1))$, which is ample on
the fibres. For $X_2$ the subbundle $\OO(1)\subset E_2=\OO(-1)\oplus
\OO\oplus\OO\oplus \OO(1)$ gives a section of $\Pee(E_2)\to\Pee^1$ which
lies in $X_2$, is of lowest degree $(-1)$ with respect to $t$, and is
unique. Thus we get a Casson invariant of one, counting a locally free sheaf.

For the moduli space that corresponds to this one
under the diffeomorphism between
the $X_i$, however, we must consider sections of $X_1\to\Pee^1$ of degree
$-1$ with respect to the K\"ahler form of $\Pee^3$ on $X\subset\Pee^1\x
\Pee^3$. There clearly are none.

(The Casson invariant of $X_1$ that corresponds to lowest degree
sections is for degree 0 sections,
i.e. those of the form $\Pee^1\x\{x\}\subset X_1$. Write the equation
defining the (2,4) divisor $X_1\cong\M_1$ as $x^2f+xyg+y^2h=0$, where
$f,\,g,\,h$ are quartic polynomials on $\Pee^3$. Then $\Pee^1\x\{x\}$ lies
in $S$ if and only if $x$ lies in the intersection of the three quartic
surfaces in $\Pee^3$. Thus we may deform to a case where there are
precisely 64 isolated sections in the moduli space, and therefore the
Casson invariant is 64, and again all the sheaves are locally free.)

Thus we have finally recovered the result that the $X_i$ are not
deformation equivalent, as polarised varieties. Although the result
is both weaker than the results using Gromov-Witten invariants
\cite{Ru} and essentially a more complicated rerun of the same proof (going
as it does via the same Gromov-Witten invariants), it is
nonetheless encouraging that there are examples of the invariant which are
calculable yet contain highly non-trivial information.

\begin{flushleft}
Harvard University, 1, Oxford Street, Cambridge MA 02138. USA. \newline
\small Email: \tt thomas@maths.ox.ac.uk
\end{flushleft}


\begin{thebibliography}{[FMW]}

\bibitem[AOS]{AOS}
Acharya, B., O'Loughlin, M. and Spence, B. (1997).
\emph{Higher Dimensional Analogues of Donaldson-Witten Theory.}
Preprint hep-th/9705138.

\bibitem[Ar]{Ar}
Artamkin, I.\,V. (1989).
\emph{On deformation of sheaves.} Math. USSR Izv. \textbf{32}
663--668.

\bibitem[BPS]{BPS}
Banica, C., Putinar, M. and Schumacher, G. (1980).
\emph{Variation der globalen Ext in deformation kompakter komplexer
r\"aume.} Math. Ann. \textbf{250}, 135--155.

\bibitem[BBH]{BBH}
Bartocci, C., Bruzzo, U., Hern\'andez Ruip\'erez, D. (1997).
\emph{A Fourier-Mukai transform for stable
bundles on $K3$ surfaces.} J. Reine Angew. Math. \textbf{486}, 1--16.

\bibitem[BKS]{BKS}
Baulieu, L. Kanno, H. and Singer, I.\,M. (1997).
\emph{Special Quantum Field Theories In Eight And Other Dimensions.}
Preprint hep-th/9704167.

\bibitem[Be]{Be}
Behrend, K. (1997).
\emph{Gromov-Witten Invariants in Algebraic Geometry.}
Invent. Math. \textbf{127}, 601--617.

\bibitem[BF]{BF}
Behrend, K. and Fantechi, B. (1997).
\emph{The intrinsic normal cone.} Invent. Math. \textbf{128}, 45--88.

\bibitem[Br1]{Br1}
Bridgeland, T. (1998).
\emph{Fourier-Mukai transforms for elliptic surfaces.}
J. Reine Angew. Math. \textbf{498}, 115--133.

\bibitem[Br2]{Br2}
Bridgeland, T. (1999).
\emph{Equivalences of triangulated categories and Fourier-Mukai transforms.}
Bull. Lond. Math. Soc. \textbf{30}, 25--34.

\bibitem[BrM]{BrM}
Bridgeland, T., and Maciocia, A. (1999).
\emph{Fourier-Mukai transforms for K3 fibrations.} Preprint
math.AG/9908022.

\bibitem[DK]{DK}
Donaldson, S.\,K. and Kronheimer, P.\,B. (1990).
\emph{The geometry of four-manifolds.} Oxford University Press.

\bibitem[DT]{DT}
Donaldson, S.\,K. and Thomas, R.\,P. (1998).
\emph{Gauge theory in higher dimensions.} In: The Geometric Universe:
Science, Geometry and the work of Roger Penrose, S.\,A. Huggett et al
(eds), Oxford University Press.

\bibitem[FK]{FK}
Frenkel, I. and Khesin, B. (1996).
\emph{Four Dimensional Realization of Two Dimensional Current Groups.}
Comm. Math. Phys. \textbf{178}, 541--562.

\bibitem[FKT]{FKT} 
Frenkel, I., Khesin, B. and Todorov, A. (1997).
\emph{Complex counterpart of Chern-Simons-Witten theory and
holomorphic linking.}
Preprint.

\bibitem[Fr]{Fr}
Friedman, R. (1990s).
\emph{Deformations of sheaves on surfaces.}
Unpublished manuscript.

\bibitem[FMW]{FMW}
Friedman, R., Morgan, J. and Witten, E. (1997).
\emph{Vector bundles over Elliptic Fibrations.}
Preprint alg-geom/9709029.

\bibitem[Fu]{Fu}
Fulton, W. (1984).
\emph{Intersection Theory.}
Springer-Verlag.

\bibitem[GH]{GH}
Griffiths, P. and Harris, J. (1978).
\emph{Principles of algebraic geometry.} Wiley, New York.

\bibitem[GH2]{GH2}
Griffiths, P. and Harris, J. (1985).
\emph{On the Noether-Lefschetz theorem and some remarks on
codimension-2 cycles.} Math. Ann. \textbf{271}, 31--51. 

\bibitem[Gr]{Gr}
Gross, M. (1997).
\emph{The deformation space of Calabi-Yau $n$-folds can be obstructed.}
Mirror symmetry II, AMS/IP Stud. Adv. Math., 1, Amer.
Math. Soc., Providence, RI, 401--411.

\bibitem[Ha1]{Ha1}
Hartshorne, R. (1977).
\emph{Algebraic Geometry.} Graduate Texts in Mathematics 52, Springer-Verlag.

\bibitem[Ha2]{Ha2}
Hartshorne, R. (1980).
\emph{Stable reflexive sheaves.} Math. Ann. \textbf{254}, 121--176.

\bibitem[Ha3]{Ha3}
Hartshorne, R. (1966).
\emph{Residues and duality.} Lecture Notes in Mathematics 20,
Springer-Verlag.

\bibitem[HL]{HL}
Huybrechts, D. and Lehn, M. (1997).
\emph{Geometry of moduli spaces of shaves.}
Aspects in Mathematics Vol. E31, Vieweg.

\bibitem[Is]{Is}
Ishii, A. (1996).
\emph{Semi-universal family of reflexive modules over a rational
double point of type $A$.} Moduli of vector bundles (Kyoto, 1994),
65--77, Lecture Notes in Pure and Appl. Math., 179, Dekker.

\bibitem[Kh]{Kh}
Khesin, B. (1997).
\emph{Informal Complexification and Poisson Structures on Moduli
spaces.} Amer. Math. Soc. Transl. (2) \textbf{180}, 147--155.

\bibitem[KR]{KR}
Khesin, B. and Rosly, A. (1998).
\emph{Canadian homology and holomorphic linking of complex manifolds.}
In preparation.

\bibitem[La]{La}
Langer, A. (2000).
\emph{Chern classes of reflexive sheaves on normal surfaces.}
Mathematische Zeitschrift, to appear.

\bibitem[L]{L}
Lehn, M. (1998).
\emph{On the cotangent sheaf of Quot-schemes.} Internat. Jour. Math.
\textbf{9}, 513--522.

\bibitem[Li]{Li}
Li, J. (1993).
\emph{Algebraic geometric interpretation of Donaldson's polynomial
invariants of algebraic surfaces.} Jour. Diff. Geom. \textbf{37} 417--466.

\bibitem[LT]{LT}
Li, J. and Tian, G. (1998).
\emph{Virtual moduli cycles and Gromov-Witten invariants of algebraic
varieties.} Jour. Amer. Math. Soc. \textbf{11} 119--174.

\bibitem[Mu1]{Mu1}
Mukai, S. (1984).
\emph{Symplectic structure of the moduli space of stable 
sheaves on an abelian or K3 surface.} Invent. Math.
\textbf{77}, 101--116.

\bibitem[Mu2]{Mu2}
Mukai, S. (1987).
\emph{On the moduli space of bundles on $K3$ surfaces, I.} Vector Bundles
on Algebraic Varieties, Oxford University Press, 341--413.

\bibitem[OSS]{OSS}
Okonek, C., Schneider, M. and Spindler, H. (1980).
\emph{Vector bundles on complex projective spaces.}
Progress in Math. 3, Birkhauser.

\bibitem[Pi]{Pi}
Pidstrigatch, V. (1991).
\emph{Deformations of instanton films.}
Izv. Akad. Nauk SSSR, Ser. Mat. \textbf{55}, 318-338. Translation in:
Math. of the USSR; Izvestya, \textbf{38}, 313--331, 1992.

\bibitem[Ru]{Ru}
Ruan, Y. (1996).
\emph{Topological sigma model and Donaldson-type invariants in Gromov
theory.} Duke Math. Jour. \textbf{83}, 461--500.

\bibitem[Si]{Si}
Siebert, B. (1997).
\emph{Global normal cones, virtual fundamental classes and Fulton's
canonical classes.} Preprint.

\bibitem[SYZ]{SYZ}
Strominger A., Yau S.-T. and Zaslow E. (1996).
\emph{Mirror Symmetry is T-Duality.} Nucl. Phys. \textbf{B479}, 243--259.

\bibitem[Ta]{Ta}
Taubes, C.\,H. (1990).
\emph{Casson's invariant and gauge theory.} Jour. Diff. Geom.
\textbf{31}, 547--599.

\bibitem[T1]{T1}
Thomas, R.\,P. (1997).
\emph{Gauge Theory on Calabi-Yau manifolds.} D.\,Phil. Thesis,
University of Oxford.

\bibitem[T2]{T2}
Thomas, R.\,P. (1999).
\emph{An obstructed bundle on a Calabi-Yau 3-fold}, Adv. Theor.
and Math. Phys. \textbf{3}.

\bibitem[Ti]{Ti}
Tian, G. (2000).
\emph{Gauge Theory and Calibrated Geometry, I.} Ann. Math. \textbf{151},
193--268.

\bibitem[Ty1]{Ty1}
Tyurin, A.\,N. (1990).
\emph{The moduli spaces of vector bundles on threefolds, surfaces
and curves.} Math. Inst. Univ. Gottingen preprint.

\bibitem[Ty2]{Ty2}
Tyurin, A.\,N. (1997).
\emph{Non-abelian analogues of Abel's theorem.} Preprint.

\bibitem[VW]{VW}
Vafa, C. and Witten, E. (1994).
\emph{A strong coupling test of S-duality.} Nucl. Phys.
\textbf{B 431}, 3--77.

\end{thebibliography}
\end{document}